\documentclass[12pt]{article}
\usepackage[utf8]{inputenc}
\usepackage[T1]{fontenc}
\usepackage[english]{babel}
\usepackage[top=2cm, bottom=2cm, left=2cm, right=2cm]{geometry}
\usepackage{layout}
\usepackage{setspace}
\usepackage{charter}
\usepackage{url}
\usepackage{graphicx}
\usepackage{amsmath}
\usepackage{amssymb}
\usepackage{mathrsfs}
\usepackage{amsthm}
\usepackage{hyperref}
\usepackage{titlesec}
\titleformat{\paragraph}
{\normalfont\normalsize\bfseries}{\theparagraph}{1em}{}
\titlespacing*{\paragraph}
{0pt}{3.25ex plus 1ex minus .2ex}{1.5ex plus .2ex}
\usepackage{enumitem}
\setlist{itemsep = 3 pt}

\makeatletter
\newcommand{\specialcell}[1]{\ifmeasuring@#1\else\omit$\displaystyle#1$\ignorespaces\fi}
\makeatother

\theoremstyle{definition}
\newtheorem{lem}{Lemma}[section]
\newtheorem{prop}{Proposition}[section]
\newtheorem{defin}{Definition}[section]
\newtheorem{nota}{Notation}[section]
\newtheorem{thm}{Theorem}[section]

\newtheorem{rmk}{Remark}[section]

\let\oldbibliography\thebibliography
\renewcommand{\thebibliography}[1]{\oldbibliography{#1}
\setlength{\itemsep}{-3pt}} %Reducing spacing in the bibliography.

\title{Four families of maximal real algebraic hypersurfaces in $\mathbb{RP}^4$}
\author{Aloïs Demory}
\date{2023}

\begin{document}

%TITRE

%\vspace*{\stretch{1}}
\begin{center}
	\begin{minipage}{15cm}
		\maketitle
	\end{minipage}
\end{center}
%\vspace*{\stretch{1}}

%\newpage
\renewcommand{\contentsname}{Sommaire}
\renewcommand\labelitemi{---}

%CORPS
\section{Introduction}
\label{definitions}

\begin{sloppypar}
The first part of D. Hilbert's sixteenth problem, interpreted broadly, asks to study the topology of real algebraic varieties.
\begin{defin}
A real algebraic hypersurface $A$ of degree $m$ in $\mathbb{RP}^n$ is given by the datum of a homogeneous degree $m$ polynomial $P_A$ in $n+1$ variables with real coefficients up to multiplication by a non-zero real constant. The hypersurface $A$ has the real part \mbox{$\mathbb{R}A := \{ [x_0:...:x_n] \in \mathbb{RP}^n | P_A([x_0:...:x_n]) =0 \}$} and the complex part \mbox{$\mathbb{C}A := \{ [z_0:...:z_n] \in \mathbb{CP}^n | P_A([z_0:...:z_n]) =0 \}$}.\\
The complex conjugation on $\mathbb{CP}^n$ induces an involution on $\mathbb{C}A$, whose fixed point set is $\mathbb{R}A$. If there does not exist any point $z \in \mathbb{CP}^n$ such that $P_A(z)=\frac{\partial P_A}{\partial z_0}(z)=...=\frac{\partial P_A}{\partial z_n}(z)=0$, then the hypersurface $A$ is said to be non-singular. In this case, $\mathbb{C}A$ and $\mathbb{R}A$ are smooth manifolds.
\end{defin}
The main result regarding the topology of non-singular real algebraic curves in $\mathbb{RP}^2$, known at the time of the formulation of Hilbert's problems, was A. Harnack's bound (see \cite{harnack}), which states that the number of connected components of the real part of a non-singular real algebraic curve of given degree $m$ in $\mathbb{RP}^2$ is at most $g+1$, where $g = \frac{(m-1)(m-2)}{2}$ is the genus of the complex part of the curve. Harnack also proved that this bound is sharp for every degree by constructing real curves with the maximal number of connected components. This way, he introduced the two main directions of the subject: finding restrictions on the topology of real algebraic hypersurfaces and constructing hypersurfaces that are extremal with respect to these restrictions. 
%Since then, a lot of progress has been achieved in both directions. 
The modern generalization of Harnack's bound is the following result. The field with two elements $\mathbb{Z}/2\mathbb{Z}$ is denoted by $\mathbb{Z}_2$.
\begin{thm}
\label{smith}
\textbf{(Smith-Thom inequality and congruence, see e.g. \cite{floyd}, \cite{bredon})}
Let $A$ be a real algebraic hypersurface of degree $m$ in $\mathbb{RP}^n$. We have the following inequality relating the total $\mathbb{Z}_2$-homology of the real and complex parts of $A$:
\begin{center}$\sum_{i=0}^{n-1} dim_{\mathbb{Z}_2}H_i(\mathbb{R}A;\mathbb{Z}_2) \leq \sum_{i=0}^{2n-2} dim_{\mathbb{Z}_2}H_i(\mathbb{C}A;\mathbb{Z}_2)$.\end{center}
Furthermore, $\sum_{i=0}^{n-1} dim_{\mathbb{Z}_2}H_i(\mathbb{R}A;\mathbb{Z}_2) \equiv \sum_{i=0}^{2n-2} dim_{\mathbb{Z}_2}H_i(\mathbb{C}A;\mathbb{Z}_2)$ (mod $2$).
\end{thm}
\begin{defin}
The hypersurface $A$ is said to be an \textit{M-hypersurface} (or to be \textit{maximal}) if $\sum_{i=0}^{n-1} dim_{\mathbb{Z}_2}H_i(\mathbb{R}A;\mathbb{Z}_2) = \sum_{i=0}^{2n-2} dim_{\mathbb{Z}_2}H_i(\mathbb{C}A;\mathbb{Z}_2)$.    
\end{defin}
When $A$ is non-singular, the number $\sum_{i=0}^{2n-2} dim_{\mathbb{Z}_2}H_i(\mathbb{C}A;\mathbb{Z}_2)$ only depends on $n$ and $m$: it is equal to $\sum_{j=0}^{n-1} (-1)^j m^{j+1} {n+1 \choose j+2}$ if $n$ is odd and $2n - \sum_{j=0}^{n-1} (-1)^j m^{j+1} {n+1 \choose j+2}$ if $n$ is even (see \cite{danilkhovan}). The sharpness of an asymptotical version of the Smith-Thom inequality was proved for every integers $n$, $m>0$ by I. Itenberg and O. Viro (see \cite{itenbergviro2007}) using Viro's combinatorial patchworking method (see \textit{e.g.} \cite{viropatchwork}). The authors also announced that the sharpness of the actual inequality could be proved by similar means (see \cite{fakeitenbergviro}). Maximal real algebraic hypersurfaces turn out to have interesting properties, one of the most remarkable ones being the following congruence.
\begin{thm}
\textbf{(Rokhlin's congruence, \cite{rokhlin})} Let $n$ be an odd positive integer and let $A$ be a non-singular maximal real algebraic hypersurface in $\mathbb{RP}^n$. Then $\chi(\mathbb{R}A) \equiv \sigma(\mathbb{C}A)$ mod $16$, where $\chi(\mathbb{R}A)$ is the Euler characteristic of $\mathbb{R}A$ and $\sigma(\mathbb{C}A)$ is the signature of $\mathbb{C}A$.
\end{thm}
As noticed by E. Brugallé in \cite{brugalle}, a large class of maximal real algebraic hypersurfaces that have been constructed up to this day respect the stronger equality $\chi(\mathbb{R}A) = \sigma(\mathbb{C}A)$. There are few examples of M-hypersurfaces that do not respect this equality, namely some surfaces in $\mathbb{RP}^3$ for every degree $m\geq 4$ (see e.g. \cite{kharlamov}) and some cubic fourfolds in $\mathbb{RP}^5$ (\cite{finashinkharlamov}).
\end{sloppypar}

\begin{sloppypar}
Let $A$ be a non-singular real algebraic hypersurface in $\mathbb{RP}^n$ defined by a homogeneous polynomial $P_A$. When the dimension of $A$ is odd, the Euler characteristic of $\mathbb{R}A$ is always zero. In this case, we look at the following related even-dimensional objects, assuming that the degree of $P_A$ is even.
\begin{nota}
We denote by $A_-$ (respectively, by $A_+$) the set of solutions of the inequality $P_A \leq 0$ (respectively, $P_A \geq 0$) in $\mathbb{RP}^n$.
\end{nota}
\begin{nota}
Let $Y$ be a double covering of $\mathbb{CP}^n$ branched along $\mathbb{C}A$ (such a covering exists and is unique up to isomorphism of ramified coverings, see \textit{e.g.} \cite{ramifcover}). The complex conjugation on $\mathbb{CP}^n$ is covered by two antiholomorphic involutions on $Y$, with respective fixed point sets $\mathbb{R}Y_-$ and $\mathbb{R}Y_+$ lying over $A_-$ and $A_+$. When $A$ is maximal, we have either \mbox{$\sum_{i=0}^n dim_{\mathbb{Z}_2}H_i(\mathbb{R}Y_-;\mathbb{Z}_2) = \sum_{i=0}^{2n} dim_{\mathbb{Z}_2}H_i(Y;\mathbb{Z}_2)$} or \mbox{$\sum_{i=0}^n dim_{\mathbb{Z}_2}H_i(\mathbb{R}Y_+;\mathbb{Z}_2) = \sum_{i=0}^{2n} dim_{\mathbb{Z}_2}H_i(Y;\mathbb{Z}_2)$}.
%; the maximal fixed point set lying over a region of $\mathbb{RP}^n$ that contains a representative of the non-zero homology class of $H_{n-1}(\mathbb{RP}^n)$. 
In the rest of the paper, unless otherwise stated, we suppose that when $A$ is maximal, we have \mbox{$\sum_{i=0}^n dim_{\mathbb{Z}_2}H_i(\mathbb{R}Y_-;\mathbb{Z}_2) = \sum_{i=0}^{2n} dim_{\mathbb{Z}_2}H_i(Y;\mathbb{Z}_2)$}. 
\end{nota}
Hence, when $n-1$ is odd and the degree of $P_A$ is even, we can compare the Euler characteristic of $\mathbb{R}Y_-$ (or $\mathbb{R}Y_+$, as $\chi(\mathbb{R}Y_-) + \chi(\mathbb{R}Y_+) = 2 \chi(\mathbb{RP}^4) =2$) with the signature of $Y$. This signature depends only on $m$ and $n$. In the case that is considered in this paper, \textit{i.e.} when $Y$ is a double covering of $\mathbb{CP}^4$ branched along a non-singular real algebraic hypersurface of even degree $m$ in $\mathbb{RP}^4$, we have $\sigma(Y) = \frac{1}{24}(5m^4-20m^2+48)$ (see \cite{danilkhovan}).
\begin{defin}
A non-singular maximal real algebraic hypersurface $A$ in $\mathbb{RP}^n$ is said to be \mbox{\textit{of type} $\chi = \sigma$} if either $n$ is odd and $\chi(\mathbb{R}A) = \sigma(\mathbb{C}A)$ or $n$ is even, the degree of $P_A$ is even and $\chi(\mathbb{R}Y_-) = \sigma(\mathbb{C}Y)$.
%\begin{itemize}
%    \item either $n$ is odd and $\chi(\mathbb{R}A) = \sigma(\mathbb{C}A)$;
%    \item or $n$ is even, the degree of $P_A$ is even and $\chi(\mathbb{R}Y_-) = \sigma(\mathbb{C}Y)$.
%\end{itemize}
\end{defin}

\begin{defin}
A family $(A_m)_{m \in \mathbb{N}^*}$ of non-singular maximal real algebraic hypersurfaces in $\mathbb{RP}^n$ is said to be \textit{asymptotically of type} $\chi = \sigma$ if either $n$ is odd and $\chi(\mathbb{R}A_m) \sim_{+\infty} \sigma(\mathbb{C}A_m)$ or $n$ is even, $A_m$ is defined only for even positive integers $m$ and $\chi(\mathbb{R}Y^m_-) \sim_{+\infty} \sigma(Y^m)$ where $Y^m$ is a double covering of $\mathbb{CP}^n$ branched along $\mathbb{C}A_m$.
%\begin{itemize}
%    \item either $n$ is odd and $\chi(\mathbb{R}A_m) \sim_{+\infty} \sigma(\mathbb{C}A_m)$;
%    \item or $n$ is even, $A_m$ is defined only for even positive integers $m$ and $\chi(\mathbb{R}Y^m_+) \sim_{+\infty} \sigma(Y^m)$ where $Y^m$ is a double covering of $\mathbb{CP}^n$ ramified over $\mathbb{C}A_m$.
%\end{itemize}
\end{defin}
In \cite{Hilbert1933}, Hilbert presented one of the first families of non-singular maximal real algebraic curves in $\mathbb{RP}^2$ that are not asymptotically of type $\chi = \sigma$. 
%Considering double coverings of $\mathbb{CP}^2$ branched along the complex parts of the curves of this family, one obtains a family of non-singular maximal real algebraic surfaces in a weighted real projective $3$-space. This family of surfaces is not asymptotically of type $\chi = \sigma$. 
Numerous other such families of curves have been described since then. Some families of non-singular maximal real algebraic surfaces in $\mathbb{RP}^3$ that are not asymptotically of type $\chi = \sigma$ are also known (see e.g. \cite{Vir79}). However, no example of a family of hypersurfaces in $\mathbb{RP}^4$ or in higher dimension that is not asymptotically of type $\chi = \sigma$ has been exhibited so far.
\end{sloppypar}
\begin{sloppypar}
In this paper, we describe four families of non-singular maximal real algebraic hypersurfaces in $\mathbb{RP}^4$. They are constructed using Viro's combinatorial patchworking method (\cite{viropatchwork}). The first two families, called respectively the odd and the even Itenberg-Viro family, are constituted of maximal hypersurfaces of type $\chi = \sigma$. The hypersurfaces in these families are particular cases of the projective hypersurfaces to be described in \cite{fakeitenbergviro}. We obtained the other two families by modifying the construction of the odd and even Itenberg-Viro hypersurfaces. Those two new families of hypersurfaces in $\mathbb{RP}^4$, respectively called the Small Devation and the Asymptotical Deviation family, allow us to prove the following results.
\begin{thm}
\label{theorem1}
For any integer $k \geq 5$, there exists a non-singular maximal real algebraic hypersurface $X_{2k}$ of degree $2k$ in $\mathbb{RP}^4$ that is not of type $\chi=\sigma$ and whose number of connected components is bigger than $h^{n,0}(\mathbb{C}X_{2k})+1$. It verifies the equality \mbox{$\chi(\mathbb{R}Y^{2k}_-) = \sigma(Y^{2k}) -  8(\frac{(2k)^3}{24}-\frac{5(2k)^2}{8}+\frac{17(2k)}{6}-4+\frac{1}{4}\bar{2k}$},  where $Y_{2k}$ is a double covering of $\mathbb{CP}^4$ branched along $\mathbb{C}X_{2k}$ and $\bar{2k}$ is the remainder of the euclidean division of $2k$ by $4$.
\end{thm}
\begin{thm}
\label{theorem2}
There exists a family $(X_{2k})_{k \geq 4}$ of non-singular maximal real algebraic hypersurfaces in $\mathbb{RP}^4$ which is not asymptotically of type $\chi = \sigma$. It satisfies the relation \mbox{$\chi(\mathbb{R}Y_{2k}^-) \sim_{+\infty} \sigma(\mathbb{C}Y_{2k}) - \frac{(2k)^4}{4}$, where $Y_{2k}$} is a double covering of $\mathbb{CP}^4$ branched along $\mathbb{C}X_{2k}$.
\end{thm}
We begin the paper with a description of the combinatorial patchworking and the statement of some important related results in Section \ref{patchwork}. As a preliminary to the study of Itenberg-Viro $3$-manifolds, we introduce in Section \ref{prelim} a construction of maximal surfaces due to Itenberg and Viro and prove some of the properties of the surfaces resulting from this construction. We then describe in Section \ref{construction3} the data necessary to construct the odd and the even Itenberg-Viro hypersurfaces using the combinatorial patchworking method. Sections \ref{maximality3manif} and \ref{prooftype} are devoted to a proof of some properties of the even and odd Itenberg-Viro hypersurfaces. The tools presented in Sections \ref{prelim}, \ref{maximality3manif} and \ref{prooftype} are ubiquitous in the paper. In Section \ref{smalldeviation}, we explain how to modify the odd Itenberg-Viro datum to obtain the Small Deviation family of hypersurfaces and we prove \mbox{Theorem \ref{theorem1}}. Finally, in Section \ref{asymptoticalmodification} we present a modification of the even Itenberg-Viro datum and describe the Asymptotical Deviation family of hypersurfaces, which allows us to prove \mbox{Theorem \ref{theorem2}}.
\end{sloppypar}

\section{Combinatorial patchworking of non-singular real algebraic hypersurfaces in $\mathbb{RP}^n$}
\label{patchwork}
\begin{sloppypar}
Let $m$ be a positive integer (it would be the degree of the to-be-constructed hypersurface), and let $\Delta_m^n$ be the $n$-dimensional simplex in $\mathbb{R}^n$ with vertices $(0,...,0)$, $(m,0,...,0)$, $(0,m,0,..,0)$, ... , $(0,...,0,m)$ . Consider a rectilinear triangulation $\tau$ of $\Delta_m^n$ with vertices having integer coordinates and a distribution of signs $\mu : V \rightarrow \{+,-\} $, where $V$ is the set of vertices of $\tau$.
Let $(\Delta_m^{n})^*$ be the union of the symmetric copies of $\Delta_m^n$ under reflections with respect to the coordinate hyperplanes and compositions of such reflections. Extend $\tau$ to a symmetric triangulation $\tau^*$ of $(\Delta_m^{n})^*$. Extend the sign distribution $\mu$ to a sign distribution $\mu^*$ at the vertices of $\tau^*$ using the following rule : if two vertices are the images of each other under a reflection with respect to a certain coordinate hyperplane, then their signs are the same if their distance to the hyperplane is even, and their signs are opposite otherwise.
If an $n$-simplex of $\tau^*$ has vertices with different signs, consider the convex hull of the middle points of its edges that have vertices of opposite sign. It is a piece of hyperplane. Denote by $\Gamma$ the union of all such hyperplane pieces. It is a piecewise-linear hypersurface contained in $(\Delta_m^{n})^*$. When $m$ is even, we denote by $\Gamma_-$ (respectively, $\Gamma_+$) the union of the regions of $(\Delta_m^{n})^* \setminus \Gamma$ which contain only vertices bearing the sign $-$ (respectively, the sign $+$).
Now, identify each pair of antipodal points lying on the boundary of $(\Delta_m^{n})^*$ and denote by $\widetilde{\Gamma}$ (respectively, $\widetilde{\Gamma}_-$, $\widetilde{\Gamma}_+$) the image of $\Gamma$ (respectively, $\Gamma_-$, $\Gamma_+$) under this identification. Notice that the quotient space $(\widetilde{\Delta_m^{n})^*}$ is $PL$-homeomorphic to the projective space $\mathbb{RP}^n$.
The triangulation $\tau$ is said to be \textit{convex} if there exists a piecewise-linear convex function $\nu:\Delta_m^n \rightarrow \mathbb{R}$ whose domains of linearity coincide with the $n$-simplices of $\tau$. The following result is due to Viro.
\begin{thm}
\textbf{(\cite{vir83}, see also \cite{viropatchwork})} If the convexity of $\tau$ is certified by a function $\nu$, then there exists a non-singular real algebraic hypersurface $X$ of degree $m$ in $\mathbb{RP}^n$ and a homeomorphism $\mathbb{RP}^n \rightarrow (\widetilde{\Delta}_m^{n})^*$ mapping $\mathbb{R}X$ onto $\widetilde{\Gamma}$. Such a hypersurface $X$ can be taken to be given by a homogeneous polynomial of the form $\sum_{(i_1,...,i_n) \in V} \mu(i_1,...,i_n)x_0^{m-i_1-...-i_n}x_1^{i_1}...x_n^{i_n}t^{\nu (i_1,...,i_n)}$ for some positive and sufficiently small real number $t$.
%Moreover, when $m$ is even, there also exists an homeomophism between the pairs $(\mathbb{RP}^n, X_-)$ (respectively $(\mathbb{RP}^n, X_+)$) and $((\widetilde{\Delta}_m^{n})^*, \widetilde{\Gamma}_-)$ (respectively $((\widetilde{\Delta}_m^{n})^*, \widetilde{\Gamma}_+)$) for some positive and sufficiently small $t$.
\end{thm}
\begin{defin}
A $k$-simplex in $\mathbb{R}^n$ with integer vertices $v_0$, ..., $v_k$ is said to be \textit{primitive} if the family of vectors $(v_i-v_0)_{1\leq i \leq k}$ is a $\mathbb{Z}$-basis of the lattice $\mathbb{Z}^n \bigcap E$, where $E$ is the vector subspace of $\mathbb{R}^n$ spanned by the vectors $v_1-v_0$, ..., $v_k-v_0$. When all the simplices of $\tau$ are primitive, the triangulation $\tau$ is said to be \textit{primitive}, and $X$ is called a \textit{primitive hypersurface}.    
\end{defin}
Primitive hypersurfaces of $\mathbb{RP}^n$ have very nice properties, such as the following ones.
\begin{thm}
\label{thmbertrand}
\textbf{(\cite{bertrand_2010})} If $X$ is a primitive hypersurface in $\mathbb{RP}^n$, then $\chi(\mathbb{R}X) = \sigma(\mathbb{C}X)$.
\end{thm}
When $n$ is even, it is only a matter of adopting the convention $\sigma(\mathbb{C}X)=0$. A proof for the case of surfaces can be found in \cite{tsurf}, while the general case is dealt with in \cite{bertrand_2010} (see also \cite{brugalle}). Furthermore, if $X$ is maximal, the following equalities hold.
\begin{thm}
\textbf{(\cite{renaudineau})} If $X$ is a maximal primitive hypersurface in $\mathbb{RP}^n$, then:
\begin{itemize}
    \item if $n-1$ is even, then $b_{\frac{n-1}{2}}(\mathbb{R}X) = h^{\frac{n-1}{2},\frac{n-1}{2}}(\mathbb{C}X)$.
    \item for every integer $i$ such that $0\leq i \leq n-1$ and $i \neq \frac{n-1}{2}$, one has $b_i(\mathbb{R}X) = h^{i,n-1-i}(\mathbb{C}X)+1$.
\end{itemize}
This implies that $\chi(\mathbb{R}X) = \sigma(\mathbb{C}X)$.
\end{thm}
Primitive odd-dimensional hypersurfaces are not necessarily of type $\chi=\sigma$. The first counterexamples are given by non-singular maximal real algebraic curves of degree 6 in $\mathbb{RP}^2$. To study further the odd-dimensional case, we recall one result proved by S. Orevkov.
\begin{lem}
\label{orevkov}
\textbf{(\cite[Lemma~4.2]{orevkov})} With the above notations, $\Gamma_-$ (respectively, $\Gamma_+$) deformation retracts to the union of the simplices of $\tau^*$ whose vertices all have the sign $-$ (respectively, the sign $+$).
\end{lem}
\end{sloppypar}

\section{Preliminaries: Itenberg-Viro surfaces}

\label{prelim}

\subsection{The $\mathbb{Z}_2$-linking number of spheres in $\mathbb{RP}^n$}

\label{sectionlinkingnumber}

\begin{defin}
    \label{deflink}
    Let $n$ be an integer such that $n\geq 3$. Let $k$ be an integer such that $0 \leq k \leq n-1$. Let $c$ and $a$ be smooth disjoint submanifolds of $\mathbb{RP}^n$ respectively diffeomorphic to $\mathbb{S}^k$ and $\mathbb{S}^{n-k-1}$. Suppose that $c$ bounds a smooth $(k+1)$-disk $d$ in $\mathbb{RP}^n$. The $\mathbb{Z}_2$\textit{-linking number of} $c$ \textit{and} $a$, denoted by $\mathrm{lk}(c,a)$, is the $\mathbb{Z}_2$-intersection product of $a$ and $d$. 
\end{defin}

\begin{sloppypar}
    For any integer $i \notin \{ 1, n \}$, the homotopy group $\pi_i(\mathbb{RP}^n)$ is trivial. Hence, when $k \neq n-2$, the $\mathbb{Z}_2$-linking number of $c$ and $a$ is independent from the choice of a smooth disk $d$. Indeed, let $d_1$ and $d_2$ be two such disks. As $a$ is contractible, $a$ is null-homologous in $\mathbb{RP}^n$ and the $\mathbb{Z}_2$-intersection product of the cycle $d_1 + d_2$ and $a$ is $0$. Hence, the $\mathbb{Z}_2$-intersection products of $a$ and $d_1$ and of $a$ and $d_2$ are the same. When $k = n-2$, the cycle $d_1 + d_2$ defines a map from $\mathbb{S}^{n-1}$ to $\mathbb{RP}^n$. Hence, it is contractible and null-homologous in $\mathbb{RP}^n$ and we reach the same conclusion, even when $a$ is a non-contractible loop in $\mathbb{RP}^n$.
\end{sloppypar}

\begin{lem}
    \label{Lemmelinkindep}
    Let $k$ be an integer such that $0 \leq k \leq n-1$. Let $H$ be a smooth hypersurface in $\mathbb{RP}^n$. Let $c$ and $a$ be disjoint smooth spheres in $\mathbb{RP}^n$ of respective dimensions $k$ and $n-k-1$. Suppose that $c$ bounds a smooth $(k+1)$-disk $d$ in $\mathbb{RP}^n$, that $c$ is contained in $H$ and that $a$ is contained in $\mathbb{RP}^n \setminus H$. If $\mathrm{lk}(c,a) = 1$, then $c$ represents a non-zero homology class in $H_k(H;\mathbb{Z}_2)$.
\end{lem}

\begin{sloppypar}
    \textbf{Proof:} Suppose that $c$ represents the trivial class in $H_k(H;\mathbb{Z}_2)$, \textit{i.e.} that there is a smooth $(k+1)$-chain $E$ in $H$ such that $\partial E = c$. Let $\hat{E} := E + d$. It is a $(k+1)$-cycle in $\mathbb{RP}^n$. Suppose that $a$ and $\hat{E}$ are transverse. 
    If $a$ is null-homologous in $\mathbb{RP}^n$, then the $\mathbb{Z}_2$-intersection product of $\hat{E}$ and $a$ is $0$. However, as $a \bigcap E = \emptyset$ and $\mathrm{lk}(c,a)=1$, the number of intersection points of $a$ and $\hat{E}$ is odd. We obtain a contradiction. 
    If $a$ is not null-homologous in $\mathbb{RP}^n$, then $a$ is a non-contractible loop in $\mathbb{RP}^n$ and $c$ is of dimension $n-2$. We consider the $(k+1)$-chain $\widetilde{E} = H + E + d$. We can slightly perturb $\hat{E}$ or $\widetilde{E}$ to make them transverse. If neither $\hat{E}$ nor $\widetilde{E}$ is null-homologous in $\mathbb{RP}^n$, then their intersection represents the non-trivial element in $H_{n-2}(\mathbb{RP}^n;\mathbb{Z}_2)$. However, $\hat{E} \bigcap \widetilde{E}$ lies in a neighbourhood of the disk $d$. Hence, $\hat{E} \bigcap \widetilde{E}$ is null-homologous in $\mathbb{RP}^n$ and either $\hat{E}$ or $\widetilde{E}$ is null-homologous in $\mathbb{RP}^n$. We obtain that the $\mathbb{Z}_2$-intersection product of either $\hat{E}$ or $\widetilde{E}$ with $a$ is $0$, which is again a contradiction. \hfill \qedsymbol \\
\end{sloppypar}

\begin{sloppypar}
    Definition \ref{deflink} can be extended bilinearly to the case in which $c$ is a disjoint union of smooth $k$-spheres bounding smooth $(k+1)$-disks in $\mathbb{RP}^n$ and $a$ is a disjoint union of smooth \mbox{$(n-k-1)$-spheres} in $\mathbb{RP}^n$. Lemma \ref{Lemmelinkindep} and its proof can be easily adapted to this situation. This generalization allows one to prove the following proposition.
\end{sloppypar}

\begin{prop}
    \label{indepmatrice}
    Let $k$ be an integer such that $0 \leq k \leq n-1$. Let $H$ be a smooth hypersurface in $\mathbb{RP}^n$. Let $(c_i)_{i=1,...,m}$ (respectively $(a_i)_{i=1,...,m}$) be smooth spheres of dimension $k$ (respectively, of dimension $n-k-1$) in $\mathbb{RP}^n$. Suppose that each $k$-sphere $c_i$ bounds a smooth $(k+1)$-disk in $\mathbb{RP}^n$. Suppose also that for every integer $i$ such that $1\leq i \leq m$, the $k$-sphere $c_i$ is contained in $H$ and the $(n-k-1)$-sphere $a_i$ is contained in $\mathbb{RP}^n \setminus H$. Let $\mathcal{A}$ be the square matrix whose coefficient in the $i$-th row and $j$-th column is $\mathrm{lk}(c_i,a_j)$. If the rank of $\mathcal{A}$ is maximal, then the family $(c_i)_{1\leq i \leq m}$ represents $m$ linearly independent classes in $H_k(H;\mathbb{Z}_2)$.
\end{prop}

\begin{sloppypar}
    \textbf{Proof:} Suppose that there exist some indices $i_1$,...,$i_p$ such that $\sum_{j=1}^p c_{i_j}$ is null-homologous in $H$. The sum of the corresponding rows of $\mathcal{A}$ is a non-zero vector, as the rank of $\mathcal{A}$ is maximal. Thus there exists an integer $q$ such that $1\leq q \leq m$ and $\mathrm{lk}(\sum_{j=1}^p c_{i_j}, a_q)=1$. Using the generalization of Lemma \ref{Lemmelinkindep}, we obtain a contradiction. \hfill \qedsymbol\\
\end{sloppypar}

\begin{sloppypar}
    In the rest of the paper, we use the notion of $\mathbb{Z}_2$-linking number in a slightly different setting. We consider a piecewise linear hypersurface $H$ in $\widetilde{(\Delta_m^n)^*}$ for some positive integers $m$ and $n$. The space $\widetilde{(\Delta_m^n)^*}$ is $PL$ homeomorphic to $\mathbb{RP}^n$. We also consider piecewise linear submanifolds of  $\widetilde{(\Delta_m^n)^*}$ that are $PL$ homeomorphic to spheres and that may bound piecewise linear disks in $\widetilde{(\Delta_m^n)^*}$.
    It is possible to extend Definition \ref{deflink}, Lemma \ref{Lemmelinkindep} and Proposition \ref{indepmatrice} to piecewise linear objects, as such objects are piecewise smooth. Transversality of a piecewise smooth disk $d$ and a piecewise smooth sphere $a$ in $\widetilde{(\Delta_m^n)^*}$ only makes sense when $d \bigcap a$ lies away from the non-smooth locus of $d$, $a$ and $\widetilde{(\Delta_m^n)^*}$. Whenever we describe objects that do not respect this condition, it is understood they can be slightly perturbed to obtain objects that do.
\end{sloppypar}

\subsection{Triangulations and sign distributions}

\label{triangulation}

\begin{sloppypar}
Let $x_1$, $x_2$, $x_3$ be the standard coordinates in $\mathbb{R}^3$ and let $m$ be a positive integer. We describe a subdivision process that can produce a variety of primitive convex triangulations of $\Delta_m^3$. A triangulation obtained through this process is called an \textit{IV triangulation of} $\Delta_m^3$. We subdivide $\Delta_m^3$ following the steps below.
\begin{itemize}
    \item Choose a $2$-face $T_m$ of $\Delta_m^3$. Let $k$ be an integer such that $1 \leq k \leq m-1$. Let $T_k$ be the triangle parallel to $T_m$, whose vertices are integer points contained in edges of $\Delta_m^3$ and whose edges are of lattice length $k$. Subdivide $\Delta_m^3$ using the triangles $T_k$. Denote by $T_0$ the vertex of $\Delta_m^3$ such that $T_0$ is opposite to $T_m$.
    \item Choose two different edges $E_1$ and $E_2$ of $\Delta_m^3$ that are not faces of $T_m$. For each integer $k$ such that $1 \leq k \leq m$ and $m-k$ is even (respectively, odd), consider the cone with basis $T_k$ and with vertex $E_1 \bigcap T_{k-1}$ (respectively $E_2 \bigcap T_{k-1}$). For each integer $k$ such that $0 \leq k \leq m-1$ and $m-k$ is even (respectively, odd), consider the cone with basis $T_k$ and with vertex $E_1 \bigcap T_{k+1}$ (respectively $E_2 \bigcap T_{k+1}$).
    \item Let $k$ be an integer such that $1 \leq k \leq m$. If $m-k$ is even (respectively, odd), we denote by $S_k$ the edge of $T_k$ that has no extremity lying on $E_1$ (respectively, on $E_2$). Divide $T_k$ into quadrangles using the segments parallel to $S_k$, whose extremities are integer points lying on edges of $T_k$. Choose a diagonal of one of the quadrangles. Its extremities are points of parity $p_1$ and $p_2$ respectively. Divide each quadrangle using its diagonal whose extremities are integer points of parity $p_1$ and $p_2$ respectively. Complete this subdivision of $T_k$ into a primitive triangulation in the only possible way.
    \item Let $k$ be an integer such that $1 \leq k \leq m$. Each cone over $T_k$ is naturally triangulated by the cone over the triangulation of $T_k$.
    \item For each $k$ such that $1 \leq k \leq m-1$, triangulate the join of $S_k$ and $S_{k+1}$ by taking the join of the respective triangulation of $S_k$ and $S_{k+1}$.
\end{itemize}
The IV triangulations of $\Delta_m^3$ are convex and primitive (see \cite{tsurf}).  A possible subdivision of $\Delta_4^3$ obtained after the first two subdivision steps is shown on Figure $1$.
\end{sloppypar}
\begin{sloppypar}
    Let $V(\Delta_m^3)$ be the set of integer points lying in $\Delta_m^3$. We identify the set $\{+,-\}$ with $\mathbb{Z}_2$. An \textit{IV sign distribution} on the vertices of $\Delta_m^3$ is a map $\mu_{IV}$ : $V(\Delta_m^3) \rightarrow \mathbb{Z}_2$ such that the signs of any two vertices of the same parity are the same. Let $k$ be an integer such that $1 \leq k \leq m$. 
    %Notice that the sum in $\mathbb{Z}_2$ of the signs associated to the four parities of integer points lying in $T_k$ is invariant under the procedure presented in Section \ref{patchwork} that extends $\mu_{IV}$ to the set of points lying in $(\Delta_m^3)^*$.
    When the sum in $\mathbb{Z}_2$ of the signs associated to the four parities of integer points lying in $T_k$ for a certain $k$ is $-$, we say that the restriction of $\mu_{IV}$ to $T_k$ is a \textit{Harnack distribution of signs on} $T_k$.
\end{sloppypar}

\begin{center}
    \includegraphics[scale=0.4]{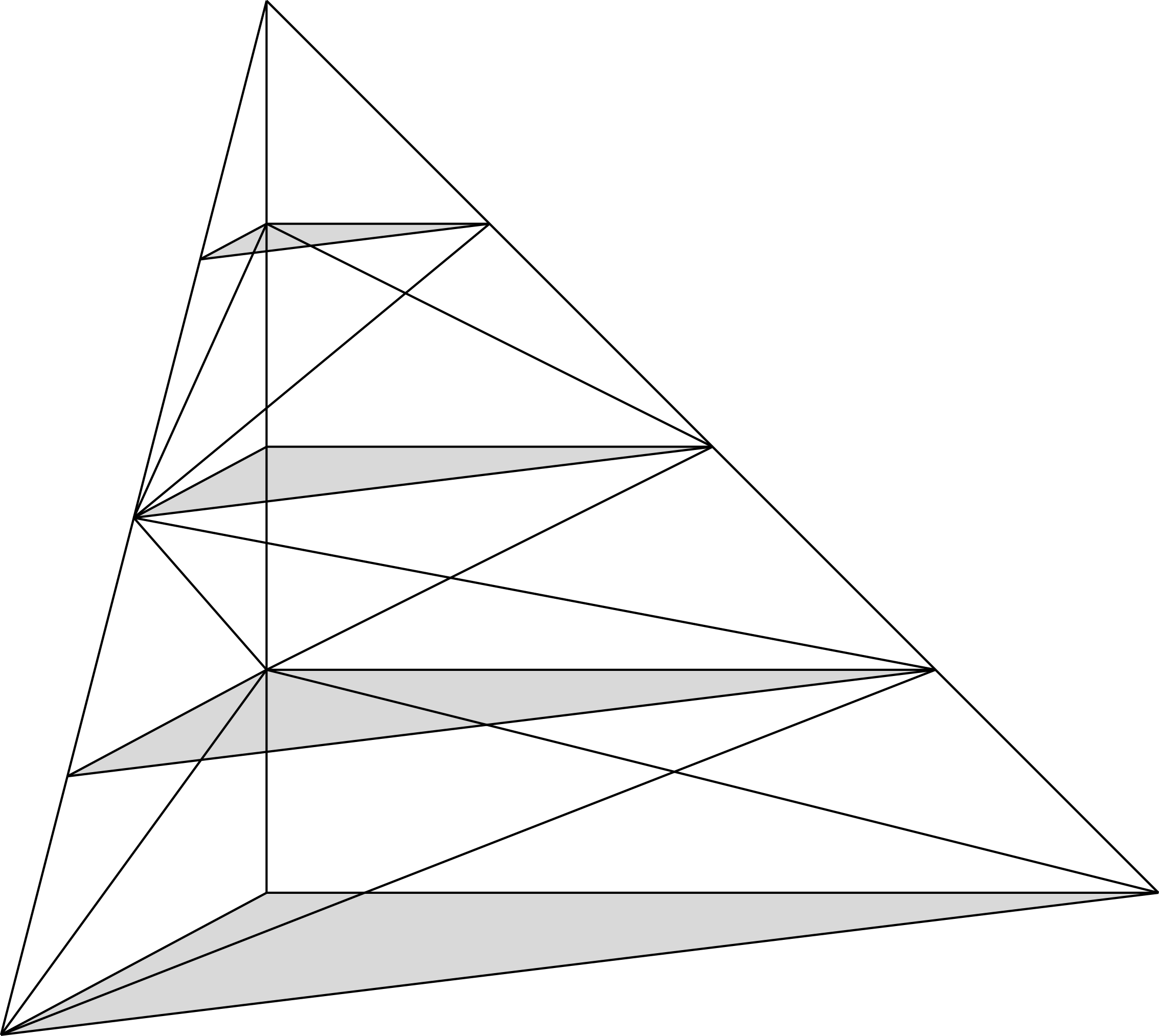}
    \begin{minipage}{16cm}
    \begin{it}
    \small{Fig. 1: A possible subdivision of $\Delta_4^3$ obtained after the first two subdivision steps described in \mbox{Section \ref{triangulation}}}
    \end{it}
    \end{minipage}
\end{center}

\subsection{Patchworking of IV-surfaces and cycle-axis pairs}

\label{result}

\begin{sloppypar}
    The following proposition is a special case of \cite[Proposition~5.1]{tsurf}. 
    \begin{prop}
    \label{resultmax}
        Applying the combinatorial patchworking theorem to an IV triangulation of $\Delta_m^3$ endowed with an IV sign distribution, we obtain a piecewise-linear surface $\widetilde{\Gamma}_m^{IV}$ in $\widetilde{(\Delta_m^3)^*}$, a non-singular real algebraic surface $X_m^{IV}$ in $\mathbb{RP}^3$ and a homeomorphism between the pairs $(\widetilde{(\Delta_m^3)^*}, \widetilde{\Gamma}_m^{IV})$ and $(\mathbb{RP}^3, \mathbb{R}X_m^{IV})$. The surface $X_m^{IV}$ is maximal and of type $\chi=\sigma$.
    \end{prop}
    A rather short proof of this proposition can be found in \cite{tsurf}. The primitivity of the triangulation allows one to compute the Euler characteristic of the surface $\widetilde{\Gamma}_m^{IV}$ directly. Indeed, given a simplex $Q$ of a primitive triangulation of $\Delta_m^3$, the number of symmetric copies of $Q$ whose vertices all bear the same sign (also referred to as \textit{the empty copies of} $Q$) only depends on the dimension of $Q$ and on its \textit{sedentarity}, \textit{i.e.} the number of $2$-faces of $\Delta_m^3$ the relative interior of $Q$ is contained in. It then suffices to count the number of connected components of the surface to deduce its Betti numbers. Unfortunately, this proof is not applicable to $3$-manifolds obtained through primitive combinatorial patchworking.\end{sloppypar}
    \begin{sloppypar}
    We give a much more intricated proof of the above statement using the tools presented in Section \ref{sectionlinkingnumber}. This proof is inspired by \cite{itenbergviro2007}. We exhibit explicit cycles that represent independent homology classes in the first and second homology groups of the surface $\widetilde{\Gamma}_m^{IV}$. These cycles are used later in Sections \ref{susprk} and \ref{suspother} to exhibit independent homology classes of some $3$-manifolds.
\end{sloppypar}

\label{properties}

\begin{sloppypar}
    We describe a collection of cycles contained in the piecewise-linear surface $\widetilde{\Gamma}^{IV}_m$. Each $j$-cycle $c$ is \mbox{$PL$ homeomorphic} to a $j$-sphere and comes with a \textit{dual axis} $a$, \textit{i.e.} a $(3-j-1)$-cycle $PL$ homeomorphic to a $(3-j-1)$-sphere and contained in $\widetilde{(\Delta_m^3)^*} \setminus \Tilde{\Gamma}^{IV}_m$. The \textit{cycle-axis pair} $(c,a)$ is called a $(j, 3-j-1)$\textit{-pair}, or a \textit{pair of type} $(j,3-j-1)$. We denote by $\mathcal{C}_{j}$ our collection of $(j,3-j-1)$-pairs and we define a partial order on $\mathcal{C}_{j}$. We call \textit{an orthant} a symmetric copy of $\mathbb{R}_{\geq 0}^n$. We want the cycle-axis pairs to satisfy the following conditions.
    \begin{itemize}
        \item Let $(c,a)$ be a $(j,3-j-1)$-pair from the collection. Then, the $\mathbb{Z}_2$-linking number of $c$ and $a$ is $1$.
        \item Let $(c,a)$ and $(c',a')$ be two distinct $(j,3-j-1)$-pairs such that $(c,a) < (c',a')$. Then the $\mathbb{Z}_2$-linking number of $c$ and $a'$ is $0$.
        \item All the cycle-axis pairs of type $(2,0)$ are such that the cycle is contained in a single orthant and the axis is contained in the union of at most two orthants adjacent to each other.
        \item Up to at most $\lfloor \frac{m}{2} \rfloor$ exceptions, all the cycle-axis pairs of type $(1,1)$ are such that the cycle or the axis is contained in the union of at most two orthants adjacent to each other.
    \end{itemize}
    
    To describe the cycle-axis pairs of the surface $\widetilde{\Gamma}^{IV}_m$, we need the following technical lemma, which is a reformulation of \cite[Proposition~3.1]{tsurf}. An orthant $O$ of $\mathbb{R}^n$ is identified with an $n$-tuple $b \in (\mathbb{Z}_2)^n$ as follows: the $i$-th coordinate of $b$ is $1$ if and only if the $i$-th coordinate of any point in $O$ is less than or equal to $0$. For any $a$, $b \in \mathbb{Z}_2^n$, we denote by $a.b$ the standard inner product of $a$ and $b$.
    \begin{lem}
        \label{linsystem}
        Let $Q \subset (\mathbb{R}_{\geq 0})^n$ be a primitive $k$-simplex with vertices $v_0$, ..., $v_{k}$ of respective parity $p_0$,...,$p_k$ and whose relative interior is contained in $(\mathbb{R}_{>0})^n$. Let $\mu$ be a sign distribution on the vertices of $Q$. Let $\epsilon \in \mathbb{Z}_2^{k+1}$. Then, there exists a unique collection $\mathcal{O} \subset \mathbb{Z}_2^n$ of cardinality $2^{n-k}$ such that for any $b \in \mathcal{O}$, either for each $i$ such that $0\leq i \leq n$, one has $\mu(v_i) +  p_i.b = \epsilon_i$ or for each $i$ such that $0\leq i \leq n$, one has $\mu(v_i) +  p_i.b = 1 + \epsilon_i$. \hfill \qedsymbol
        %\begin{itemize}
        %    \item either $\forall \: i$ such that $0\leq i \leq n$, one has $\mu(v_i) +  p_i.b = \epsilon_i$,
        %    \item or $\forall \: i$ such that $0\leq i \leq n$, one has $\mu(v_i) +  p_i.b = 1 + \epsilon_i$.     \hfill \qedsymbol
        %\end{itemize}
    \end{lem}
    In this lemma, $\epsilon$ represents a list of prescribed signs at the vertices of $Q$. For any $i$, in the context of combinatorial patchworking, the sign of the symmetric copy of $v_i$ lying in an orthant $b$ is given by $\mu(v_i) + p_i.b$. The lemma means that there are exactly $2^{n-k}$ symmetric copies of the primitive $k$-simplex $Q$ whose vertices bear the prescribed signs or their opposite. The proof of the lemma consists of simple linear algebra observations. It can be found in \cite{tsurf}.
\end{sloppypar}

\subsection{Pairs inside suspensions}

\label{surfsusp}

\begin{sloppypar}
Let $m$ be a positive integer, let $\tau$ be an IV triangulation of $\Delta_m^3$, and let $\mu$ be an IV distribution of signs on the vertices of $\tau$. We use the notations introduced in Section \ref{triangulation}. Let $k$ be an integer such that $1 \leq k\leq m-1$. The triangulation $\tau$ contains a union of simplices which is a suspension over $T_k$. We describe cycle-axis pairs lying in the symmetric copies of this suspension.
%An orthant $O$ of $\mathbb{R}^3$ is identified with an element $b \in \mathbb{Z}_2^3$ as follows: the $i$-th coordinate of $b$ is equal to $1$ if and only if the $i$-th coordinate of any point in $O$ is less than or equal to $0$. 
The symmetric copy of a simplex $Q$ in an orthant $O$ is called \textit{the} $O$-\textit{copy of} $Q$.
\end{sloppypar}

\subsubsection{Suspension pairs of type $(2,0)$}

\label{susp20}

\begin{sloppypar}
    Let $v$ be an integer point of parity $p_0$ lying in the relative interior of $T_k$. The star of $v$ in $\tau$ contains $6$ integer points different from $v$ of $3$ different parities $p_1$, $p_2$, $p_3$. Each tetrahedron in the star of $v$ has one vertex of each of these parities. We apply Lemma \ref{linsystem} to any tetrahedron in the star of $v$: up to inversion of all the signs, there exists a unique orthant $O$ such that the $O$-copy of $v$ bears the sign $+$ and the $O$-copies of the vertices of parities $p_1$, $p_2$ and $p_3$ bear the sign $-$. Hence, the star of the $O$-copy of $v$ contains a cell complex $PL$ homeomorphic to a sphere and belonging to $\Tilde{\Gamma}^{IV}_m$. This sphere is a $2$-cycle.
    The dual axis is the pair of points composed of the $O$-copy of $v$ and the $O$- or $O'$-copy of a triangulation vertex $v'$ lying in the boundary on $T_k$ but not in $S_k$, where $O'$ is an orthant adjacent to $O$ and such that the suspension over the $O$- and $O'$-copies of $T_k$ share the same vertices. We choose $v'$ and the orthant ($O$ or $O'$) such that the chosen copy of $v'$ bears the same sign as the $O$-copy of $v$.
    The cycle of such a $(2,0)$-pair is entirely contained in the orthant $O$. The axis may be contained in a union of two adjacent orthants. Given a pair $(c,a)$ obtained from the description above, the segment whose extremities are the points in $a$ meets $c$ in exactly one point. This segment meets any other described $2$-cycle in an even number of points. Hence the $\mathbb{Z}_2$-linking number of $a$ and $c$ is $1$ and, given another pair $(c',a')$ obtained from the description above, the $\mathbb{Z}_2$-linking number of $c$ and $a'$ and the $\mathbb{Z}_2$-linking number of $c'$ and $a$ are both $0$.
    \end{sloppypar}
    \begin{sloppypar}
    For each integer $k$ such that $1 \leq k \leq m-1$, we have described $\frac{(k-1)(k-2)}{2}$ pairs of type $(2,0)$. In total, this amounts to $\frac{(m-1)(m-2)(m-3)}{6}$ pairs of type $(2,0)$.
\end{sloppypar}

\subsubsection{Suspension pairs of type $(1,1)$ - a}

\label{11a}

\begin{sloppypar}
    Retain the notations from the previous Section \ref{susp20}. Two of the integer points in the star of $v$ do not lie in $T_k$. We assume that their parity is $p_3$. Let $Q$ be a tetrahedron in the star of $v$. The vertices $v$, $v_1$, $v_2$, $v_3$ of $Q$ are of respective parities $p_0$, $p_1$, $p_2$, $p_3$. We apply Lemma \ref{linsystem} to $Q$: up to inversion of all the signs, there exists a unique orthant $O$ in which the $O$-copies of $v$ and $v_3$ bear the sign $+$, while the $O$-copies of $v_1$ and $v_2$ bear the sign $-$. In the solution orthant $O$, take the union of the symmetric copies of the two segments of the triangulation whose extremities are $v$ and the considered points of parity $p_3$. This union is called the \textit{essential part of the axis}. We complete this union into a $1$-axis by taking its union with two segments whose extremities are the considered copies of the points of parity $p_3$ and a vertex of the triangulation chosen as in the axis description of the previous Section \ref{susp20}. In particular, the chosen copy of this vertex bears the same sign as the considered copies of the points of parity $p_3$ and lies on a copy of an edge of $T_k$ different from $S_k$. As $O$ is the solution orthant, the intersection of $\Tilde{\Gamma}^{IV}_m$, the star of the $O$-copy of $v$ and the $O$-copy of $T_k$ is a circle bounding a $2$-disk in the $O$-copy of $T_k$ containing the $O$-copy of $v$. The $O$-copy of $v$ is the only intersection point of this $2$-disk and the $1$-axis described above. Hence, the $\mathbb{Z}_2$-linking number of the described cycle and its dual axis is $1$. The cycle of such a $(1,1)$-pair is contained in a single orthant.
    \end{sloppypar}
    \begin{sloppypar}
    In this section, we have described $\frac{(k-1)(k-2)}{2}$ pairs of type $(1,1)$. In total, this amounts to $\frac{(m-1)(m-2)(m-3)}{6}$ pairs of type $(1,1)$.
\end{sloppypar}

\subsubsection{Suspension pairs of type $(1,1)$ - b}

\label{11b}

\begin{sloppypar}
    Retain the notations from the previous Sections \ref{susp20} and \ref{11a}. The integer point $v$ lies in the relative interior of a segment $S_v$ parallel to $S_k$ and such that the extremities of $S_v$ are integer points lying on the boundary of $T_k$. Suppose that the lattice length of $S_v$ is odd. The star of $v$ contains $2$ points that lie in $T_k$ but not in $S_v$. We assume that their parity is $p_2$. Let $Q$ be a tetrahedron in the star of $v$. The vertices $v$, $v_1$, $v_2$, $v_3$ of $Q$ are of respective parities $p_0$, $p_1$, $p_2$, $p_3$. We apply Lemma \ref{linsystem} to $Q$: up to inversion of all the signs, there exists a unique orthant $O$ such that the $O$-copies of $v$ and $v_2$ bear the sign $+$ while the $O$-copies of $v_1$ and $v_3$ bear the sign $-$. Take the union of
    \begin{itemize}
        \item the $O$-copies of the segments of the triangulation whose extremities are $v$ and the two points of parity $p_2$ lying in the star of $v$ (this union in the essential part of the axis);
        \item the $O$-copies of the segments of the triangulation whose extremities are the two points of parity $p_2$ lying in the star of $v$ and the copy of the extremity of $S_k$ that has the same parity as $v$.
    \end{itemize}
    This union is a $1$-axis. Lemma \ref{linsystem} ensures that there is a $1$-cycle of $\Tilde{\Gamma}^{IV}_m$ linked with the described axis and lying in the star of the $O$-copy of $v$. For instance, this $1$-cycle can be taken as the suspension having as basis the pair of points lying in the $O$-copy of the intersection of the star of $v$ with $S_v \bigcap \widetilde{\Gamma}_m^{IV}$ and with vertices the middle points of the $O$-copies of the triangulation segments whose extremities are $v$ and the vertices of the suspension over $T_k$. It bounds an obvious $2$-disk that meets the essential part of the dual axis in exactly one point, from which we deduce that the linking number of the described cycle with its dual axis is $1$. The shape of generic cycle-axis pairs described in Sections \ref{11a}, \ref{susp20}, \ref{11b} is represented in Figure $2$.
\end{sloppypar}

\begin{center}
    \includegraphics[scale=1]{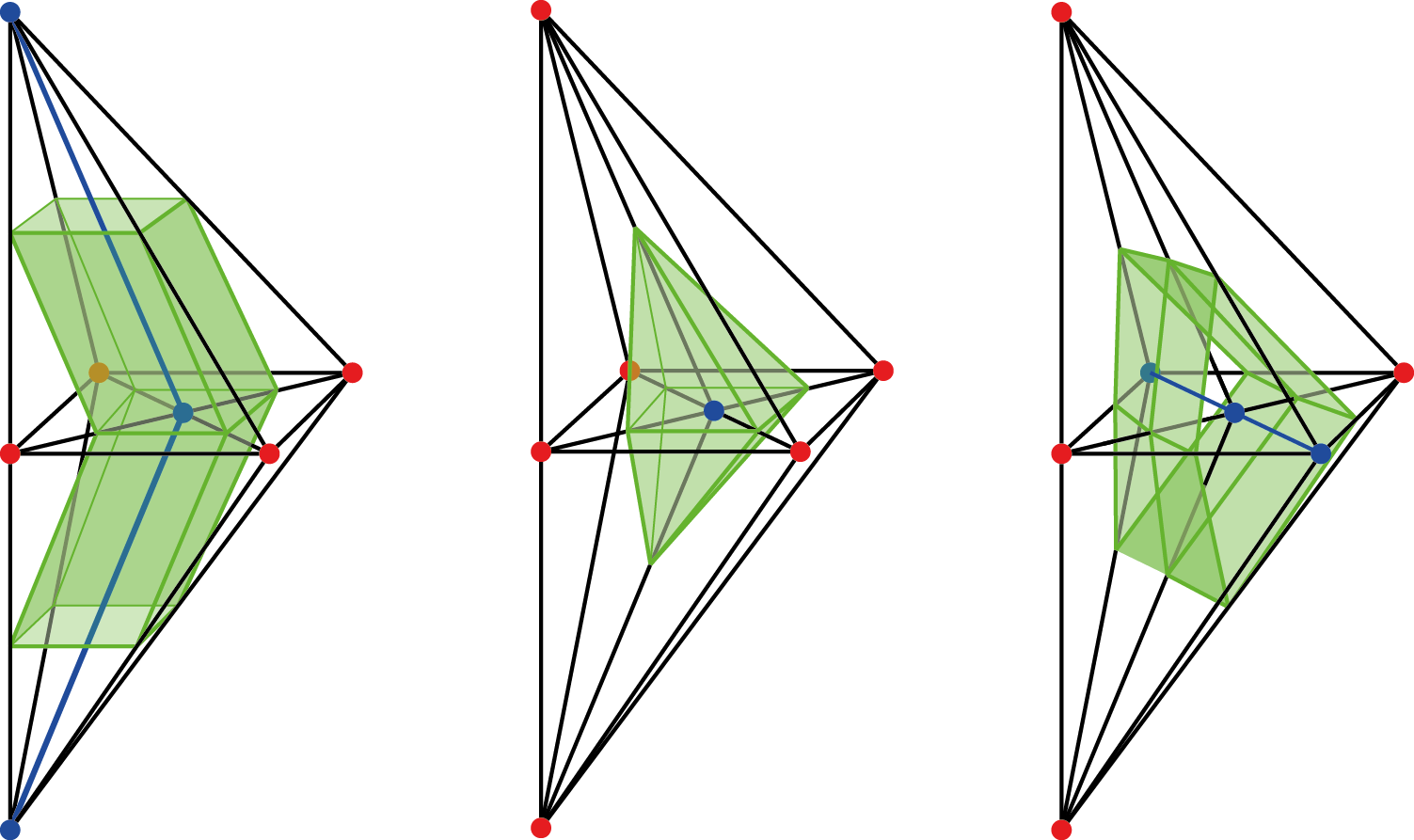}
    \begin{minipage}{16cm}
    \begin{it}
    \small{Fig. 2: In green, the intersection of the patchworked surface with the union of the stars of the essential parts of the axes of pairs from Sections \ref{11a}, \ref{susp20}, \ref{11b} (from left to right). The colors of the vertices represent their signs.}
    \end{it}
    \end{minipage}
\end{center}

\subsubsection{Suspension pairs of type $(1,1)$ - c}

\label{11c}

\begin{sloppypar}
    Let $v$ be an integer point of parity $p_0$ lying in the relative interior of $T_k$. It lies in the relative interior of a segment $S_v$ parallel to $S_k$. Suppose that the lattice length of this segment is even. We denote by $v_0$ the integer point on $S_v$ that is adjacent to the extremity of $S_v$ whose star contains exactly four $3$-simplices. We suppose that $v$ is distinct from $v_0$. The star of $v$ contains $2$ points that lie in $T_k$ but not in $S_v$. We assume that their parity is $p_2$. Let $Q$ be a simplex in the star of $v$. The vertices $v$, $v_1$, $v_2$, $v_3$ of $Q$ are of respective parity $p_0$, $p_1$, $p_2$, $p_3$. We apply Lemma \ref{linsystem} to $Q$: up to inversion of all the signs, there exists a unique orthant $O$ such that the $O$-copy of $v$ and $v_2$ bear the sign $+$ while $v_1$ and $v_3$ bear the sign $-$. Take the union of
    \begin{itemize}
        \item the $O$-copies of the segments of the triangulation whose extremities are $v$ and the two points of parity $p_2$ that lie in the star of $v$ (the union of these segments is the essential part of the axis);
        \item if the parity of $v$ is the same as the parity of the extremities of $S_v$, the $O$-copies of the segments whose extremities are the points of parity $p_2$ that lie in the star of $v$ and the extremity of $S_v$ whose star contains exactly four $3$-simplices;
        \item if the parity of $v$ is different from the parity of the extremities of $S_v$, the $O$-copies of the segments whose extremities are $v_0$ and the points of parity $p_2$ that lie in the star of $v$.
    \end{itemize}
    This union is a $1$-axis. Lemma \ref{linsystem} ensures that there is a $1$-cycle of $\Tilde{\Gamma}^{IV}_m$ linked with the described axis and lying in the star of the suitable copy of $v$. As in the previous section, this $1$-cycle can be taken as the suspension having as basis the pair of points lying in the $O$-copy of the intersection of the star of $v$ with $S_v \bigcap \widetilde{\Gamma}_m^{IV}$ and with vertices the middle points of the $O$-copies of the triangulation segments whose extremities are $v$ and the vertices of the suspension over $T_k$.\\
\end{sloppypar}

\subsubsection{Suspension pairs of type $(1,1)$ - d }

\label{11d}

\begin{sloppypar}
    Let $v$ be an extremity of a segment $S_v$ of even lattice length contained in $T_k$, parallel to $S_k$ and different from $S_k$ such that the star of $v$ does not contain exactly four $3$-simplices. Let $A$ be the triangulation segment whose relative interior lies in the relative interior of $T_k$, that has $v$ as an extremity and whose other extremity is an integer point $f$ that is the extremity of a segment with integer extremities contained in the edges of $T_k$, parallel to $S_k$ and of lattice length less than the lattice length of $S_v$. We denote by $v'$ the triangulation vertex that is the closest to $v$ in $S_v$. Let $Q$ be a tetrahedron in the star of $v'$ that also has $v$ as a vertex. Let $p_0$, $p_1$, $p_2$ and $p_3$ denote the respective parities of the vertices $v'$, $v$, $v_2$, $v_3$ of $Q$, where $v_3$ is the vertex of $Q$ that does not belong to $T_k$. We apply Lemma \ref{linsystem} to $Q$: up to inversion of all the signs, there exists a unique orthant $O$ such that the $O$-copies of $v'$ and $v_3$ bear the sign $+$ while the $O$-copies of $v$ and $v_2$ bear the sign $-$. We complete the $O$-copy of $A$ into a $1$-axis by taking its union with the copies lying in $O$ of:
    \begin{itemize}
        \item the triangulation segment $A'$ whose relative interior lies in the relative interior of $T_k$ and whose extremities are $v$ and the extremity of a segment with integer extremities contained in $T_k$, parallel to $S_k$, of lattice length greater than the lattice length of $S_v$ (the $O$-copy of $A'$ is the essential part of the axis);
        \item the segment whose extremities are the extremities of $A$ and $A'$ that are different from $v$. This segment is the union of two triangulation segments.
    \end{itemize}
    There is exactly one edge $E_{p_0}$ of $T_k$ that contains triangulation vertices of the same parity as $v'$. It is different from the edge of $T_k$ that contains $f$. Take the union of the following segments:
    \begin{itemize}
        \item the $O$-copy of $[v',f]$;
        \item the $O$-copy of the segment whose extremities are $f$ and the extremity $e$ of $E_{p_0}$ lying on the same edge of $T_k$ as $f$;
        \item the $O$-copy of the segment whose extremties are $e$ and the closest triangulation vertex whose parity is the same as the parity of $v'$.
    \end{itemize}
    This broken line $C$ meets $\Tilde{\Gamma}^{IV}_m$ in at least two points. Take the broken line $\Bar{C}$ that is a subset of $C$ and whose extremities are the intersection points of $C$ with $\Tilde{\Gamma}^{IV}_m$ that are the closest to the $O$-copy of $v'$ on $C$. As the sign of the $O$-copies of the vertices of the suspensions over $T_k$ are different from the sign of the integer points lying in $\Bar{C}$, there is a double cover of the relative interior of $\Bar{C}$ contained in the intersection of $\Tilde{\Gamma}^{IV}_m$ with the union of the suspensions over the simplices of $C$. The sheets of the cover meet at the extremities of $\Bar{C}$, forming the $1$-cycle dual to the $1$-axis described above. From the fact that $\Bar{C}$ meets the $O$-copy of $A \bigcup A'$ in exactly one point, we deduce that the $\mathbb{Z}_2$-linking number of a described cycle with its dual axis is $1$.\\
    Figure 3 shows the form of the described broken line $C$ and axis in the case $k=4$ and with two different sign distributions.
    \begin{center}
    \includegraphics[scale=0.35]{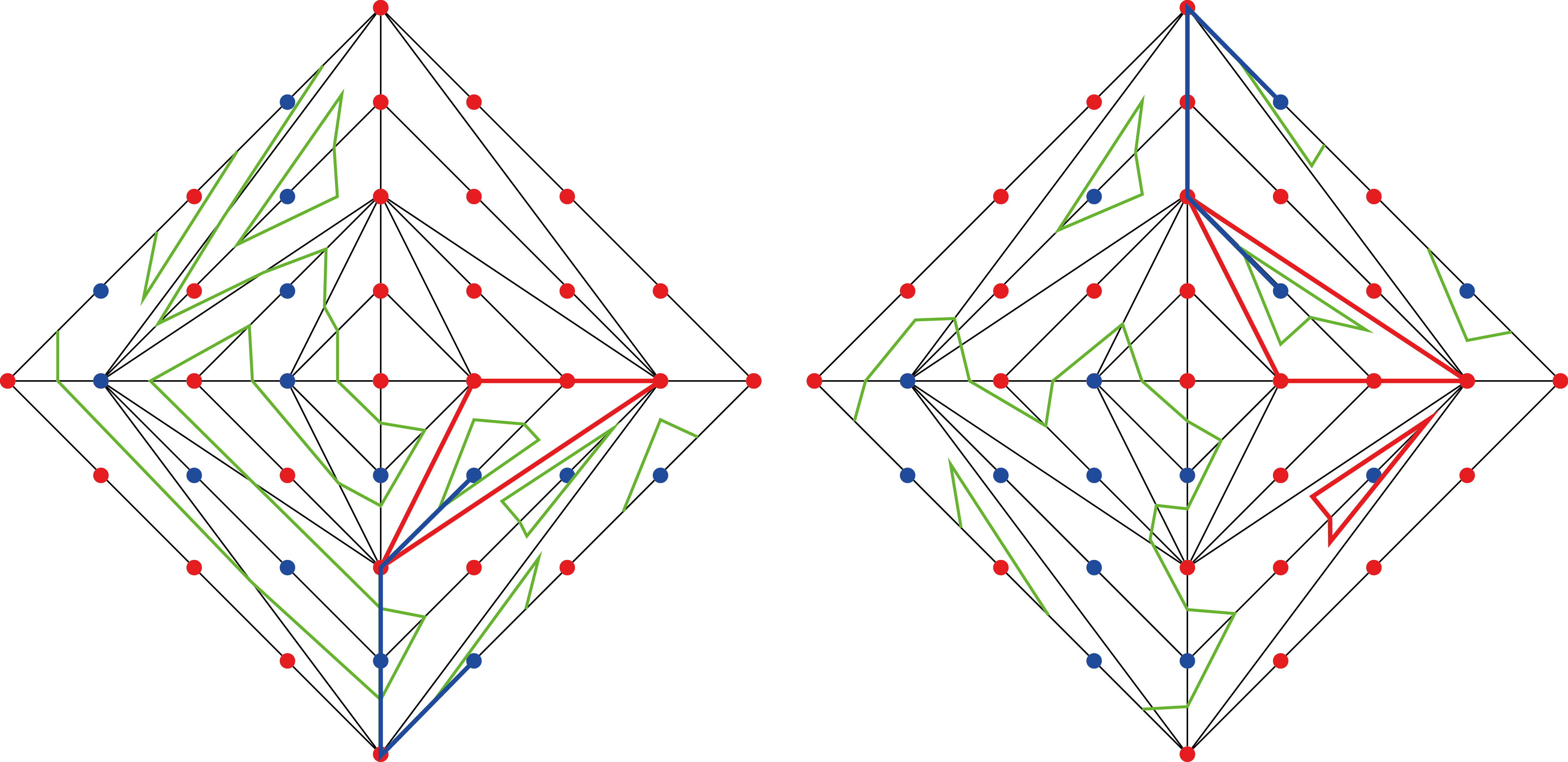}
    \begin{minipage}{16cm}
    \begin{it}
    \small{Fig. 3: Four copies of $T_4$ are represented. We suppose that the vertices of the suspensions over these four copies are the same. On the left, the integer points in $T_4$ are endowed with a non-Harnack sign distribution. On the right, the sign distribution is a Harnack one. The colors of the vertices represent their sign. The intersection of the patchworked surface with the considered copies of $T_4$ appears in green. The red broken lines are axes corresponding to the description of Section \ref{11d}. Each blue broken line corresponds to a broken line named $C$ from Section \ref{11d}.}
    \end{it}
    \end{minipage}
    \end{center}
    In Sections \ref{11b} \ref{11c} and \ref{11d}, for a fixed integer $k$ such that $1\leq k \leq m$, we have exhibited one $(1,1)$-pair for each integer point in the relative interior of $T_k$, that is $\frac{(k-1)(k-2)}{2}$ pairs of type $(1,1)$. In total, this amounts to $\frac{(m-1)(m-2)(m-3)}{6}$ pairs of type $(1,1)$. All $(1,1)$-pairs presented up to this point are called \textit{regular suspension pairs}.
\end{sloppypar}

\subsubsection{Special suspension $(1,1)$-pair}

\label{suspspecial}

\begin{sloppypar}
    These pairs are somewhat similar to the $(1,1)$-pairs described in Section \ref{11d}. Let $k$ be an even integer such that $1 \leq k \leq m-1$. Let $v$ be the extremity of $S_k$ whose star does not contain exactly two $3$-simplices. Let $A$ be the triangulation segment whose extremities are $v$ and an integer point that lies on another edge of $T_k$ and that is the extremity of a segment with integer extremities contained in $T_k$, parallel to $S_k$, of lattice length less than the lattice length of $S_k$. We denote by $v'$ the triangulation vertex that is the closest to $v$ in $S_k$. Let $Q$ be a tetrahedron in the intersection of the star of $v'$ with the suspension over $T_k$. We assume that $v$ is a vertex of $Q$. Let $p_0$, $p_1$, $p_2$ and $p_3$ denote the respective parities of the vertices $v'$, $v$, $v_2$, $v_3$ of $Q$, where $v_3$ is a vertex of the suspension over $T_k$. We apply Lemma \ref{linsystem} to $Q$: up to inversion of all signs, there exists a unique orthant $O$ such that the $O$-copies of $v'$ and $v_3$ bear the sign $+$, while the $O$-copies of $v$ and $v_2$ bear the sign $-$. To form a $1$-axis, take the union of the four copies of $A$ with extremities bearing the same sign. To form the dual $1$-cycle, we have to dissociate two cases.
    \end{sloppypar}
    \begin{sloppypar}
    First, when the restriction to $T_k$ of the sign distribution $\mu_{IV}$ is not a Harnack distribution, take the union of the $O$-copy of $[v',v]$ and the segment whose extremities are the $O$-copy of $v$ and the $O$-copy of the point closest to $v$ on the edge of $T_k$ that is different from $S_k$ and has $v$ as an extremity.
    %the following segments:
    %\begin{itemize}
    %    \item the $O$-copy of $[v',v]$;
    %    \item the segment whose extremities are the $O$-copy of $v$ and the $O$-copy of the point closest to $v$ on the boundary segment that is different from $S_k$ and has $v$ as an extremity.
    %\end{itemize}
    This broken line $C$ meets $\Tilde{\Gamma}^{IV}_m$ in exactly two points. Remove from this broken line the segments whose extremities are an extremity of $C$ and the closest intersection point of $C$ and $\Tilde{\Gamma}^{IV}_m$. Denote by $\Bar{C}$ the resulting broken line. As in the previous section, the intersection of $\Tilde{\Gamma}^{IV}_m$ with the $O$-copy of the suspension over $T_k$ contains a double cover of the relative interior of $\Bar{C}$ whose sheets meet at the extremities of $\Bar{C}$, forming the $1$-cycle dual to the $1$-axis described above. 
    %This $1$-cycle is contained in a single orthant. 
    From the fact that $\Bar{C}$ meets the $O$-copy of $A$ in exactly one point we deduce that the $\mathbb{Z}_2$-linking number of the described cycle with its dual axis is $1$.
    \end{sloppypar}
    \begin{sloppypar}
    When the restriction to $T_k$ of the sign distribution $\mu_{IV}$ is a Harnack distribution, take the union of the following segments:
    \begin{itemize}
        \item the $O$-copy of $[v,v']$;
        \item the $O$-copy of the edge of $T_k$ that is different from $S_k$ and that contains $v$;
        \item the segment whose extremities are the $O$-copy of the vertex $v''$ of $T_k$ that does not lie on $S_k$ and the $O'$-copy of a triangulation vertex adjacent to $v''$ in the boundary of $T_k$ where $O'$ is an orthant adjacent to $O$ and such that the suspension over the $O$- and the $O'$-copy of $T_k$ share the same vertices.
        %(possibly after identifying the antipodal points on the boundary of $\Delta_m^3$) ??
    \end{itemize}
    The broken line $C$ meets $\Tilde{\Gamma}^{IV}_m$ in exactly two points. Remove from this broken line the segments whose extremities are an extremity of $C$ and the closest intersection point of $C$ and $\Tilde{\Gamma}^{IV}_m$. Denote by $\Bar{C}$ the resulting broken line. As in the previous section, the intersection of $\Tilde{\Gamma}^{IV}_m$ with the union of the $O$- and $O'$-copies of the suspension over $T_k$ contains a double cover of the relative interior of $\Bar{C}$ whose sheets meet at the extremities of $\Bar{C}$, forming the $1$-cycle dual to the $1$-axis described above. From the fact that $\Bar{C}$ meets the $O$-copy of $A$ in exactly one point we deduce that the $\mathbb{Z}_2$-linking number of the described cycle with its dual axis is $1$. 
    %Notice that in this case, neither the cycle nor the axis are contained in a single orthant. This default is responsible for at most $\lfloor \frac{m-1}{2} \rfloor$ of the possible exceptions mentioned in Section \ref{properties}. 
    The pairs described in the present section are called \textit{special suspension pairs}. When the sign distribution on $T_k$ is a Harnack one, neither the axes nor the cycles of these pairs are contained in a single orthant. This default is responsible for at most $\lfloor \frac{m-1}{2} \rfloor$ of the possible exceptions mentioned in Section \ref{properties}. Figure 4 shows the form of the described broken line $C$ and of the described axis in the case $k=4$, with two different sign distributions. 
    \begin{center}
    \includegraphics[scale=0.35]{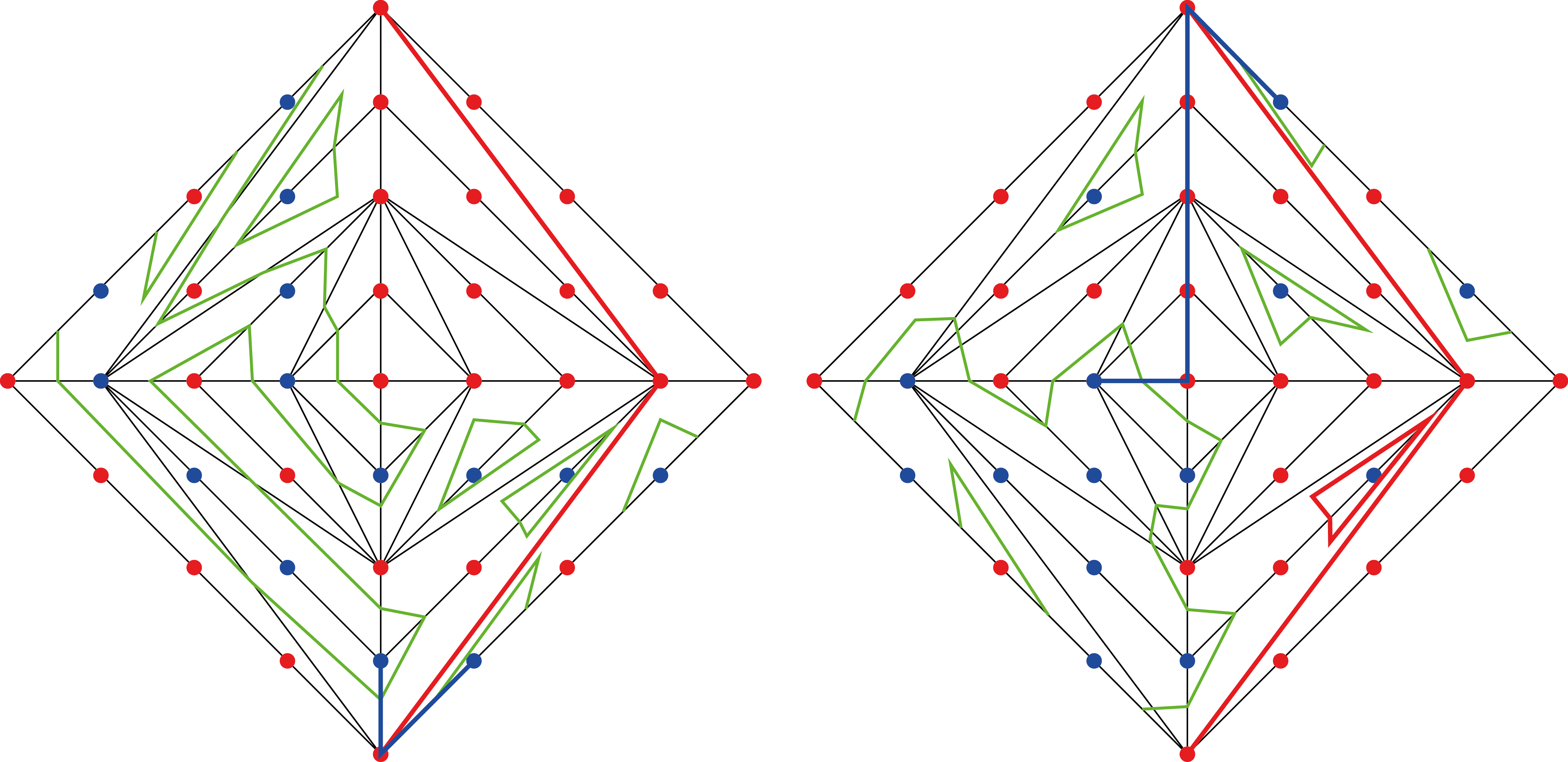}
    \begin{minipage}{16cm}
    \begin{it}
    \small{Fig. 4: Four copies of $T_4$ are represented. We suppose that the vertices of the suspensions over these four copies are the same. On the left, the integer points in $T_4$ are endowed with a non-Harnack sign distribution. On the right, the sign distribution is a Harnack one. The colors of the vertices represent their sign. The intersection of the patchworked surface with the considered copies of $T_4$ appears in green. The red broken lines are two segments from axes corresponding to the description of Section \ref{suspspecial} (the two other segments of such axis lie in the other four copies of $T_4$). Each blue broken line corresponds to a broken line named $C$ from Section \ref{suspspecial}.}
    \end{it}
    \end{minipage}
    \end{center}
    \begin{rmk}
    \label{rmksuspspecial}
        One sees from the description of the broken line $C$ that when the restriction to $T_k$ of the sign distribution $\mu_{IV}$ is not a Harnack distribution, the described cycle is contained in a single orthant; otherwise, the described cycle is contained in the union of two adjacent orthants.
        %\begin{itemize}
        %    \item when the restriction to $T_k$ of the sign distribution $\mu_{IV}$ is not a Harnack distribution, the described cycle is contained in a single orthant;
        %    \item otherwise, the described cycle is contained in the union of two adjacent orthants.
        %\end{itemize}
    \end{rmk}
\end{sloppypar}

\subsection{Pairs of type $(1,1)$ inside joins of two segments}

\label{join}

\begin{sloppypar}
    Let $k$ be an integer such that $1\leq k \leq m-1$. The join of $S_k$ and $S_{k+1}$ appears as a union of simplices of the IV triangulation $\tau$. Let $a_0$ (respectively $a_1$) be an integer point lying in the relative interior of $S_k$ (respectively $S_{k+1}$). Let $p_0$ and $p_1$ be the respective parities of $a_0$ and $a_1$. Let $A$ be the triangulation segment with extremities $a_0$ and $a_1$. Aside from $a_0$ and $a_1$, there are $4$ integer points in the star of $A$ (the neighbours of $a_0$ and $a_1$ on the segments $S_k$ and $S_{k+1}$ respectively), of two different parities $p_2$ and $p_3$. Let $Q$ be a tetrahedron in the star of $A$ with vertices $a_0$, $a_1$, $v_2$, $v_3$ of respective parities $p_0$, $p_1$, $p_2$, $p_3$. We apply Lemma \ref{linsystem} to $Q$: up to inversion of all the signs, there exists a unique orthant $O$ such that the $O$-copies of $a_0$ and $a_1$ bear the sign $+$ while the $O$-copies of $v_2$ and $v_3$ bear the sign $-$. To form a $1$-axis, first take the union of
    \begin{itemize}
        \item the $O$-copy of the segment $A$ (it is the essential part of the axis);
        \item if $k$ is even (respectively, odd), the $O$-copy of the segment whose extremities are $a_0$ (respectively, $a_1$) and the extremity $e_1$ of $S_{k+1}$ of parity $p_1$ (respectively, the extremity of $S_{k}$ of parity $p_0$);
        %\item if $k$ is odd, the copy in $O$ of the segment whose extremities are $a_1$ and the extremity of $S_{k}$ of parity $p_0$;
        \item if $k$ is even (respectively, odd), the copy either in $O$ or in the orthant adjacent to $O$ and containing the $O$-copy of $S_{k+1}$ (respectively, $S_k$) of two triangulation segments lying in $T_{k+1}$ (respectively, $T_k$) with extremities $a_1$, $e_1$ and the vertex $v$ of the cone over $S_{k+1}$ in $T_{k+1}$ (respectively, $a_0$, $e_0$ and the vertex $v$ of the cone over $S_{k}$ in $T_{k}$.). The reflection or composition of reflections that interchanges the two distinguished orthants changes the sign of the corresponding copy of $v$. We take the orthant in which the copies of $a_1$, $e_1$ and $v$ (respectively, $a_0$, $e_0$ and $v$) have the same sign;
        %\item if $k$ is odd, the copy either in $O$ or in the orthant adjacent to $O$ and containing the $O$-copy of $S_{k}$ of two triangulation segments lying in $T_{k}$ with extremities $a_0$, $e_0$ and the vertex $v$ of the cone over $S_{k}$ in $T_{k}$. The reflection or composition of reflections that interchanges the two distinguished orthants changes the sign of the corresponding copy of $v$. We take the orthant in which the copies of $a_0$, $e_0$ and $v$ have the same sign.
    \end{itemize}
    Lemma \ref{linsystem} ensures that there is a $1$-cycle contained in the patchworked surface, lying in the star of the $O$-copy of $A$ and linked with the described $1$-axis. For instance, we can take it to be the join of two pairs of points of $\Tilde{\Gamma}^{IV}_m$. The first pair of points is the union of the middle points of the $O$-copies the segments adjacent to $a_0$ in $S_k$ and the second one is the union of the middle points of the $O$-copies of the segments adjacent to $a_1$ in $S_{k+1}$. Such a cycle is contained in a single orthant. It bounds an obvious $2$-disk that meets the essential segment of the axis in exactly one point, ensuring that the $\mathbb{Z}_2$-linking number of the described cycle with its dual axis is $1$.\\
    For a fixed $k$, we can construct $(k-1)k$ pairs of type $(1,1)$ associated with the join of $S_k$ and $S_{k+1}$ following this description. These pairs are called \textit{join pairs}. In total, this amounts to $\frac{m^3}{3} -m^2 + \frac{2m}{3}$ pairs of type $(1,1)$.
    The shape of a generic cycle-axis pair described in the current Section is represented in Figure $5$.
\end{sloppypar}

\begin{center}
    \includegraphics[scale=1]{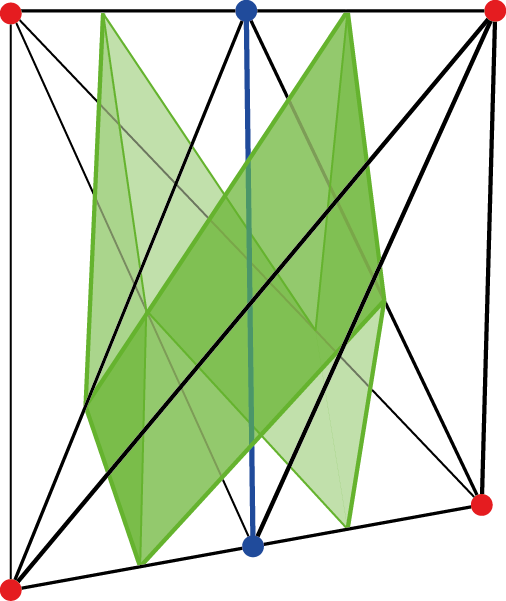}
    \begin{minipage}{16cm}
    \begin{it}
    \small{Fig. 5: In green, the intersection of the patchworked surface with the star of the essential segment of an axis from a pair described in \ref{join}. The colors of the vertices represent their signs.}
    \end{it}
    \end{minipage}
\end{center}

\subsection{Pairs of type $(1,1)$ in a $2$-dimensional wall between a suspension and a join}

\subsubsection{Regular wall pairs}

\label{wallregu}

\begin{sloppypar}
    Let $k$ be an integer such that $2 \leq k \leq m$. The IV triangulation $\tau$ of $\Delta_m^3$ contains a union of simplices forming a cone with basis $T_k$ and with vertex an integer point lying on $T_{k-1}$. If $k\leq m-1$, the triangulation also contains a union of simplices that is a cone with basis $T_k$ and with vertex an integer point lying on $T_{k+1}$. Consider such a cone and denote by $v$ its vertex. Let $e$ be an integer point lying in the relative interior of $S_k$ and at lattice distance $1$ from one of the extremities of $S_k$. From now on, we suppose that $k$ is greater than 2. Choose an integer point $f$ lying in the relative interior of the segment $S_k$ and distinct from $e$ if $k$ is even. The segment $A:=[v,f]$ is a segment of the triangulation. Let $p_0$ and $p_1$ be the respective parities of $v$ and $f$. Aside from $v$ and $f$, there are $3$ integer points of $2$ different parities $p_2$ and $p_3$ in the intersection of the star of $A$ with the cone over $T_k$. Let $Q$ be a tetrahedron in this intersection with vertices $v$, $f$, $v_2$, $v_3$ of respective parities $p_0$, $p_1$, $p_2$, $p_3$. We apply Lemma \ref{linsystem} to $Q$: up to inversion of all signs, there exists a unique orhtant $O$ such that the $O$-copy of $v$ and $f$ bear the sign $+$ while the $O$-copies of $v_2$ and $v_3$ bear the sign $-$. To form a $1$-axis, consider the union of:
    \begin{itemize}
        \item the $O$-copy of $A$ (the essential part of the axis);
        \item if $k$ is odd, the $O$-copy of the segments with vertices $v$ and the extremity $e_1$ of $S_k$ of parity $p_1$;
        \item if $k$ is even and the extremities of $S_k$ are of parity $p_1$, the $O$-copy of the segments with vertices $v$ and a fixed extremity $e_1$ of $S_k$;
        \item if $k$ is even and the extremities of $S_k$ are not of parity $p_1$, the $O$-copy of the segment with vertices $v$ and $e$.
        \item if $k$ is odd or $k$ is even and the extremities of $S_k$ are of parity $p_1$, the copy in the orthant $O'$ adjacent to $O$ and containing the $O$-copy of $S_{k}$ of two triangulation segments lying in $T_{k}$ with extremities $f$, $e_1$ and the vertex $w$ of the cone over $S_{k}$ in $T_{k}$. The copies of $f$, $e_1$ and $w$ lying in this orthant bear the same sign. Indeed, by Lemma \ref{linsystem}, the sign of $O$-copy of $w$ is opposite to the signs of the $O$-copies of $f$ and $e_1$ and the composition of reflections that exchanges $O$ and $O'$ changes the sign of $w$ but not the signs of $f$ and $e_1$.
        \item if $k$ is even and the extremities of $S_k$ are not of parity $p_1$, the copy in the orthant $O'$ adjacent to $O$ and containing the $O$-copy of $S_{k}$ of two triangulation segments lying in $T_{k}$ with extremities $f$, $e_1$ and the vertex $w$ of the cone over $S_{k}$ in $T_{k}$. As above, the $O'$-copies of $a$, $e$ and $w$ bear the same sign.
    \end{itemize}
    Now, we describe the $1$-cycle dual to such a $1$-axis. Lemma \ref{linsystem} ensures that the intersection of the star of the $O$-copy of $A$ with the $O$-copy of the cone over $T_k$ with vertex $v$ contains an arc $\alpha$ of the patchworked surface whose intersection with the $O$-copy of the cone of vertex $v$ over $S_k$ is reduced to the middle points of the $O$-copies of the triangulation edges adjacent to $f$ and contained in $S_k$. There is only one triangulation vertex $u$ lying in the closure of the star of $A$ but not in the cone over $T_k$. Depending on the sign of this vertex, there are two ways of completing this arc into a $1$-cycle linked with the $1$-axis described above.
    If the sign of the $O$-copy of $u$ is different from the sign of the $O$-copy of $v$, then there is an arc of the patchworked surface contained in the $O$-copy of the star of $A$ and disjoint from the interior of the $O$-copy of the cone over $T_k$ with vertex $v$ whose intersection with the $O$-copy of the cone of vertex $v$ over $S_k$ is reduced to the middle points of the $O$-copies of the triangulation segments adjacent to $f$ and contained in $S_k$. The union of this arc with $\alpha$ is the desired $1$-cycle. In this case, the cycle is entirely contained in a single orthant.
    If the sign of the $O$-copy of $u$ is the same as the sign of the $O$-copy of $v$, we consider the orthant $O'$ adjacent to $O$ and containing the $O$-copy of $S_k$. The signs of the $O'$-copies of the integer points in $S_{k+1}$ alternate. Let $p$ be the middle point of the segment adjacent to $v$ in $S_{k+1}$. The cone over the $O'$-copies of the middle of the triangulation segments adjacent to $f$ and contained in $S_k$ with vertex the $O'$-copy of $p$ is contained in the patchworked surface. The union of this arc with $\alpha$ is the desired $1$-cycle. In this case, the described cycle is not contained in a single orthant.
    %\begin{itemize}
    %    \item If the sign of the $O$-copy of $u$ is different from the sign of the $O$-copy of $v$, then there is an arc of the patchworked surface contained in the $O$-copy of the star of $A$ and disjoint from the interior of the $O$-copy of the cone over $T_k$ with vertex $v$ whose intersection with the $O$-copy of the cone of vertex $v$ over $S_k$ is reduced to the middle points of the $O$-copies of the triangulation segments adjacent to $a$ and contained in $S_k$. The union of this arc with $\alpha$ is the desired $1$-cycle. In this case, the cycle is entirely contained in a single orthant.
    %    \item If the sign of the $O$-copy of $u$ is the same as the sign of the $O$-copy of $v$, we consider the orthant $O'$ adjacent to $O$ and containing the $O$-copy of $S_k$. The signs of the $O'$-copies of the integer points in $S_{k+1}$ alternate. Let $p$ be the middle point of the segment adjacent to $v$ in $S_{k+1}$. The cone over the $O'$-copies of the middle of the triangulation segments adjacent to $a$ and contained in $S_k$ with vertex the $O'$-copy of $p$ is contained in the patchworked surface. The union of this arc with $\alpha$ is the desired $1$-cycle. In this case, the described cycle is not contained in a single orthant.
    %\end{itemize}
    \end{sloppypar}
    \begin{sloppypar}
    Notice that the axes described in this section are never contained in a single orthant. Moreover, as noticed above, some of the cycles are not contained in a single orthant. For each cycle-axis pair such that the cycle is not contained in a single orthant, we modify the axis as follows. First, delete from the axis the segments that have the $O$-copy of $v$ as an extremity. There are two such segments, let us denote by $v_1$, $v_2$ their respective extremities that are different from $v$. Then, take the union of the rest of the axis with the $O$-copy of the triangulation segments with extremities $v_1$, $v_2$ and the extremity of $S_{k+1}$ that is different from $v$. These two segments lie in the intersection of the orthant $O$ with the orthant that contains the rest of the axis.
    %\begin{itemize}
    %    \item Delete from the axis the segments that have the $O$-copy of $v$ as an extremity. There are two such segments, let us denote by $v_1$, $v_2$ their respective extremities that are different from $v$.
    %    \item Take the union of the rest of the axis with the $O$-copy of the triangulation segments with extremities $v_1$, $v_2$ and the extremity of $S_{k+1}$ that is different from $v$. These two segments lie in the intersection of the orthant $O$ with the orthant that contains the rest of the axis.
    %\end{itemize}
    A cycle obtained from the description above bounds an obvious $2$-disk that meets the essential segment of its axis in exactly one point (the $O$-copy of $f$), ensuring the desired linking relations. When $k$ is odd, we have described one $(1,1)$-pair for each triangulation segment whose relative interior is contained in the relative interior of the cone of vertex $v$ over $S_k$. When $k$ is even, we have described one $(1,1)$-pair for each triangulation segment whose relative interior is contained in the relative interior of the cone of vertex $v$ over $S_k$ except one. The pairs presented in the present section are called \textit{regular wall pairs}. Figures 6 and 7 represent different examples of such pairs, with $k$ odd and even.
    \begin{center}
    \includegraphics[scale=0.7]{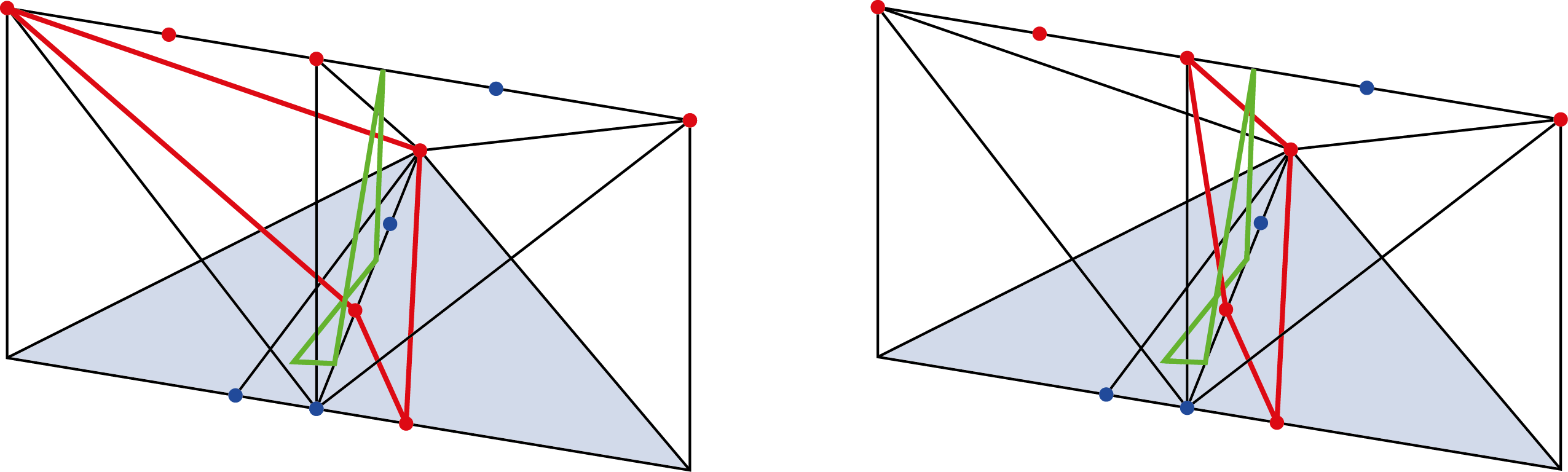}
    \begin{minipage}{16cm}
    \begin{it}
    \small{Fig. 6: On the left and on the right, two adjacent copies of a cone over $T_3$ and of a join of $S_3$ and $S_2$ are represented. The colors of the integer points (blue or red) represent their sign. On the left, a cycle axis pair as described in Section \ref{wallregu}. The cycle is green, the axis is red. Neither the axis nor the cycle is contained in a single orthant. On the right, the pair after the axis modification explained in Section \ref{wallregu}. The axis is now contained in a single orthant.}
    \end{it}
    \end{minipage}
\end{center}
    \begin{center}
    \includegraphics[scale=0.7]{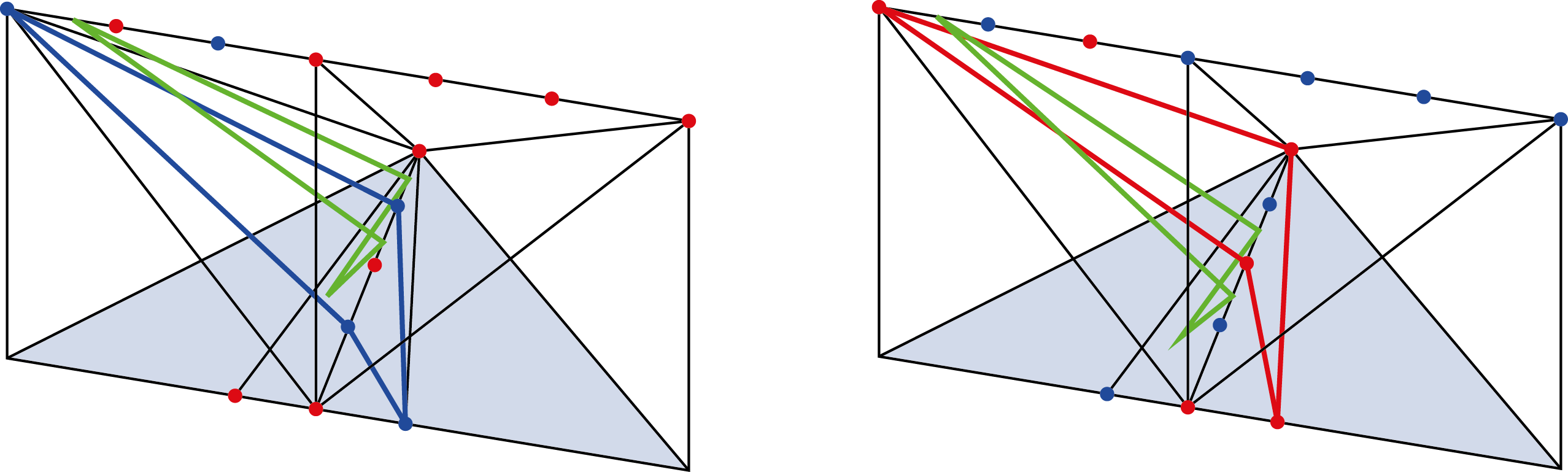}
    \begin{minipage}{16cm}
    \begin{it}
    \small{Fig. 7: On the left and on the right, two adjacent copies of a cone over $T_4$ and of a join of $S_4$ and $S_3$ are represented. The colors of the integer points (blue or red) represent their sign. On the left and on the right, a cycle-axis pair as described in Section \ref{wallregu}. Cycles are green, axes are red or blue. On the right, one extremity of the essential part of the axis is of the same parity as the extremities of $S_4$. On the left, the parities of both extremities of the essential part of the axis are different from the parity of the extremities of $S_4$.}
    \end{it}
    \end{minipage}
\end{center}
\end{sloppypar}

\subsubsection{Special wall pairs}

\label{wallspecial}

\begin{sloppypar}
    Let $k$ be an even integer such that $2 \leq k \leq m-1$. The IV triangulation $\tau$ of $\Delta_m^3$ contains a union of simplices that form a cone with basis $T_k$ and with vertex an integer point lying on $T_{k-1}$. The triangulation also contains a cone with basis $T_k$ and with vertex an integer point lying on $T_{k+1}$. Consider such a cone and denote by $v$ its vertex. In the previous section \ref{wallregu}, we chose an integer point $e$ lying in the relative interior of $S_k$ and at lattice distance $1$ from one of the extremities of $S_k$. For each point $f$ lying in the interior of $S_k$ and different from $e$, we constructed an axis starting with the segment $[f,v]$.\\
    We apply Lemma \ref{linsystem} as in the previous section \ref{wallregu} to a $3$-simplex $Q$ lying in the intersection of the star of $[v,e]$ with the cone over $T_k$ and denote by $O$ the solution orthant. We construct a $1$-axis by taking the union of the following segments:
    \begin{itemize}
        \item the $O$-copy of $[e,v]$ (the essential part of the axis);
        \item the $O$-copy of the segment $[e,v']$, where $v'$ is the vertex of the other cone over $T_k$;
        \item the segments whose extremities are the $O$-copy of $v$ (respectively the $O$-copy of $v'$) and the $O$- or $O'$-copy of a triangulation vertex $a'$ lying in the boundary on $T_k$ but not in $S_k$, where $O'$ is an orthant adjacent to $O$ and such that the suspension over the $O$- and $O'$-copies of $T_k$ share the same vertices. We choose $a'$ and the orthant ($O$ or $O'$) such that the chosen copy of $a'$ bears the same sign as the $O$-copy of $a$.
    \end{itemize}
    The described $1$-axis is contained in the union of two adjacent copies of $\Delta_m^3$ (without identifying pairs of antipodal points lying on the boundary of the union of the symmetric copies of $\Delta_m^3$). The description of the dual cycle is the same as in Section \ref{wallregu}.\\
    Now suppose that $m$ is even and let $k=m$. Take the union of the two copies of $[e,v]$ whose extremities are the $O$-copies of $e$ and $v$ (after identifying antipodal points of the boundary of $(\Delta_m^3)^*$). This forms a $1$-axis. It is not contained in a single orthant or in a union of two adjacent orthants. The description of a dual cycle is the same as in the previous section. Notice that it is possible that the cycle is not contained in a union of two adjacent orthants. This default is responsible for at most one of the possible exceptions mentioned in Section \ref{properties}\\
    For every even $k$, we obtain one $(1,1)$-pair, called a \textit{special wall pair}. An example of such a pair is represented in Figure 8.
    \begin{center}
    \includegraphics[scale=0.7]{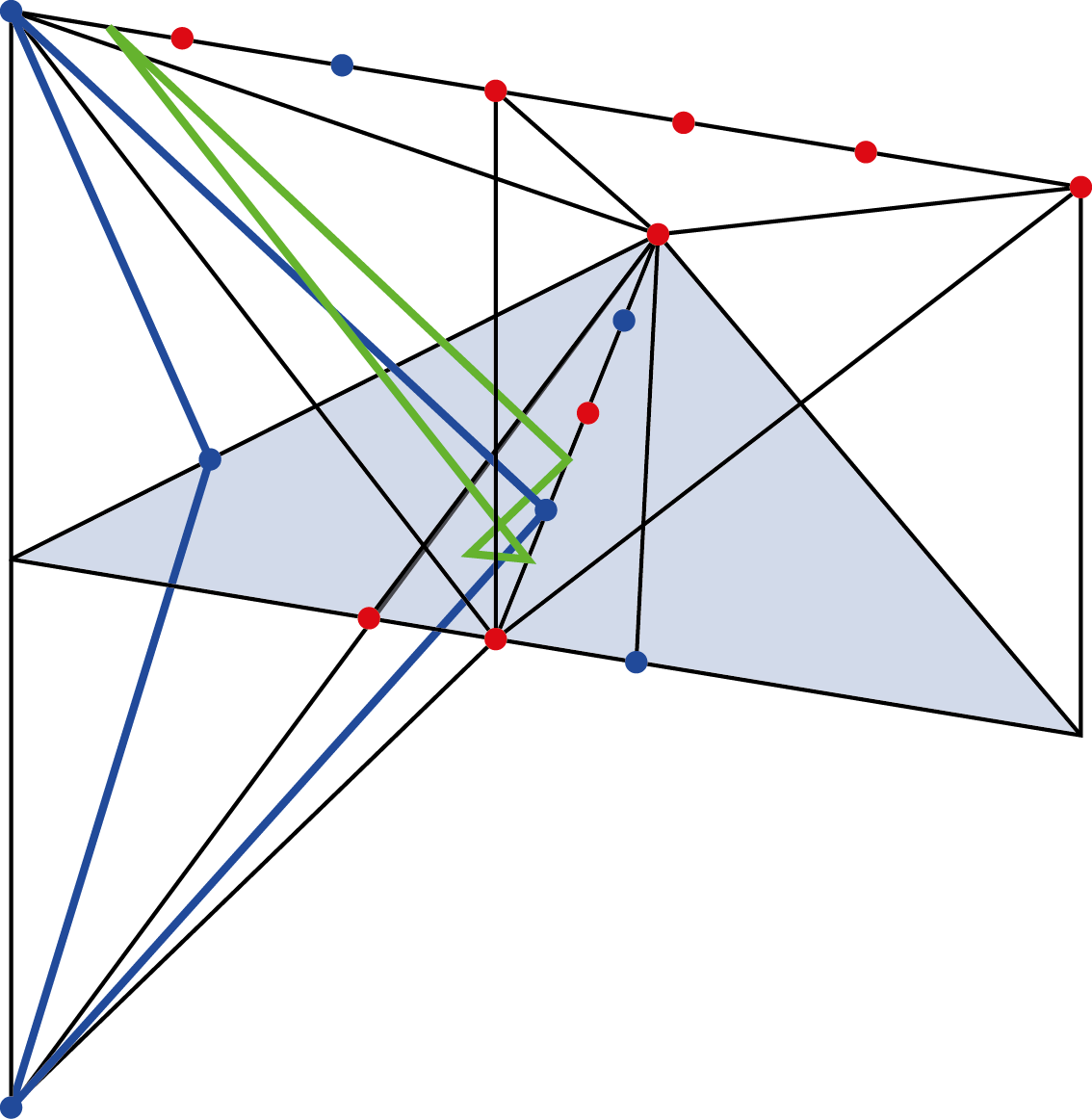}
    \begin{minipage}{16cm}
    \begin{it}
    \small{Fig. 8: Two adjacent copies of a cone over $T_4$ and of a join of $S_4$ and $S_3$ are represented, as well as the copy of another cone over $T_4$. The colors of the integer points (blue or red) represent their sign. A typical special wall pair cycle appears in green. The dual axis appears in blue. In this case, the sign distribution on $T_4$ is not a Harnack one, hence the axis can be chosen to be contained in a single orthant.}
    \end{it}
    \end{minipage}
\end{center}
    In Sections \ref{suspspecial} \ref{wallregu} \ref{wallspecial}, we have exhibited two $(1,1)$-pairs for each integer point in the relative interior of $S_k$ with $1\leq k \leq m-1$ and one $(1,1)$-pair for each integer point in the relative interior of $S_m$. In total, this amounts to $m^2-2m+1$ pairs of type $(1,1)$.
\end{sloppypar}

\subsection{Defining the order on the collection of $(j, n-j)$-pairs}

\subsubsection{Linking relations}

\label{link}

\begin{sloppypar}
    Every described cycle is linked with its axis. In the case of $(1,1)$-pairs, there may be other linking relations. When considering two different $(1,1)$-pairs $(c,a)$ and $(c',a')$ such that $c$ (respectively $c'$) lies in the union of the stars of the segments in the essential part of $a$ (respectively $a'$), their associated linking numbers can be obtained by checking whether or not one of the axes contains the essential part of the other. There are also relations involving the suspension pairs whose cycle does not lie in a union of stars of  segments of the dual axis. We describe below all possible linking relations involving two $(1,1)$-pairs. Whether the corresponding  $\mathbb{Z}_2$-linking numbers actually are equal to $1$ or not depends on the chosen sign distribution.
    \begin{itemize}
        \item The axes of some join pairs contain a segment contained in the copy of a $2$-dimensional wall. This segment may appear in the axis of a regular or special wall pair. Hence these join axes may be linked with some regular or special wall cycles.
        \item The axes of some wall pairs contain a segment that also belongs to the  axis of a special wall pair. These axes may be linked with a special wall cycle.
        \item The axis of any regular wall pair contains two segments lying in the symmetric copy of some triangle $T_k$. If $k \neq m$, one of these segments may appear in the axis of the special suspension pair associated with this copy of $T_k$. Hence it may be linked with the special suspension cycle associated to the suspension over $T_k$.
        \item The cycles from some suspension pairs (those described in \ref{11d}) are linked with the axis of a special suspension pair.
    \end{itemize}
\end{sloppypar}

\subsubsection{Partial order on $\mathcal{C}_1$}

\label{order}

\begin{sloppypar}
    On $\mathcal{C}_1$, we take the partial order generated by the following relations.
    \begin{itemize}
        \item Let $(c,a)$ be a join pair. Let $(c',a')$ be a regular wall pair. Then we have the relation $(c,a)<(c',a')$.
        \item Let $(c,a)$ be a regular wall pair. Let $(c',a')$ be a special wall pair. Then we have the relation $(c,a)<(c',a')$.
        \item Let $(c,a)$ be a special wall pair. Let $(c',a')$ be a special suspension pair. Then we have the relation $(c,a)<(c',a')$.
        \item Let $(c,a)$ be a special suspension pair. Let $(c',a')$ be a regular suspension pair. Then we have the relation $(c,a)<(c',a')$.
    \end{itemize}
\end{sloppypar}

\subsection{Independence of the classes represented by the described cycles}

\label{indep}

\begin{sloppypar}
    The classes represented by the cycles from the pairs in $\mathcal{C}_2$ are clearly independent, as the $2$-cycles represent the classes of distinct connected components of the surface $\widetilde{\Gamma}_m^{IV}$.
    Let \mbox{$i:\mathcal{C}_1 \rightarrow \{1,...,\# \mathcal{C}_1\}$} be an injective map inducing a total order on $\mathcal{C}_1$ that refines the partial order defined in Section \ref{order}. Here $\# \mathcal{C}_1$ denotes the cardinal of the set $\mathcal{C}_1$. We want to check whether the dimension of the subspace of $H_1(\widetilde{\Gamma}_m^{IV};\mathbb{Z}_2)$ spanned by the classes represented by the $1$-cycles from the pairs of $\mathcal{C}_1$ is equal to the number of described cycles. We define the matrix $\mathcal{A}_1$ whose coefficient in the $i$-th row and $j$-th column is the $\mathbb{Z}_2$-linking number of the cycle from the $i$-th pair with the axis from the $j$-th pair.
    The order defined Section \ref{order} has been chosen by taking into account the linking relations presented in Section \ref{link} to ensure the following property: if $(c,a)$ and $(c',a')$ are two distinct pairs in $\mathcal{C}_1$ such that $(c,a)<(c',a')$, then the $\mathbb{Z}_2$-linking  number of $c$ and $a'$ is $0$. Thus the matrix $\mathcal{A}_1$ is lower triangular. Hence, using \mbox{Proposition \ref{indepmatrice}}, we obtain that the classes represented by the cycles from the pairs in $\mathcal{C}_1$ form a basis of the subspace of $H_1(\Tilde{\Gamma}^{IV}_m)$ they span. In particular, these classes are independent.
\end{sloppypar}

\subsection{Proof of Proposition \ref{resultmax}}

\label{maxsurf}
\label{cyclessupp}

\begin{sloppypar}
    \textbf{Proof:} We have exhibited $\frac{(m-1)(m-2)(m-3)}{6}$ independent classes in $H_2(\Tilde{\Gamma}^{IV}_m; \mathbb{Z}_2)$. There is at least one more $2$-cycle that represents a homology class independent from those represented by a cycle from a pair in $\mathcal{C}_2$. Indeed, there are hyperplane pieces from the piecewise linear surface that do not belong to any $2$-cycle from a pair in $\mathcal{C}_2$. Hence, we have $b_0(\Tilde{\Gamma}^{IV}_m) = b_2(\Tilde{\Gamma}^{IV}_m) \geq \frac{(m-1)(m-2)(m-3)}{6}+1$.
    We have exhibited $\frac{2m^3}{3}-2m^2+\frac{7m}{3}-1$ independent classes in $H_1(\Tilde{\Gamma}^{IV}_m; \mathbb{Z}_2)$. Hence, we have $b_1(\Tilde{\Gamma}^{IV}_m) \geq$ $\frac{2m^3}{3}-2m^2+\frac{7m}{3}-1$. Hence, $b_0(\Tilde{\Gamma}^{IV}_m) + b_1(\Tilde{\Gamma}^{IV}_m) + b_2(\Tilde{\Gamma}^{IV}_m) \geq 2 \times (\frac{(m-1)(m-2)(m-3)}{6}+1) + \frac{2m^3}{3}-2m^2+\frac{7m}{3}-1 = m^3-4m^2+6m-1$. It remains to use Smith's congruence to obtain that the surface $X_m^{IV}$ is maximal and that \mbox{$b_1(\mathbb{R}X_m^{IV}) =$ $\frac{2m^3}{3}-2m^2+\frac{7m}{3}$}.
    Now, \mbox{$\chi(\mathbb{R}X_m^{IV}) = 2 \times (\frac{(m-1)(m-2)(m-3)}{6}+1) - (\frac{2m^3}{3}-2m^2+\frac{7m}{3}) = \frac{-m^3}{3} + \frac{4m}{3}$}. Hence, $X_m^{IV}$ is of type $\chi = \sigma$. \hfill \qedsymbol\\
\end{sloppypar}

\begin{sloppypar}
    Notice that all the $1$-cycles we described in Sections \ref{11a}, \ref{11b}, \ref{11c}, \ref{11d}, \ref{suspspecial}, \ref{join}, \ref{wallregu} and \ref{wallspecial} realize the trivial class in $H_1(\widetilde{(\Delta_m^3)^*};\mathbb{Z}_2)$. As $X_m^{IV}$ is maximal, the inclusion morphism $in: \mathbb{R}X_m^{IV} \rightarrow \mathbb{RP}^3$ induces a surjective morphism $in^*:H_1(\mathbb{R}X_m^{IV};\mathbb{Z}_2) \rightarrow H_1(\mathbb{RP}^3;\mathbb{Z}_2)$ (this fact is due to V. Kharlamov, see e.g. \cite{kharlamov}). Hence, there exists a $1$-cycle $c_{NT}$ of $\Tilde{\Gamma}^{IV}_m$ that realizes the non-trivial homology class in $H_1(\widetilde{(\Delta_m^3)^*};\mathbb{Z}_2)$. The properties of this cycle are used in \mbox{Section \ref{cyclessupp3var}}.
\end{sloppypar}

\section{The Itenberg-Viro construction for $3$-manifolds}

\label{construction3}

\begin{sloppypar}
In this section, we describe two variants of a primitive triangulation and a sign distribution that produce non-singular maximal real algebraic hypersurfaces in $\mathbb{RP}^4$ through combinatorial patchworking. This construction is due to Itenberg and Viro (see \cite{fakeitenbergviro}). In the following, $m$ is a positive even integer and $x_1, x_2, x_3, x_4$ denote the standard coordinates in $\mathbb{R}^4$. 
%We denote by $\Delta_m^4$ the $4$-simplex with vertices $(0,0,0,0)$, $(m,0,0,0)$, $(0,m,0,0)$, $(0,0,m,0)$, $(0,0,0,m)$.
\end{sloppypar}

\subsection{Triangulations and sign distributions}

\begin{sloppypar}
    First, recall the standard triangulation of the product of two triangles. Let $T_1$ and $T_2$ be non-degenerate triangles in two copies of the plane with vertices labelled from $0$ to $2$. The cartesian product $T_1 \times T_2$ is a four-dimensional polytope with $9$ vertices and $12$ facets which are all prisms. The labelling of the vertices of $T_1$ and $T_2$ induces a two-digit labelling of the vertices of $T_1 \times T_2$. The standard triangulation of $T_1 \times T_2$ with respect to the choice of labelling of their vertices consists of the following six $4$-simplices, represented by the \mbox{$5$-tuples} of their vertices: $(00,01,02,12,22)$,  $(00,01,11,12,22)$, $(00,01,11,21,22)$, $(00,10,11,12,22)$, $(00,10,11,21,22)$, $(00,10,20,21,22)$. There are six $3$-simplices of this triangulation that are common faces of two of the above $4$-simplices : $(00,01,12,22)$, $(00,01,11,22)$, $(00,11,12,22)$, $(00,11,21,22)$, $(00,10,11,22)$, $(00,10,21,22)$.\\
    We describe successive subdivisions of $\Delta_m^4$ that produce what we call the \emph{odd} (respectively, the \emph{even}) \emph{Itenberg-Viro triangulation} (or IV triangulation) of $\Delta_m^4$.
    \begin{itemize}
        \item Let $k$ be an integer such that $1\leq k \leq m$. Let $R_k := \Delta_m^4 \bigcap \{x_4=m-k\}$. Subdivide $\Delta_m^4$ using the tetrahedra $R_k$. For every integer $k$ such that $1\leq k \leq m$, let $T_k$ be the $2$-face of $R_k$ that lies in the hyperplane $\{x_1+x_2+x_3+x_4=m\}$.
        \item Let $k$ be an odd (respectively, even) integer such that $1 \leq k \leq m-1$. Take the cones over $R_k$ with vertices $(0,0,0,m-k+1)$ and $(0,0,0,m-k-1)$. In the case of the \emph{even} IV construction, take the cone over $R_m$ with vertex $(0,0,0,1)$.
        \item Let $k$ be an even (respectively, odd) integer such that $1 \leq k \leq m$. For every integer $l$ such that $1 \leq l \leq k$, let $T_k^l$ be the triangle parallel to $T_k$, whose vertices are integer points lying on the edges of $R_k$ and whose edges have lattice length $l$. Subdivide $R_k$ using these triangles. Let $T_k^0$ be the vertex of $R_k$ opposite to $T_k^k = T_k$.
        \item Let $k$ be an odd (respectively, even) integer such that $1 \leq k \leq m$ (respectively, such that $1 \leq k \leq m$). Let $l$ be an integer such that $0 \leq l \leq k-2$. We label the vertices of $T_k$, $T_{k-1}^l$ and $T_{k-1}^{l+1}$ as follows.
        \begin{itemize}
            \item The triangle $T_k$ is labeled by $1$.
            \item The triangle $T_{k-1}^l$ is labeled by $0$ if $l$ is odd (respectively, even) and by $2$ otherwise.
            \item The triangle $T_{k-1}^{l+1}$ is labeled by $0$ if $l$ is even (respectively, odd) and by $2$ otherwise.
            \item A vertex of any of the three triangles is labeled by two numbers. The first one is the is the label of the triangle the vertex belongs to. The second one is $0$ (respecively, $1$, $2$) if the vertex lies in the plane $\{ x_2=0 \} \bigcap \{ x_3=0 \}$ (respectively, in the plane $\{ x_1=0 \} \bigcap \{ x_3=0 \}$, in the plane $\{ x_1=0 \} \bigcap \{ x_2=0 \}$).
        \end{itemize}
        The convex hull $H_{k,-}^l$ of the triangles $T_k$, $T_{k-1}^l$ and $T_{k-1}^{l+1}$ has the same combinatorial structure as a product of two triangles. Using the labeling of the vertices of $H_{k,-}^l$ described above, subdivide $H_{k,-}^l$ using the triangulation of a product of two triangles presented at the beginning of the paragraph. Repeat this step with the convex hull $H_{k,+}^l$ of the triangles $T_k$, $T_{k+1}^l$ and $T_{k+1}^{l+1}$, for every odd (respectively, even) integer $k$ such that $1\leq k \leq m-1$ and for every integer $l$ such that $0\leq l \leq k$. We denote by $\tau'$ the triangulation of $\Delta_m^4$ obtained after this subdivision step.
        \item Let $k$ be an even (respectively, odd) integer such that $1 \leq k \leq m$. Complete the subdivision of the tetrahedron $R_k$ into an Itenberg-Viro triangulation of $R_k$ (seeing it as a copy of $\Delta_k^3$). Take the induced subdivision of $\Delta_m^4$.
        \item Let $k$ be an odd (respectively, even) integer such that $1 \leq k \leq m$. In the case of the odd IV construction, triangulate $R_k$, seen as a copy of $\Delta_k^3$, using an Itenberg-Viro triangulation that satisfies the following conditions.
        \begin{itemize}
            \item The $2$-face chosen at the first step of the subdivision in Section \ref{triangulation} is $R_k \bigcap \{ x_3=0 \}$.
            \item The segment $E_1$ chosen at the second step of the subdivsion in Section \ref{triangulation} is \mbox{$R_k \bigcap  \{ x_1=0, x_1+x_2+x_3=k \}$}.
            \item The segment $E_2$ chosen at the second step of the subdivsion in Section \ref{triangulation} is \mbox{$R_k \bigcap  \{ x_2=0, x_1+x_2+x_3=k \}$}.
        \end{itemize}
        In the case of the even IV construction, triangulate $R_k$ using an Itenberg-Viro triangulation that satisfies the following conditions.
        \begin{itemize}
            \item The $2$-face chosen at the first step of the subdivision in Section \ref{triangulation} is \mbox{$R_k \bigcap \{ x_1+x_2+x_3=k \}$}.
            \item The segment $E_1$ chosen at the second step of the subdivsion in Section \ref{triangulation} is \mbox{$R_k \bigcap  \{ x_2=0,x_3=0 \}$}.
            \item The segment $E_2$ chosen at the second step of the subdivsion in Section \ref{triangulation} is \mbox{$R_k \bigcap  \{ x_1=0, x_2=0 \}$}.
        \end{itemize}
        Take the induced subdivision of $\Delta_m^4$.
        %\item Take the primitive triangulation of $\Delta_m^4$ induced by the subdivision described above.
    \end{itemize}
    Notice that by construction, the odd and even IV triangulations of $\Delta_m^4$ are primitive. The vertices of the odd IV triangulation are endowed with the sign $+$ if they are of parity $(1,0,0,1)$ or $(0,0,1,1)$; otherwise they are endowed with the sign $-$. Likewise, the vertices of the even IV triangulation are endowed with the sign $+$ if they are of parity $(0,1,1,0)$, $(1,1,0,0)$ or $(0,1,0,0)$; otherwise, they are endowed with the sign $-$.
\end{sloppypar}

\subsection{Convexity of the even and odd Itenberg-Viro triangulations of $\Delta_m^4$}

\begin{sloppypar}
To be able to apply Viro's patchworking theorem to $\Delta_m^4$ endowed with the odd or even IV datum, we need to verify that the odd and even IV triangulations of $\Delta_m^4$ are convex.
\begin{prop}
    There exists a convex piecewise linear function $\nu:\Delta_m^4 \rightarrow \mathbb{R}$ whose domains of linearity coincide with the simplices of the odd (respectively, of the even) IV triangulation of $\Delta_m^4$.
\end{prop}
\textbf{Proof:} We describe the required function $\nu:\Delta_m^4 \rightarrow \mathbb{R}$ in the case of the odd IV construction. The case of the even IV triangulation can be dealt with in a similar way, inverting the role of the even and odd objects. We start with a simple piecewise linear convex function and we refine it in order to follow the different subdivision steps described above.
\begin{itemize}
    \item Denote by $\nu_1:\Delta_m^4 \rightarrow \mathbb{R}$ the pointwise maximal piecewise linear convex function whose value on every integer point $(x_1,x_2,x_3,x_4)$ is $e^{x_4}$.
    \item Let $\bar{\mu}_1(0,0,0,k):= 1$ for any even integer $k$ such that $0 \leq k \leq m$. Let $\bar{\mu}_1(v):= 0$ for any integer point $v$ that is a vertex of $R_k$ for some integer $k$ such that $0\leq k \leq m$ but that is not of the form $(0,0,0,k)$ for some even integer $k$.
    Let $\mu_1:\Delta_m^4 \rightarrow \mathbb{R}$ be a piecewise linear function such that, for any integer $k$ such that $1 \leq k \leq m$, the graph of the restriction of $\mu_1$ to $\{ m-k \leq x_4 \leq m-k+1 \}$ is the lower part of the convex hull of the points from the set $\{ (v, \bar{\mu}_1(v)) \in \mathbb{R}^5 |$ $v$ is a vertex of $R_{k}$ or $R_{k-1} \}$. For a sufficiently small positive real number $\epsilon_1$, the piecewise linear function $\nu_2 = \nu_1 + \epsilon_1 \mu_1$ is convex.
    \item Let $\bar{\mu}_2(x_1,x_2,x_3,x_4):=e^{x_1+x_2+x_3}$ for any integer point $(x_1,x_2,x_3,x_4)$ with $x_4$ even. Let $\bar{\mu}_2(x_1,x_2,x_3,x_4):=0$ for any integer point $(x_1,x_2,x_3,x_4)$ with $x_4$ odd. Let $\mu_2:\Delta_m^4 \rightarrow \mathbb{R}$ be a piecewise linear function such that for any integer $k$ such that $1 \leq k \leq m$, the graph of the restriction of $\mu_2$ to $\{ m-k \leq x_4 \leq m-k+1 \}$ is the lower part of the convex hull of the points from the set $\{ (v, \bar{\mu}_2(v)) \in \mathbb{R}^5|v$ is an integer point lying in $R_k$ or $R_{k-1} \}$. For a sufficiently small positive real number $\epsilon_2$, the piecewise linear function $\nu_3 = \nu_2 + \epsilon_2 \mu_2$ is convex.
    \item Let $v=(x_1,x_2,x_3,x_4)$ be an integer point lying in $\Delta_m^4$. Let $E$, $F$, $G$ be three positive real numbers such that $E>F>G$. We let
    \begin{itemize}
        \item $\bar{\mu}_3(v)= E(2x_1+x_2)$ if $x_4 \equiv 0$(mod $2$) and $x_1+x_2+x_3 \equiv 1$(mod $2$);
        \item $\bar{\mu}_3(v)=G(2x_1 + x_2)$ if $x_4 \equiv 0$(mod $2$) and $x_1+x_2+x_3 \equiv 0$(mod $2$);
        \item $\bar{\mu}_3(v)=F(2x_1+x_2)$ if $x_4 \equiv 1$(mod $2$) and $x_1+x_2+x_3+x_4 = m$;
        \item $\bar{\mu}_3(v)=0$ if $x_4 \equiv 1$(mod $2$) and $x_1=x_2=x_3=0$.
    \end{itemize}
    Let $\mu_3: \Delta_m^4 \rightarrow \mathbb{R}$ be a piecewise linear function such that for any odd integer $k$ such that $1\leq k \leq m$, we have that
    \begin{itemize}
        \item for any integer $l$ such that $0 \leq l \leq k-2$ (respectively, such that $0 \leq l \leq k$), the graph of the restriction of $\mu_3$ to the convex hull $H_{k,+}^l$ of the triangles $T_k$, $T_{k+1}^l$ and $T_{k+1}^{l+1}$ (respectively, to  the convex hull $H_{k,-}^l$ of the triangles $T_k$, $T_{k-1}^l$ and $T_{k-1}^{l+1}$) is the lower part of the convex hull of the points from the set $\{ (v, \bar{\mu}_3(v)) \in \mathbb{R}^5|$ $v$ is an integer point lying in $H_{k,+}^l \}$ (respectively, the set $\{ (v, \bar{\mu}_3(v)) \in \mathbb{R}^5|$ $v$ is an integer point lying in $H_{k,-}^l \}$);
        \item the graph of the restriction of $\mu_3$ to the suspension over $R_k$ with vertices \mbox{$(0,0,0,m-k+1)$} and \mbox{$(0,0,0,m-k-1)$} is the lower part of the convex hull of the points from the set 
        $\{ (v,\bar{\mu}_3(v)) \in \mathbb{R}^5|$ $v\in R_k$ or $v = (0,0,0,m-k+1)$ or $v=(0,0,0,m-k-1) \}$.
    \end{itemize}
    For a sufficiently small real number $\epsilon_3$, the piecewise linear function $\nu_4 = \nu_3 + \epsilon_3 \mu_3$ is convex. The linearity domains of  $\nu_4$ are the simplices of the triangulation $\tau'$.
    \item For any odd integer $k$ such that $1 \leq k \leq m-1$, there exists a convex piecewise linear function $\nu^k: R_k \rightarrow \mathbb{R}$ that certifies the convexity of the odd IV triangulation restricted to $R_k$ (see \cite{tsurf}).
    For any even integer $k$ such that $0 \leq k \leq m$, for any integer $l$ such that \mbox{$1 \leq l \leq k$}, there exists a convex piecewise linear function $\nu^k_l$ that certifies the convexity of the odd IV triangulation restricted to $T_k^l$.
    Any element $\mu$ of one of the families $(\nu^k)_{1\leq k \leq m, \: k \; odd}$ and $(\nu^k_l)_{1\leq k \leq m, \: k \; even, \: 1\leq l \leq k}$ can be extended to a piecewise linear function $\widetilde{\mu}: \Delta_m^4 \to \mathbb{R}$ such that
    \begin{itemize}
        \item the value of $\widetilde{\mu}$ on any integer point $v$ that does not belong to the domain of $\mu$ is zero;
        \item the vertices of the graph of $\widetilde{\mu}$ are integer points;
        \item any linearity domain of $\widetilde{\mu}$ is a subset of a linearity domain of $\nu_4$;
        \item the restriction of $\widetilde{\mu}$ to any linearity domain of $\nu_4$ is convex.
    \end{itemize}
    %Let $\widetilde{\nu}^k(x_1,x_2,x_3,x_4) := \nu_k(x_1,x_2,x_3)$ for any integer point $(x_1,x_2,x_3,x_4)$ lying in $R_k$. Let $\widetilde{\nu}^k(x_1,x_2,x_3,x_4) := 0$ for any other integer point lying in $\Delta_m^4$. Let $\widetilde{\nu}^k: \Delta_m^4 \rightarrow \mathbb{R}$ be a piecewise linear function such that for any integer $k$ such that $1 \leq k \leq m$, the graph of the restriction of $\widetilde{\nu}^k$ to $\{ m-k \leq x_4 \leq m-k+1 \}$ is the lower part of the convex hull of the points from the set $\{ (v, \widetilde{\nu}^k(v)) \in \mathbb{R}^5|v$ is an integer point lying in $R_k$ or $R_{k-1} \}$. 
    To perform the last refinement, we add to $\nu_4$ a linear combination of elements of the families $(\widetilde{\nu}^k)_{1\leq k \leq m, \: k \; odd}$ and $(\widetilde{\nu}^k_l)_{1\leq k \leq m, \: k \; even, \: 1\leq l \leq k}$ multiplied by a sufficiently small positive real number $\epsilon_4$. The resulting function $\nu: \Delta_m^4 \rightarrow \mathbb{R}$ is convex and its linearity domains are the simplices of the odd IV triangulation of $\Delta_m^4$. \hfill \qedsymbol
\end{itemize}
\begin{rmk}
\label{rmkconvexity}
The last step of the construction of the convex function is not specific to the functions described in \cite{tsurf}. Any family of convex piecewise linear functions $(\nu^k)_{1\leq k \leq m, \: k \; odd}$ with $\nu^{k}: R_k \rightarrow \mathbb{R}$ can be used to perform the final refinement. This remark also holds in the case of the even IV construction.
\end{rmk}
\end{sloppypar}
\begin{sloppypar}
The hypersurfaces obtained by applying the combinatorial patchworking to the odd and even IV data are interesting in the sense of the following propositions.
\begin{prop}
    \label{maximal}
    The real algebraic hypersurfaces of degree $m$ in $\mathbb{RP}^4$ produced by applying the combinatorial patchworking theorem to the signed odd IV triangulation and to the signed even IV triangulation of $\Delta_m^4$ are maximal.
    \end{prop}
    \begin{prop}
    \label{type}
    The real algebraic hypersurfaces of even degree $m$ in $\mathbb{RP}^4$ produced by applying the combinatorial patchworking theorem to the signed odd IV triangulation and to the signed even IV triangulation are of type $\chi = \sigma$.
    \end{prop}
The goal of Sections \ref{maximality3manif} and \ref{prooftype} is respectively to prove Proposition \ref{maximal} and Proposition \ref{type}. 
\end{sloppypar}

\section{Maximality of the Itenberg-Viro $3$-manifolds}

\label{maximality3manif}

\begin{sloppypar}
    From now on, $m$ is an even integer greater than $2$. We describe cycle-axis pairs for piecewise linear $3$-manifolds produced by applying the combinatorial patchworking to the data described in Section \ref{construction3}. We describe a collection $\mathcal{C}_3$ of $3$-cycles, a collection $\mathcal{C}_2$ of $(2,1)$-pairs and a collection $\mathcal{C}_1$ of $(1,2)$-pairs. We then endow $\mathcal{C}_1$ and $\mathcal{C}_2$ with a total order. To prove the maximality of the patchworked hypersurfaces, we need the final collection of pairs to satisfy the following properties.
    \begin{itemize}
        \item Let $k=1,2$. Let $(c,a) \in \mathcal{C}_k$ be a $(k,3-k)$-pair. Then, the $\mathbb{Z}_2$-linking number of $c$ and $a$ is 1. 
        \item Let $k=1,2$. Let $(c,a)$ and $(c',a') \in \mathcal{C}_k$ be two distinct $(k,3-k)$-pairs such that $(c,a)<(c',a')$. Then, the $\mathbb{Z}_2$-linking number of $c$ and $a'$ is $0$.
        \item Let $(c,a)\in \mathcal{C}_2$. Let $(c',a')\in \mathcal{C}_1$. Then, the intersection of the supports of $c$ and $c'$ is empty.
    \end{itemize}
    %We use the first two conditions in Section \ref{independence3} to prove that the homology classes represented by the $2$-cycles from pairs in $\mathcal{C}_2$ (repectively the $1$-cycles from pairs in $\mathcal{C}_1$) are linearly independent over $\mathbb{Z}_2$. The last condition ensures that a class represented by a cycle from a pair in $\mathcal{C}_1$ is not Poincaré dual to any class represented by a cycle from a pair in $\mathcal{C}_2$. 
    We describe cycle-axis pairs simultaneously for the odd and even IV construction of degree $m$. The piecewise linear $3$-manifold under consideration is denoted by $\Tilde{\Gamma}_m^{IV}$.
\end{sloppypar}

\subsection{Cycle-axis pairs in the suspension over a tetrahedron $R_k$}

\label{susprk}

\begin{sloppypar}
    Let $k$ be an integer such that $1 \leq k \leq m-1$. If $k$ is odd, the odd IV triangulation of $\Delta_m^4$ contains a union of simplices forming a suspension over $R_k$. If $k$ is even, the even IV triangulation contains a union of simplices forming a suspension over $R_k$. In this section, we suppose that either $k$ is odd and $\Delta_m^4$ is endowed with the odd IV datum or that $k$ is even and $\Delta_m^4$ is endowed with the even IV datum.
\end{sloppypar}

\subsubsection{Pairs of type $(3,0)$}

\label{susp30}

\begin{sloppypar}
    Let $v$ be an integer point lying in the relative interior of $R_k$. The star of $v$ contains $8$ integer points distinct from $v$ and of $4$ different parities. Let $Q$ be a $4$-simplex lying in the star of $v$. We apply Lemma \ref{linsystem} to $Q$: up to inversion of all signs, there exists a unique orthant $O$ such that the $O$-copy of $v$ bears the sign $+$ while the $O$-copies of the other vertices of $Q$ bear the sign $-$. Hence, the star of the $O$-copy of $v$ contains a cell complex homeomorphic to a $3$-sphere and belonging to $\Tilde{\Gamma}_m^{IV}$. This sphere is a $3$-cycle. We do not need to describe dual axes, as different connected components clearly represent independent classes in $H_3(\Tilde{\Gamma}_m^{IV};\mathbb{Z}_2)$.
\end{sloppypar}

\subsubsection{Pairs of type $(2,1)$ - a}

\label{susp21a}

\begin{sloppypar}
    Let $v$ be an integer point lying in the relative interior of $R_k$. The star of $v$ contains $8$ integer points distinct from $v$ and of $4$ different parities. Let $Q$ be a $4$-simplex lying in the star of $v$ with vertices $v$, $v_1$, $v_2$, $v_3$, $v_4$, where $v_4$ is the vertex of $Q$ that does not lie in $R_k$. We apply Lemma \ref{linsystem} to $Q$: up to inversion of all signs, there exists a unique orthant $O$ such that the $O$-copies of $v$ and $v_4$ bear the sign $+$ while the $O$-copies of $v_1$, $v_2$ and $v_3$ bear the sign $-$. As described in Section \ref{susp20}, there is a $2$-cycle $c$ of $\Tilde{\Gamma}_m^{IV} \bigcap \Tilde{(R_k)^*}$ homeomorphic to a $2$-sphere in the intersection of the $O$-copy of $R_k$ with the star of the $O$-copy of $v$. The cycle $c$ comes with an axis $a$ which is a pair of points which bear the same sign as the $O$-copy of $v$. Using the $(2,0)$-pair $(c,a)$, we form a $(2,1)$-pair $(c,\hat{a})$ by taking $\hat{a}$ to be the suspension over $a$ whose vertices are the $O$-copies of $v_4$ and the other vertex of the star of $v$ that does not lie in $R_k$. The $2$-cycle $c$ bounds a $3$-disk in the $O$-copy of $R_k$. This disk intersects $\hat{a}$ in exactly one point: the $O$-copy of $v$. Hence, the $\mathbb{Z_2}$-linking number of $\hat{a}$ and $c$ is equal to $1$.
\end{sloppypar}

\subsubsection{Pairs of type $(2,1)$ - b}

\label{susp21b}

\begin{sloppypar}
    Let $v$ be an integer point lying in the relative interior of $R_k$. The star of $v$ contains $8$ integer points distinct from $v$ and of $4$ different parities. Let $Q$ be a $4$-simplex lying in the star of $v$ with vertices $v$, $v_1$, $v_2$, $v_3$, $v_4$, where $v_4$ is the vertex of $Q$ that does not lie in $R_k$ and $v_3$ is the vertex of $Q$ that does not lie in the triangle that is parallel to a face of $R_k$, that contains $v$ and that was used to define the IV triangulation on $R_k$. We apply Lemma \ref{linsystem} to $Q$: up to inversion of all signs, there exists a unique orthant $O$ such that the $O$-copies of $v$ and $v_3$ bear the sign $+$ while the $O$-copies of $v_1$, $v_2$ and $v_4$ bear the sign $-$. As described in \ref{11a}, there is a $1$-cycle of $\Tilde{\Gamma}_m^{IV} \bigcap \Tilde{(R_k)*}$ in the $O$-copy of the intersection of the star of $v$ with $R_k$ linked with a $1$-axis $a$ lying in the $O$-copy of $R_k$ and containing the $O$-copy of the segment $[v,v_3]$. Using the $(1,1)$-pair $(c,a)$, we form a $(2,1)$-pair $(\hat{c},a)$ by taking $\hat{c}$ to be the suspension over $c$ with vertices the $O$-copies of the middle points of the segments $[v,v_4]$ and $[v,v'_4]$, where $v_4$ and $v'_4$ are the points in the star of $v$ that do not lie on $R_k$. The linking properties of $c$ and $a$ ensure that the $\mathbb{Z}_2$-linking number of $\hat{c}$ and $a$ is $1$.
    The $(2,1)$-pairs described in this Section are obtained by taking a suspension over a cycle from a pair described in Section \ref{11a} and contained in a symmetric copy of $R_k$. The same can be done starting with a pair $(c,a)$ described in Sections \ref{11b}, \ref{11c} or \ref{11d} to form other $(2,1)$-pairs.
\end{sloppypar}

\subsubsection{Pairs of type $(2,1)$ - c}

\label{susp21c}

\begin{sloppypar}
    Let $A$ be a triangulation segment joining two vertices $p$ and $p'$ lying in the relative interior of two different $2$-faces of $R_k$. The vertices $p$ and $p'$ necessarily have different parities. The segment $A$ lies in the join of two segments $S$ and $S'$ contained in two different $2$-faces of $R_k$ and that are unions of triangulation edges. Suppose that $p$ lies in $S$ and that $p'$ lies in $S'$. Aside from $p$ and $p'$, the star of $A$ contains $6$ integer vertices of $3$ different parities: the extremities $v_2$ and $v_2'$ (respectively $v_3$ and $v_3$') of the triangulation segments adjacent to $p$ (respectively $p'$) and lying in $S$ (respectively $S'$) and the vertices $v_4$ and $v_4'$ of the suspension over $R_k$. Let $Q$ be any $4$-simplex lying in the star of $A$. We apply Lemma \ref{linsystem} to $Q$: up to inversion of all signs, there exists a unique orthant $O$ such that the $O$-copies of $p$ and $p'$ bear the $+$ while the $O$-copies of the other vertices bear the sign $-$. As described in \ref{join}, there is a $1$-cycle of $\Tilde{\Gamma}_m^{IV} \bigcap \Tilde{(R_k)^*}$ in the $O$-copy of the intersection of the star of $A$ with $R_k$ linked with a $1$-axis $a$ containing the $O$-copy of the segment $A$. Using the $(1,1)$-pair $(c,a)$, we form a $(2,1)$-pair $(\hat{c},a)$ by taking $\hat{c}$ to be a suspension over $c$ lying in the suspension over the $O$-copy of the star of $A$ and belonging to $\Tilde{\Gamma}_m^{IV}$. The linking properties of $c$ and $a$ ensure that the $\mathbb{Z}_2$-linking number of $\hat{c}$ and $a$ is $1$.
\end{sloppypar}

\subsubsection{Pairs of type $(2,1)$ or $(1,2)$ - d}

\label{susp21d}

\begin{sloppypar}
    We want to apply the process described in Sections \ref{susp21b} and \ref{susp21c} to a $(1,1)$-pair described in Section \ref{wallregu} and whose cycle lies in the patchworked surface contained in the symmetric copies of $R_k$. When the cycle of the pair is contained in a single orthant, this is possible and we obtain a $(2,1)$-pair. However, the $1$-cycle of such a pair may not be contained in a single orthant. Whenever this is the case, the modification performed in \ref{wallregu} ensures that the dual $1$-axis is contained in a single orthant and, following the process described in Section \ref{susp21a}, we take a suspension over this $1$-axis to obtain a $(1,2)$-pair.
\end{sloppypar}

\subsubsection{Pairs of type $(2,1)$ - e}

\label{susp21e}

\begin{sloppypar}
    We want to apply the process described in Section \ref{susp21b} to the $(1,1)$-pairs described in Sections \ref{suspspecial} and \ref{wallspecial} associated to the patchworked surface lying in the symmetric copies of $R_k$. There is one slight irregularity to notice.
    \end{sloppypar}
    \begin{sloppypar}
    In the case of the odd IV construction (we assume $k$ to be odd), considering the surface lying in the symmetric copies of $R_k$, neither the cycles nor the axes of the pairs described in \ref{suspspecial} are contained in a single orthant. However, the choice of the segment $E_1$ at the penultimate step of the description of the odd IV triangulation in Section \ref{construction3} ensures that the cycles from these pairs are contained either in the hyperplane $\{ x_4=m-k \}$ or in the hyperplane $\{ x_4=-(m-k) \}$. Now, the choice of the vertices of the cones over $R_k$ at the second step of the description of the odd IV triangulation is such that the suspension over all the symmetric copies of $R_k$ lying in $\{ x_4=m-k \}$ (respectively in $\{ x_4=-(m-k) \}$) have the same vertices. Hence, there is a suspension over the desired $1$-cycles that belongs to $\Tilde{\Gamma}_m^{IV}$, providing $(2,1)$-pairs.
    \end{sloppypar}
    \begin{sloppypar}
    When $k$ is even and $\Delta_m^4$ is endowed with the even IV datum, the chosen triangulation and sign distribution ensure that all the cycles from pairs described in Sections \ref{suspspecial} (see Remark \ref{rmksuspspecial}) and \ref{wallspecial} are contained in a single orthant.
\end{sloppypar}

\subsubsection{Number of cycles described in the suspension over a tetrahedron}

\begin{sloppypar}
    In the case of the odd IV construction, we described $\frac{1}{48}m^4-\frac{1}{6}m^3+\frac{5}{12}m^2-\frac{1}{3}m$ pairs of type $(3,0)$ and $\frac{5}{48}m^4 - \frac{1}{2}m^3+\frac{5}{6}m^2-\frac{1}{2}m$ pairs of type $(1,2)$ or $(2,1)$.\\
    %\begin{itemize}
    %    \item $\frac{1}{48}m^4-\frac{1}{6}m^3+\frac{5}{12}m^2-\frac{1}{3}m$ pairs of type$(3,0)$;
    %    \item $\frac{5}{48}m^4 - \frac{1}{2}m^3+\frac{5}{6}m^2-\frac{1}{2}m$ pairs of type $(1,2)$ or $(2,1)$.
    %\end{itemize}
    In the case of the even IV construction, we described $\frac{1}{48}m^4 - \frac{1}{4}m^3 + \frac{25}{24}m^2 - \frac{7}{4}m + 1$ pairs of type $(3,0)$ and $\frac{5}{48}m^4 - \frac{11}{12}m^3 + \frac{71}{24}m^2 - \frac{49}{12}m + 2$ pairs of type $(1,2)$ or $(2,1)$.
    %\begin{itemize}
    %    \item $\frac{1}{48}m^4 - \frac{1}{4}m^3 + \frac{25}{24}m^2 - \frac{7}{4}m + 1$ pairs of type $(3,0)$;
    %    \item $\frac{5}{48}m^4 - \frac{11}{12}m^3 + \frac{71}{24}m^2 - \frac{49}{12}m + 2$ pairs of type $(1,2)$ or $(2,1)$.
    %\end{itemize}
\end{sloppypar}

\subsection{Other pairs whose cycle or axis is contained in a horizontal hyperplane}

\label{suspother}

\begin{sloppypar}
    Let $k$ be an integer such that $1 \leq k \leq m-1$. If $k$ is even, the odd IV triangulation does not contain a union of simplices forming a suspension over $R_k$. Likewise, if $k$ is odd, the even IV triangulation does not contain a union of simplices forming a suspension over $R_k$. However, we can still use the cycle-axis pairs described in Section \ref{prelim} of the patchworked surface obtained by restricting the odd or even IV data to $R_k$ to form cycle-axis pairs of $\Tilde{\Gamma}_m^{IV}$. This is how we treat the parts of type $(00,01,02,12,22)$, $(00,01,11,21,22)$ and $(00,10,20,21,22)$ as well as the walls of type $(00,01,11,22)$ and $(00,11,21,22)$. In this section,  we consider two cases: either $k$ is even and $\Delta_m^4$ is endowed with the odd IV datum, or $k$ is odd and $\Delta_m^4$ is endowed with the even IV datum.
\end{sloppypar}

\subsubsection{Pairs in parts of type $(00,01,02,12,22)$ or $(00,10,20,21,22)$}

\label{0001021222}

\begin{sloppypar}
    Let $k$ be an integer such that $1 \leq k \leq m-1$. In the odd or even IV triangulation, the union of four parts $P_1$, $P_2$, $P_3$ and $P_4$ of type $(00,01,02,12,22)$ (respectively, four parts of type $(00,10,20,21,22)$) sharing a triangle of type $(00,01,02)$ (respectively, a triangle of type $(20,21,22)$) lying in $R_k$ is a double suspension over that triangle. As $R_k$ is a copy of $\Delta_k^3$, Section \ref{surfsusp} describes $(2,0)$-pairs and $(1,1)$-pairs lying in the intersection of $\Tilde{\Gamma}_m^{IV}$ with the union of the symmetric copies of the $3$-faces of $P_1$, $P_2$, $P_3$ and $P_4$ that are contained in $R_k$. We use these pairs to form $(3,0)$- and $(2,1)$-pairs lying in the union of the symmetric copies of $\bigcup_{i=1}^4 P_i$.
\end{sloppypar}
\begin{sloppypar}
    As in \ref{susp30}, in the star of each integer point lying in the relative interior of $\bigcup_{i=1}^4 P_i$, we exhibit a $3$-cycle homeomorphic to a $3$-sphere.
\end{sloppypar}
\begin{sloppypar}
    As in \ref{susp21a}, from a $(2,0)$-pair associated to $\bigcup_{i=1}^4 P_i$ described in \ref{susp20}, we can form a $(2,1)$-pair by taking a suspension over the axis of the pair.
\end{sloppypar}
\begin{sloppypar}
    As in \ref{susp21b}, we can form $(2,1)$-pairs by taking a suspension over a cycle from a $(1,1)$-pair described in Sections \ref{11a}, \ref{11b}, \ref{11c}, \ref{11d}.
\end{sloppypar}
\begin{sloppypar}
    As in \ref{susp21e}, we can form $(2,1)$-pairs by taking a suspension over a cycle from a $(1,1)$-pair described in Section \ref{suspspecial}. Notice that in the presently considered cases, the restriction of the odd (respectively even) IV sign distribution to $R_k$ is constant, \textit{i.e.} it is non-Harnack on every of the triangles used to slice $R_k$ to endow it with an IV triangulation. Hence, as explained in Section \ref{suspspecial}, the $1$-cycles from the $(1,1)$-pairs described in Section \ref{suspspecial} are contained in a single orthant (see Remark \ref{rmksuspspecial}).
\end{sloppypar}

%\paragraph{Pairs of type $(3,0)$}

%\begin{sloppypar}
%    As in \ref{susp30}, in the star of each integer point lying in the relative interior of $P \bigcup P'$, we exhibit a $3$-cycle homeomorphic to a $3$-sphere.
%\end{sloppypar}

%\paragraph{Pairs of type $(2,1)$ - a}

%\begin{sloppypar}
%    As in \ref{susp21a}, from a $(2,0)$-pair associated to $P \bigcup P'$ described in \ref{susp20}, we can form a $(2,1)$-pair by taking a suspension over the axis of the pair.
%\end{sloppypar}

%\paragraph{Pairs of type $(2,1)$ - b}

%\begin{sloppypar}
%    As in \ref{susp21b}, we can form $(2,1)$-pairs by taking a suspension over a cycle from a $(1,1)$-pair described in Sections \ref{11a}, \ref{11b}, \ref{11c}, \ref{11d}.
%\end{sloppypar}

%\paragraph{Pairs of type $(2,1)$ - c}
%\begin{sloppypar}
%    As in \ref{susp21e}, we can form $(2,1)$-pairs by taking a suspension over a cycle from a $(1,1)$-pair described in Section \ref{suspspecial}. Notice that in the presently considered cases, the restriction of the odd (respectively even) IV sign distribution to $R_k$ is constant, \textit{i.e.} it is non-Harnack on every of the triangles used to slice $R_k$ to endow it with an IV triangulation. Hence, as explained in Section \ref{suspspecial}, the $1$-cycles from the $(1,1)$-pairs described in Section \ref{suspspecial} are contained in a single orthant (see Remark \ref{rmksuspspecial}).
%\end{sloppypar}

\subsubsection{Pairs of type $(2,1)$ in parts of type $(00,01,11,21,22)$}

\label{0001112122}

\begin{sloppypar}
    Let $k$ be an integer such that $1 \leq k \leq m-1$. 
    %We consider two cases: either $k$ is even and $\Delta_m^4$ is endowed with the odd IV datum, or $k$ is odd and $\Delta_m^4$ is endowed with the even IV datum. 
    In the odd or even IV triangulation, the union of two parts $P$ and $P'$ of type $(00,01,11,21,22)$ sharing a tetrahedron of type $(00,01,21,22)$ lying in $R_k$ is a suspension over this tetrahedron. Considering $R_k$ as a copy of $\Delta_k^3$, Section \ref{join} describes $(1,1)$-pairs whose cycles lie in the intersection of $\Tilde{\Gamma}_m^{IV}$ with the symmetric copies of $R_k \bigcap (P \bigcup P')$. As in \ref{susp21c}, we can take a suspension over the cycles from these pairs to form $(2,1)$-pairs.
\end{sloppypar}

\subsubsection{Pairs associated to walls of type $(00,01,11,22)$ or $(00,11,21,22)$}

\label{wall1}
\begin{sloppypar}
    Let $k$ be an integer such that $1\leq k \leq m-1$. 
    %We suppose that either $k$ is even and $\Delta_m^4$ is endowed with the odd IV datum, or $k$ is odd and $\Delta_m^4$ is endowed with the even IV datum.
    Let $b$ be an integer point lying in the relative interior of a segment $S$ of type $(00,01)$ (respectively $(21,22)$) lying in $R_k$. Let $e$ be the integer vertex of type $22$ (respectively of type $00$). Let $A$ be the segment $[b,e]$. Such a segment $A$ is contained in two parts of type $(00,01,02,12,22)$ (respectively $(00,10,20,21,22)$), two parts of type $(00,01,11,12,22)$ (respectively $(00,10,11,21,22)$) and two parts of type $(00,01,11,21,22)$. The segment $A$ is also contained in two walls of type $(00,01,11,22)$ (respectively $(00,11,21,22)$). The union of the six mentioned parts can be viewed as a sort of suspension over the union of the tetrahedra of type $(00,01,02,22)$ (respectively $(00,20,21,22)$) and $(00,01,21,22)$ shared by two parts of the same type and contained in $R_k$. Using Sections \ref{wallregu} and \ref{wallspecial}, we can describe two cycle-axis $(1,1)$-pairs of the patchworked surface lying in the symmetric copies of $R_k$ whose axis contains a copy of $A$ as its essential part (when the lattice length of $S$ is even, exactly two different segments with a common extremity share the same cycle-axis pairs, see Section \ref{wallspecial}). There are four triangulation vertices of two different parities in the star of $A$ that do not lie in $R_k$: two points of type $11$ and their closest neighbours in segments of type $(11,12)$ (respectively $(10,11)$). When choosing the sign distribution, we have set the signs of those points such that in the orthants where there is a $(1,1)$-pair, the copies of the four points bear the same sign. Depending on the properties of the $(1,1)$-pair, one can either take a suspension over its axis to obtain a $(1,2)$-pair or a suspension over its cycle to obtain a $(2,1)$-pair.
\end{sloppypar}

\subsubsection{Number of pairs}

\begin{sloppypar}
    In the case of the odd IV construction, in this Section, we described $\frac{1}{48}m^4 - \frac{1}{4}m^3 + \frac{25}{24}m^2 - \frac{7}{4}m + 1$ pairs of type $(3,0)$ and $\frac{5}{48}m^4 - \frac{11}{12}m^3 + \frac{71}{24}m^2 - \frac{49}{12}m + 2$ pairs of type $(1,2)$ or $(2,1)$.\\
    %\begin{itemize}
    %    \item $\frac{1}{48}m^4 - \frac{1}{4}m^3 + \frac{25}{24}m^2 - \frac{7}{4}m + 1$ pairs of type $(3,0)$;
    %    \item $\frac{5}{48}m^4 - \frac{11}{12}m^3 + \frac{71}{24}m^2 - \frac{49}{12}m + 2$ pairs of type $(1,2)$ or $(2,1)$.
    %\end{itemize}
    In the case of the even IV construction, we described $\frac{1}{48}m^4-\frac{1}{6}m^3+\frac{5}{12}m^2-\frac{1}{3}m$ pairs of type$(3,0)$ and $\frac{5}{48}m^4 - \frac{1}{2}m^3+\frac{5}{6}m^2-\frac{1}{2}m$ pairs of type $(1,2)$ or $(2,1)$.
    %\begin{itemize}
    %    \item $\frac{1}{48}m^4-\frac{1}{6}m^3+\frac{5}{12}m^2-\frac{1}{3}m$ pairs of type$(3,0)$;
    %    \item $\frac{5}{48}m^4 - \frac{1}{2}m^3+\frac{5}{6}m^2-\frac{1}{2}m$ pairs of type $(1,2)$ or $(2,1)$.
    %\end{itemize}
\end{sloppypar}

\subsection{Pairs of type $(2,1)$ in parts of type $(00,01,11,12,22)$ or $(00,10,11,21,22)$}

\label{0001111222}

\begin{sloppypar}
    Let $k$ be an integer such that $1 \leq k \leq m$. We consider two cases: either $k$ is even and $\Delta_m^4$ is endowed with the odd IV datum, or $k$ is odd and $\Delta_m^4$ is endowed with the even IV datum. The union of two parts $P$ and $P'$ of type $(00,01,11,12,22)$ (respectively two parts of type $(00,10,11,21,22)$) sharing a tetrahedron of type $(00,01,11,12)$ (respectively of type $(10,11,21,22)$) is a suspension over this tetrahedron. The polytope $P \bigcup P'$ contains a single segment $S$ of type $(00,01)$ (respectively of type $(21,22)$) and a single segment $S'$ of type $(11,12)$ (respectively of type $(10,11)$). Let $b$ and $b'$ be integer points lying respectively in the relative interior of $S$ and $S'$. The segment $A:=[b,b']$ is a triangulation edge. Aside from $b$ and $b'$, the star of $A$ contains $6$ integer vertices of $3$ different parities: the extremities $v_2$ and $v_2'$ of parity $p_2$ (respectively $v_3$ and $v_3'$ of parity $p_3$) of the triangulation segments adjacent to $b$ (respectively $b'$) and lying in $S$ (respectively $S'$) and the vertices $v_4$ and $v_4'$ of parity $p_4$ of the suspension over the join of $S$ and $S'$. Let $Q$ be any $4$-simplex in the star of $A$. Its vertices are $b$, $b'$, and one point of each parity $p_2$, $p_3$, $p_4$. We apply Lemma \ref{linsystem} to $Q$: up to inversion of all signs, there exists a unique orthant $O$ such that the sign of the $O$-copies of $b$ and $b'$ is $+$ while the sign of the $O$-copies of the other vertices of $Q$ is $-$. This means that the star of the $O$-copy of $A$ contains a cell complex $c$ included in $\Tilde{\Gamma}_m^{IV}$, homeomorphic to a $2$-sphere and bounding an obvious $3$-disk that meets the $O$-copy of $A$ in exactly one point. To form an axis $a$ dual to $c$, we distinguish several cases.
    \end{sloppypar}
    \begin{sloppypar}
    If $S$ is of odd lattice length, we take $a$ to be the union of the $O$-copy of $A$, the $O$-copy of the segment whose extremities are $b'$ and the extremity of $S$ that is of the same parity as $b$ and the copies (in $O$ or in an adjacent orthant) of two segments lying in the triangle included in $R_k$, parallel to $\{x_1+x_2+x_3+x_4=m\}$ and containing $S$.
    \end{sloppypar}
    \begin{sloppypar}
    If $S$ is of even lattice length and if $\Delta_m^4$ is endowed with the odd IV triangulation (in this case, $S'$ is of type $(10,11)$ and is of odd lattice length), we take $a$ to be the union of the $O$-copy of $A$, the $O$-copy of the segment whose extremities are $b$ and the extremity of $S'$ that is of the same parity as $b'$ and the copies of two triangulation segments lying in the $O$- or $O'$-copy of $R_{k+1} \bigcup \{x_1+x_2+x_3+x_4=m\}$, where $O'$ is an orthant adjacent to $O$. The common extremity of the last two segments lies in the copy of an edge of type $(10,12)$ and is of the same sign as the other extremities of the segments.
    \end{sloppypar}
    \begin{sloppypar}
    If $S$ is of even lattice length, if $\Delta_m^4$ is endowed with the even IV triangulation (which means that $S'$ is also of even lattice length) and if $b$ is of the same parity as the extremities of $S$ (alternatively, if $b'$ is of the same parity as the extremities of $S'$), we take $a$ to be the union of the $O$-copy of $A$, the $O$-copy of a segment whose extremities are $b$ (alternatively, $b'$) and an extremity $e$ of $S$ (alternatively, $S'$), and the copies in $O$ or in an adjacent orthant of two segments with a common extremity lying on the horizontal hyperplane that contains $S$ (alternatively $S'$) and whose other extremities are respectively $b$ (alternatively $b'$) and $e$. The signs of the chosen copies of the extremities of the two last segments bear the same sign.
    \end{sloppypar}
    \begin{sloppypar}
    If $S$ is of even lattice length, if $\Delta_m^4$ is endowed with the even IV triangulation (\textit{i.e.} $S'$ is also of even lattice length) and if $b$ is not of the same parity as the extremities of $S$ and $b'$ is not of the same parity as the extremities of $S'$, the description of the axis is more involved. Choose an integer point $b''$ in the segment of type $(10,11)$. Let $T$ (respectively $T'$) be the triangle lying in $\{ x_4=m-k \}$ (respectively $\{ x_4=m-k+1 \}$ or $\{ x_4=m-k-1 \}$), parallel to $\{ x_1+x_2+x_3+x_4= m \}$ and containing $S$ (respectively $S'$) as a face. We then consider a path between the $O$-copy of $b$ and the copy of $b$ lying in the same horizontal hyperplane as the the $O$-copy of $b$ but distinct from the $O$-copy of $b$. We call $O'$ the orthant containing the desired copy of $b$. The path can be taken as a broken line formed of triangulation segments whose extremities bear the same sign and lying in copies of the suspension over $T$ contained in $R_k$. Likewise, we consider a path between the $O$-copy of $b'$ and the $O'$-copy of $b''$. It can be taken as a broken line formed of triangulation segments whose extremities bear the same sign and lying in copies of the cones over $T'$ contained in $R_{k-1}$ or $R_{k+1}$. The $O'$-copy of the point $p$ of type $21$ bears the same sign as the $O'$-copies of $b$ and $b''$. Moreover, $[e,b]$ and $[e,b'']$ are triangulation segments. The union of the $O$-copy of $A$, of the two broken lines described above and of the $O'$-copies of $[e,b]$ and $[e,b'']$ is the desired axis. The $O$-copy of $A$ is the essential part of this axis.
    \end{sloppypar}
    \begin{center}
    \includegraphics[scale=0.3]{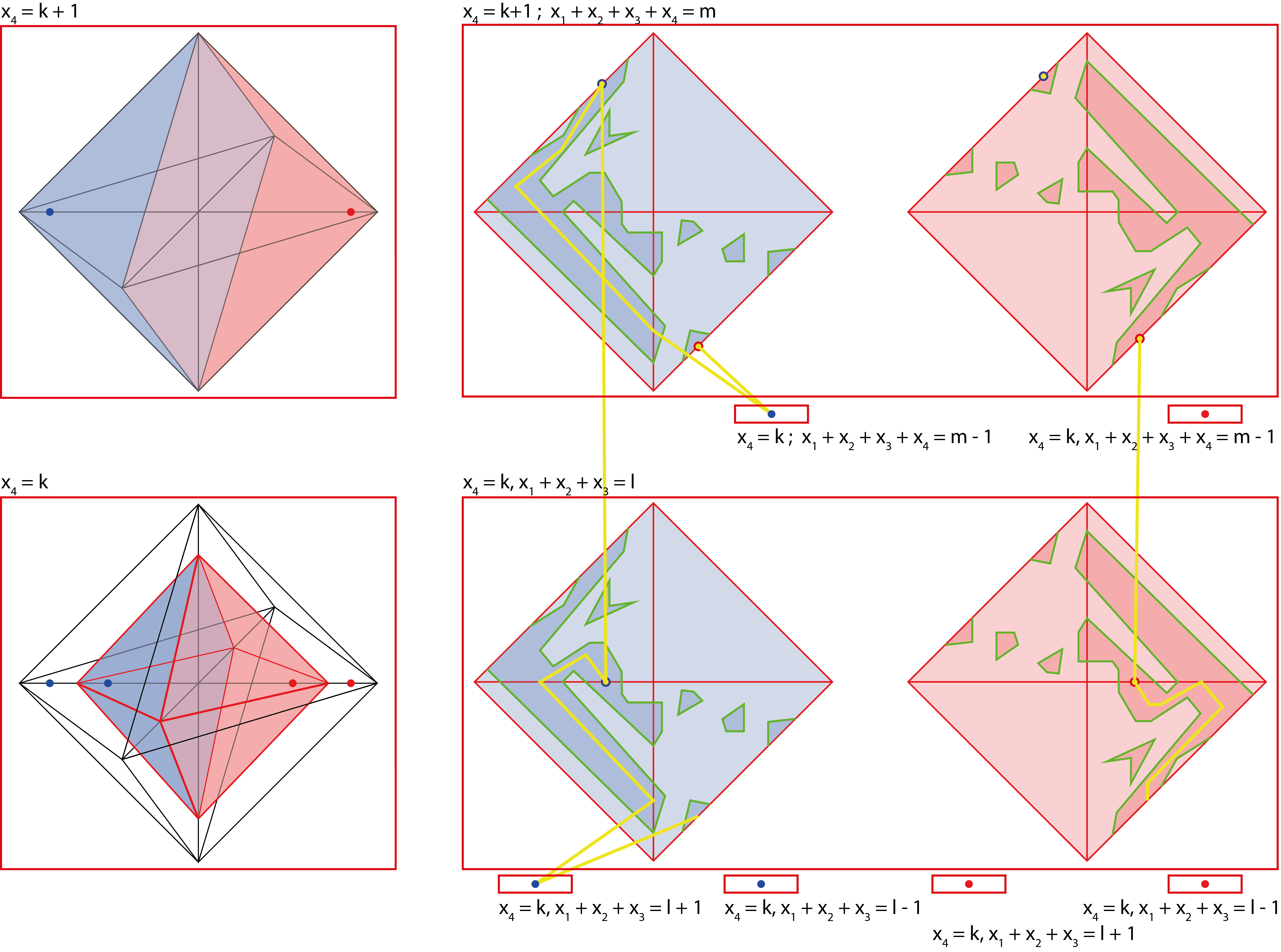}
    \begin{minipage}{16cm}
    \begin{it}
    \small{Fig. 9: Illustration of an axis (in yellow) corresponding to the last set of axes described in Section \ref{0001111222}. These axes pass through several orthants. The equations indicate the hyperplane or intersection of hyperplane in a copy of which the objects lie. The left hand side shows where the considered triangles are located. On the right hand side, the intersection of the patchworked $3$-manifold with the represented triangles appears in green. The blue (respectively red) dots represent vertices of cones over blue (respectively red) triangles. The white dots with red or blue contour are copies of extremities of the essential part of the axis.}
    \end{it}
    \end{minipage}
\end{center}
    \begin{sloppypar}
    The idea behind the description of the last set of axes in the list above is the following: if an axis cannot be completed in the orthant that contains its essential part or in an adjacent orthant, then we travel inside copies of triangles that are contained in a horizontal hyperplane and that are parallel to $\{ x_1+x_2+x_3+x_4=m \}$ until we reach an orthant in which the signs allow us to complete the axis. Figure 9 represents an axis completed in this way. The cycle-axis pairs with such axes are called \textit{special pairs}. Other instances of this process appear in Sections \ref{wall2} and \ref{wall3}.
    \end{sloppypar}
    \begin{sloppypar}
        For parts of type $(00,01,11,12,22)$, we obtain $\frac{1}{24}m^4-\frac{1}{3}m^3+\frac{5}{6}m^2-\frac{2}{3}m$ pairs of type $(2,1)$ in the case of the odd IV datum and $\frac{1}{24}m^4-\frac{1}{6}m^3+\frac{1}{12}m^2+\frac{1}{6}m$ pairs of type $(2,1)$ in the case of the even IV datum.
    \end{sloppypar}
    \begin{sloppypar}
    For parts of type $(00,10,11,21,22)$, we obtain $\frac{1}{16}m^4 - \frac{5}{12}m^3+m^2-\frac{5}{6}m$ pairs of type $(2,1)$ in the case of the odd IV datum and $\frac{1}{16}m^4-\frac{5}{12}m^3+\frac{3}{4}m^2-\frac{1}{3}m$ pairs of type $(2,1)$ in the case of the even IV datum.
\end{sloppypar}

\subsection{Pairs of type $(2,1)$ in parts of type $(00,10,11,12,22)$}

\label{0010111222}

\begin{sloppypar}
    Let $k$ be an integer such that $1 \leq k \leq m$. We consider two cases: either $k$ is even and $\Delta_m^4$ is endowed with the odd IV datum or $k$ is odd and $\Delta_m^4$ is endowed with the even IV datum. We consider the union of two parts $P$ and $P'$ of type $(00,10,11,12,22)$ sharing a tetrahedron of type $(00,10,11,12)$ (respectively of $(10,11,12,22)$). Let $A$ be a segment whose extremities are an integer point $v$ lying in the relative interior of the triangle of type $(10,11,12)$ and the vertex  $v'$ of $P$ and $P'$ of type $00$ (respectively of type $22$). The star of $v$ is contained in $P \bigcup P'$. Aside from $v$ and $v'$, the star of $A$ contains $6$ integer points of $3$ different parities. Let $Q$ be a $4$-simplex in the star of $A$. We apply Lemma \ref{linsystem} to $Q$: up to inversion of all signs, there exists a unique orthant $O$ such that the $O$-copies of $v$ and $v'$ bear the sign $+$, while the other vertices of $Q$ bear the sign $-$. This means that in the star of the $O$-copy of $A$, there is cell complex $c$ contained in $\Tilde{\Gamma}_m^{IV}$, homeomorphic to a $2$-sphere and bounding an obvious $3$-disk disjoint from $\Tilde{\Gamma}_m^{IV}$ that intersects the $O$-copy of $A$ in exactly one point. For instance, we can take $c$ to be the join of the circle that is the intersection of $\Tilde{\Gamma}_m^{IV}$ with the $O$-copy of the star of $v$ in the triangle of type $(10,11,12)$ and of the pair of points formed by the middles of the $O$-copies of the two segments of type $(00,22)$. We form an axis dual to $c$ by taking the union of the following segments:
    \begin{itemize}
        \item the $O$-copy of $A$;
        \item the $O$-copy of the segment whose extremities are $v'$ and an integer vertex $v''$ lying on the boundary of the triangle of type $(10,11,12)$, but not on the edge of type $(11,12)$ (respectively of type $(10,11)$) and that is of the same parity as $v$;
        \item the copies in $O$ or in an adjacent orthant $O'$ of two segments lying in $R_k$ with one common extremity $e$ and whose other extremity is respectively $v$ and $v'$. We take the copies of the segments in the orthant in which their extremities the same sign.
    \end{itemize}
    In the case of the odd IV construction, we obtain $\frac{1}{8}m^4 - \frac{2}{3}m^3 + m^2- \frac{1}{3}m$ pairs of type $(2,1)$.\\
In the case of the even IV construction, we obtain $\frac{1}{8}m^4-\frac{2}{3}m^3+\frac{5}{4}m^2-\frac{5}{6}m$ pairs of type $(2,1)$.
\end{sloppypar}

\subsection{Pairs associated to the wall of the cone over a tetrahedron}

\label{wallcone}

\begin{sloppypar}
    Let $k$ be an integer such that $1 \leq k \leq m$. If $k$ is odd, the odd IV triangulation of $\Delta_m^4$ contains a union of simplices forming a suspension over $R_k$. If $k$ is even and $k<m$, the even IV triangulation contains a union of simplices forming a suspension over $R_k$. As $m$ is even, the even IV triangulation contains a cone over $R_m$. In this section, we suppose that either $k$ is odd and $\Delta_m^4$ is endowed with the odd IV datum or that $k$ is even and $\Delta_m^4$ is endowed with the even IV datum. The cycle-axis $(2,1)$-pairs we are about to describe are similar to those described in the case of surfaces in Section \ref{wallregu}. Choose a cone over $T_k$ that is the boundary of a cone over $R_k$. Let $A$ be the segment whose extremities are an integer point lying in the interior of $T_k$ and the vertex of the considered cone over $T_k$. Let $Q$ be a $4$-simplex in the intersection of the star of $A$ with the cone over $R_k$. We apply Lemma \ref{linsystem} to $Q$: up to inversion of all signs, there exists an orthant $O$ such that the $O$-copies of $v$ and $v'$ bear the sign $+$, while the other vertices of the $O$-copy of $Q$ bear the sign $-$. This means that in the $O$-copy of $T_k$, there is a circle belonging to $\Tilde{\Gamma}_m^{IV}$ surrounding the $O$-copy of $v$ and that there is a cone over this circle lying in the $O$-copy of the intersection of the star of $v$ with $R_k$ and included to $\Tilde{\Gamma}_m^{IV}$. To complete this into a $2$-cycle, we proceed as in Section \ref{wallregu}: either there is a cone over the circle contained in $\Tilde{\Gamma}_m^{IV}$ in the $O$-copy of the star of $A$ minus the cone over $R_k$, or this cone can be found in an adjacent orthant. The union $c$ of the two cones over the circle is a cell complex homeomorphic to a $2$-sphere. To form a $1$-axis dual to the $2$-cycle $c$, take the union of
    \begin{itemize}
        \item the $O$-copy of $A$ (the essential segment of the axis);
        \item the $O$-copy of a segment with vertices $v'$ and an integer point $v''$ lying in a $1$-face of $T_k$ that is of the same parity as $v$;
        \item the copies in an orthant $O'$ adjacent to $O$ of the two segments that respectively have $v''$ and $v$ as a vertex, that lie in the $O'$-copy of $R_k$ and whose other (common) vertex is an integer point whose $O'$-copy bears the same sign as $v$.
    \end{itemize}
    \end{sloppypar}
    \begin{center}
    \includegraphics[scale=0.35]{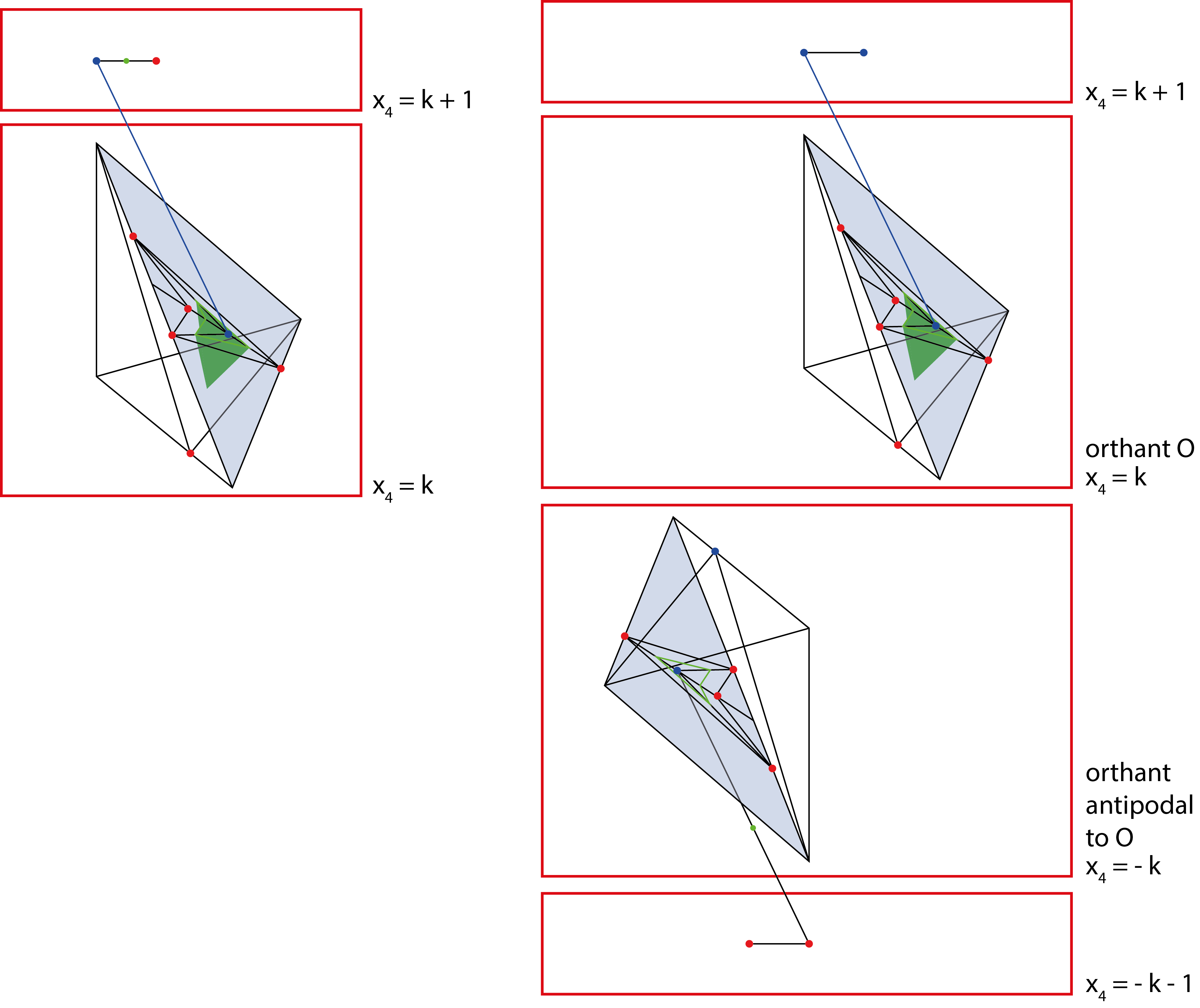}
    \begin{minipage}{16cm}
    \begin{it}
    \small{Fig. 10: On the left, an illustration of a cycle as described in Section \ref{wallcone} that is contained in a single orthant. On the right, an illustration of a cycle as described in Section \ref{wallcone} that is not contained in a single orthant. The equations indicate the hyperplane or intersection of hyperplane in a copy of which the objects lie. The colors of the points represent their sign. The cycle partially appears in green. In both cases, the missing part of the cycle is the cone over the green circle lying in the blue triangle with vertex the green point.}
    \end{it}
    \end{minipage}
\end{center}
    \begin{sloppypar}
    A cycle obtained from the description above bounds an obvious $3$-disk that meets the essential segment of its axis in exactly one point (the $O$-copy of $v$), ensuring the desired linking relations. Figure 10 represents a cycle corresponding to the description above and the essential part of its dual axis.
    
    In the case of the odd IV construction, we obtain $\frac{1}{3}m^3-\frac{3}{4}m^2+\frac{5}{6}m$ pairs of type $(2,1)$.\\
    In the case of the even IV construction, we obtain $\frac{1}{6}m^3 - \frac{3}{4}m^4+\frac{4}{3}m-1$ pairs of type $(2,1)$.
\end{sloppypar}

\subsection{Pairs associated to walls of type $(00,01,12,22)$ or $(00,10,21,22)$}
\label{wall2}

\begin{sloppypar}
    A wall of type $(00,01,12,22)$ (respectively of type $(00,10,21,22)$) is the common face of a part of type $(00,01,02,12,22)$ (respectively $(00,10,20,21,22)$) and a part of type $(00,01,11,12,22)$ (respectively $(00,10,11,21,22)$). We suppose that the segment $S$ of type $(00,01)$ (respectively $(21,22)$) does not lie on $\{ x_1+x_2+x_3+x_4=m \}$. If $S$ is of even lattice length, choose an integer point $p$ in the relative interior of $S$. Let $b$ be another integer point in the relative interior of $S$. Let $e$ be the point of type $12$ (respectively $10$). Let $A$ be the segment $[b,e]$. The segment $A$ is contained in two walls of type $(00,01,12,22)$ (respectively $(00,10,21,22)$), two parts of type $(00,01,02,12,22)$ (respectively $(00,10,20,21,22)$) and two parts of type $(00,01,11,12,22)$ (respectively $(00,10,11,21,22)$). The description of the cycles is as in Section \ref{wallcone}: we apply \mbox{Lemma \ref{linsystem}} to a $4$-simplex lying in the intersection of the star of $A$ with the parts of type $(00,01,02,12,22)$ (respectively $(00,10,20,21,22)$). We denote by $O$ the solution orthant. We obtain a $2$-disk that can be completed into a $2$-cycle either in $O$ or in an adjacent orthant (the arguments are the same as in Section \ref{wallcone}). The dual axis is the union of the following segments:
    \begin{itemize}
        \item the $O$-copy of $A$;
        \item if $S$ is of odd lattice length, the $O$-copy of the segment whose extremities are $e$ and the extremity of $S$ whose parity is the same as the parity of $b$;
        \item if $S$ is of even lattice length and $b$ is of the same parity as the extremities of $S$ (respectively, if $S$ is of even lattice length and the parity of $b$ is different from the parity of the extremities of $S$), the $O$-copy of the segment whose extremities are $e$ and an extremity of $S$ (respectively, the $O$-copy of $[e,p]$);
        %\item if $S$ is of even lattice length and the parity of $a$ is different from the parity of the extremities of $S$, the $O$-copy of $[e,p]$;
        \item two segments with extremities bearing the same sign and lying in a copy of the triangle of type $(00,01,02)$ (respectively $(20,21,22)$) with a common extremity and whose other extremities are $v$ and $v'$ respectively, where $v$ and $v'$ are the extremities of the two previously described segments.
    \end{itemize}
    We now suppose that $S$ is of even lattice length and we consider the segment $[p,e]$. The process described below allows one to find a cycle contained in the symmetric copies of the star of this segment. The axis, however, differs. It is the union of the segments described below.
    \begin{itemize}
        \item Take the $O$-copy of $[p,e]$.
        \item Choose an integer point $v$ lying in the relative interior of the segment of type $(01,02)$ (respectively $(20,21)$) and whose parity is different from the parity of the extemities of $S$ and whose $O$-copy bears the same sign as the $O$-copy of $e$. Take the $O$-copy of $[e,v]$. 
        \item We use the trick mentionned in Section \ref{0001111222}: there exists a path from the $O$-copy of $p$ to the $O$-copy of $v$. It can be taken to be a broken line formed of triangulation edges whose extremities bear the same sign and which lie in the cones over the symmetric copies of the triangle of type $(00,01,02,12,22)$ (respectively $(00,10,20,21,22)$) that are in a single copy of the hyperplane $\{ x_4=m-k \}$.
    \end{itemize}
    This last pair is called a special pair. The axis from such a pair is represented in Figure 11.
    \begin{center}
    \includegraphics[scale=0.4]{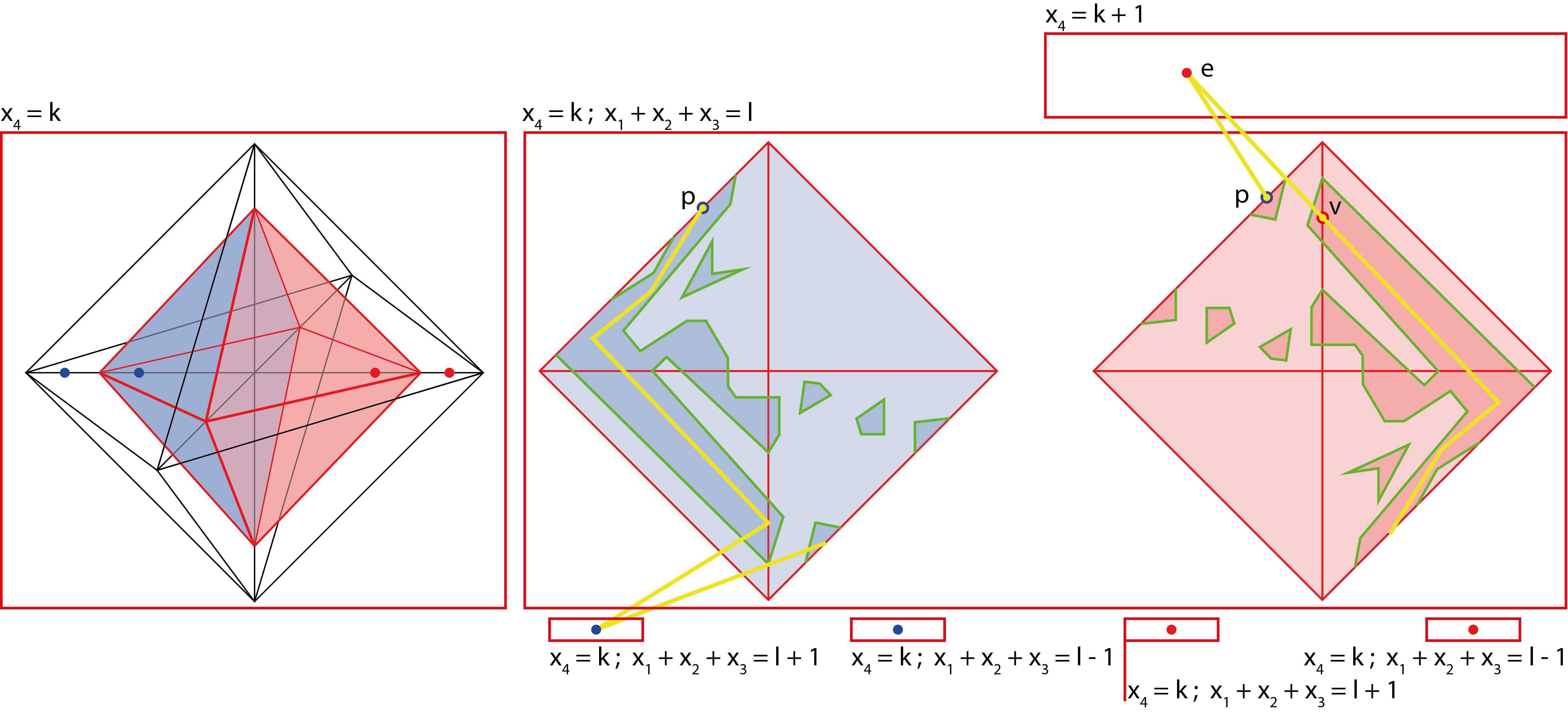}
    \begin{minipage}{16cm}
    \begin{it}
    \small{Fig. 11: Illustration of an axis (in yellow) corresponding to the description of the last set of axes in Section \ref{wall2}. This axis pass through several orthants. The equations indicate the hyperplane or intersection of hyperplane in a symmetric copy of which the objects lie. The left hand side shows where the considered objects are located. On the right hand side, the intersection of the patchworked $3$-manifold with the represented triangles appears in green. The blue (respectively red) points represent vertices of cones over blue (respectively red) triangles. A copy of any integer point that was named in Section \ref{wall2} bears the name of this point.}
    \end{it}
    \end{minipage}
\end{center}
    For walls of type $(00,01,12,22)$, we obtain $\frac{1}{12}m^3-\frac{1}{4}m^2+\frac{1}{6}m$ pairs of type $(2,1)$ in the case of the odd IV datum and $\frac{1}{12}m^3-\frac{1}{4}m^2+\frac{1}{6}m$ pairs of type $(2,1)$ in the case of the even IV datum.
    For walls of type $(00,10,21,22)$, we obtain $\frac{1}{12}m^3 - \frac{1}{2}m^2+\frac{7}{6}m-1$ pairs of type $(2,1)$ in the case of the odd IV datum and $\frac{1}{2}m^2+\frac{2}{3}m$ pairs of type $(2,1)$ in the case of the even IV datum.    
\end{sloppypar}

\subsection{Pairs associated to walls of type $(00,11,12,22)$ or $(00,10,11,22)$}
\label{wall3}

\begin{sloppypar}
    \textit{Mutatis mutandis}, we use the same arguments as in the previous Section \ref{wall2}. We look at two adjacent walls of type $(00,11,12,22)$ (respectively $(00,10,11,22)$) sharing a triangle of type $(00,11,12)$ (respectively $(10,11,22)$). Let $b$ be an integer point in the relative interior of the segment $S$ of type $(11,12)$ (respectively $(10,11)$). Let $e$ be the point of type $00$ (respectively $22$). We look for cycles in the copies of the star of $A:=[b,e]$. The star of $A$ is contained in parts of type $(00,01,11,12,22)$ (respectively $(00,10,11,21,22)$) and $(00,10,11,12,22)$.  This time, the axes have to be completed using segments lying in the star of the copies of the triangle of type $(10,11,12)$. One of the obtained pairs is a special pair.
    \end{sloppypar}
    \begin{sloppypar}
    For walls of type $(00,10,11,22)$, we obtain $\frac{1}{6}m^3-\frac{3}{4}m^2+\frac{5}{6}m$ pairs of type $(2,1)$ in the case of the odd IV datum and $\frac{1}{6}m^3-\frac{1}{2}m^2+\frac{1}{3}m$ pairs of type $(2,1)$ in the case of the even IV datum.
    For walls of type $(00,11,12,22)$, we obtain $\frac{1}{6}m^3-\frac{1}{4}m^2-\frac{1}{6}m$ pairs of type $(2,1)$ in the case of the odd IV datum and $\frac{1}{6}m^3-\frac{1}{2}m^2+\frac{1}{3}$ pairs of type $(2,1)$ in the case of the even IV datum.
\end{sloppypar}

\subsection{Defining an order on the collections of pairs}

\subsubsection{Linking relations}

\label{link3manif}

\begin{sloppypar}
    We consider the pairs in $\mathcal{C}_2$ described above. Each $2$-cycle is linked with its dual axis. There are three types of other linking relations.
    \begin{itemize}
        \item Let $(c,a)$ and $(c',a')$ be two different $(2,1)$-pairs described in any of the Sections \ref{0001021222}, \ref{0001112122}, \ref{wall1}, \ref{susprk}. Then the $\mathbb{Z}_2$-linking numbers of $c$ and $a'$ and of $c'$ and $a$ can be deduced from the linking relations described in Section \ref{link}, as $(c,a)$ and $(c',a')$ are obtained by taking the suspension over cycles from $(1,1)$-pairs described in \ref{prelim}.
        \item Let $(c,a)$ be a $(2,1)$-pair described in any of the Sections \ref{0001021222}, \ref{0001112122}, \ref{wall1}, \ref{susprk}. Let $(c',a')$ be a $(2,1)$-pair described in any of the Sections \ref{wall2}, \ref{wall3}, \ref{wallcone}, \ref{0010111222}, \ref{0001111222}. Then the linking number of $a'$ and $c$ may be $1$, as the segments used to form the axis $a'$ may contain segments in the essential part of the axis $a$.
        \item Let $(c,a)$ and $(c',a')$ be two different $(2,1)$-pairs described in the same Section among the following ones: \ref{0010111222}, \ref{wall2}, \ref{wall3}. If $(c,a)$ is a special pair and $(c',a')$ is not, then the $\mathbb{Z}_2$-linking number of $a$ and $c'$ may be $1$.
    \end{itemize}
    We now turn to the pairs in $\mathcal{C}_1$, which are all described in Sections \ref{susp21d} and \ref{wall1}. Their linking properties can be deduced from the linking relations described in Section \ref{link}, as the $(1,2)$-pairs are obtained by taking the suspension over axes from $(1,1)$-pairs described in \ref{prelim}.
\end{sloppypar}

\subsubsection{Partial order}

\label{order3manif}

\begin{sloppypar}
    On $\mathcal{C}_2$, we take the partial order generated by the following relations.
    \begin{itemize}
        \item Let $(c,a)$ be a non-special pair described in \ref{wall3}. Let $(c',a')$ be a special pair described in \ref{wall3}. Then $(c,a)<(c',a')$.
        \item Let $(c,a)$ be a non-special pair described in \ref{wall2}. Let $(c',a')$ be a special pair described in \ref{wall2}. Then $(c,a)<(c',a')$.
        \item Let $(c,a)$ be a non-special pair described in \ref{0010111222}. Let $(c',a')$ be a special pair described in \ref{0010111222}. Then $(c,a)<(c',a')$.
        \item Let $(c,a)$ be a pair described in any of the Sections \ref{wall2}, \ref{wall3}, \ref{wallcone}, \ref{0010111222}, \ref{0001111222}. Let $(c',a')$ be a pair described in any of the Sections \ref{0001021222}, \ref{0001112122}, \ref{wall1}, \ref{susprk}. Then $(c,a)<(c',a')$.
        \item Let $k$ be an integer such that $1\leq k \leq m-1$. If $k$ is odd and $\Delta_m^4$ is endowed with the odd IV datum or if $k$ is even and $\Delta_m^4$ is endowed with the even IV datum, the order defined in \ref{order} on the collection of $(1,1)$-pairs of the patchworked surface associated to $R_k$ induces an order on the pairs described in Sections \ref{susp21b}, \ref{susp21c}, \ref{susp21d}, \ref{susp21e} that are obtained from $(1,1)$-pairs of the patchworked surface in $R_k$.
        \item Let $k$ be an integer such that $1\leq k \leq m-1$. If $k$ is odd and $\Delta_m^4$ is endowed with the even IV datum or if $k$ is even and $\Delta_m^4$ is endowed with the odd IV datum, the order defined in \ref{order} on the collection of $(1,1)$-pairs of the patchworked surface associated to $R_k$ induces an order on the pairs described in Sections \ref{0001021222}, \ref{0001112122} and \ref{wall1} that are obtained from $(1,1)$-pairs associated to $R_k$.
    \end{itemize}
    On $\mathcal{C}_1$ we take the partial order induced by the order defined in \ref{order}, as all $(1,2)$-pairs are obtained by taking the suspension over axes from $(1,1)$-pairs described in \ref{prelim}.
\end{sloppypar}

\subsection{Independence of the described cycles}

\label{independence3}

\begin{sloppypar}
    The classes represented by the cycles from the pairs in $\mathcal{C}_3$ are clearly independent, as the $3$-cycles are distinct connected components of $\widetilde{\Gamma}_m^{IV}$.
    Let $i_2:\mathcal{C}_2 \rightarrow \{1,...,\# \mathcal{C}_2\}$ (respectively, \mbox{$i_1:\mathcal{C}_1 \rightarrow \{1,...,\# \mathcal{C}_1\}$}) be an injective map inducing a total order on $\mathcal{C}_2$ (respectively, on $\mathcal{C}_1$) that refines the partial order defined in Section \ref{order3manif}. Here $\# \mathcal{C}_2$ (respectively, $\# \mathcal{C}_1$) denotes the cardinal of the set $\mathcal{C}_2$ (respectively, $\mathcal{C}_1$). We want to check whether the dimension of the subspace of $H_2(\Tilde{\Gamma}_m^{IV};\mathbb{Z}_2)$ (respectively, $H_1(\Tilde{\Gamma}_m^{IV};\mathbb{Z}_2)$) spanned by the classes represented by the $2$-cycles (respectively, $1$-cycles) from the pairs of $\mathcal{C}_2$ (respectively, $\mathcal{C}_1$) is equal to $\# \mathcal{C}_2$ (respectively, $\# \mathcal{C}_1$). We define the matrix $\mathcal{A}_2$ (respectively, $\mathcal{A}_1$) whose coefficient in the $i$-th row and $j$-th column is the $\mathbb{Z}_2$-linking number of the cycle from the $i$-th pair of $\mathcal{C}_2$ (respectively, $\mathcal{C}_1$) with the axis from the $j$-th pair of $\mathcal{C}_2$ (respectively, $\mathcal{C}_1$). 
    The order defined in Section \ref{order3manif} has been chosen by taking into account the linking relations presented in Section \ref{link3manif} to ensure the following property: if $(c,a)$ and $(c',a')$ are two distinct pairs in $\mathcal{C}_2$ (respectively, in $\mathcal{C}_1$) such that $(c,a)<(c',a')$, then the $\mathbb{Z}_2$-linking  number of $c$ and $a'$ is $0$. Furthermore, any cycle is $\mathbb{Z}_2$-linked with its dual axis. Thus the matrix $\mathcal{A}_2$ (respectively, the matrix $\mathcal{A}_1$) is lower triangular and its diagonal coefficients are all non-zero. By Proposition \ref{indepmatrice}, this implies that the classes represented by the cycles from the pairs in $\mathcal{C}_2$ (respectively, in $\mathcal{C}_1$) form a basis of the subspace of $H_2(\Tilde{\Gamma}^{IV}_m;\mathbb{Z}_2)$ (respectively, $H_1(\Tilde{\Gamma}^{IV}_m;\mathbb{Z}_2)$) they span. In particular, these classes are independent.
\end{sloppypar}

\subsection{Supplementary $3$-cycle and $1$-cycle}

\label{cyclessupp3var}

\begin{sloppypar}
    There are some hyperplane pieces of the hypersurface $\Tilde{\Gamma}_m^{IV}$ that are not part of the described $3$-cycles. Hence, the hypersurface $\Tilde{\Gamma}_m^{IV}$ has at least one more connected component than the number of described $3$-cycles, which means that there is at least one more $3$-cycle representing a class independent from the others.\end{sloppypar}
    \begin{sloppypar}
    Consider the surface defined as the intersection of the $3$-manifold $\Tilde{\Gamma}_m^{IV}$ with the hyperplane $\{x_4=0\}$. The triangulation of $\Delta_m^4$ and the distribution of signs on its vertices restricted to $\{ x_4=0 \}$ coincide with an IV triangulation and an IV sign distribution for surfaces (see Section \ref{prelim}). In particular, the associated surface contains a loop that realizes the non-trivial class in $H_1(\widetilde{(\Delta_m^3)^*};\mathbb{Z}_2)$ (see Section \ref{cyclessupp}). This same loop, seen as a $1$-cycle in the $3$-manifold $\Tilde{\Gamma}_m^{IV} \subset \widetilde{(\Delta_m^4)^*}$, realizes the non-trivial class in $H_1(\widetilde{(\Delta_m^4)^*};\mathbb{Z}_2)$. All the $1$-cycles of our collection realize the trivial class in $H_1(\widetilde{(\Delta_m^4)^*};\mathbb{Z}_2)$, hence this loop realizes a class in $H_1(\Tilde{\Gamma}_m^{IV};\mathbb{Z}_2)$ that is independent from the classes realized by the $1$-cycles in our collection.
\end{sloppypar}

\subsection{Proof of Proposition \ref{maximal}}

\begin{sloppypar}
    \textbf{Proof:} At this point, we have exhibited $\frac{1}{2} (m^4 -5m^3 +10m^2 -10m +8)$ cycles representing homology classes that are independent over $\mathbb{Z}_2$. The collection contains $3$-cycles, $2$-cycles, \mbox{$1$-cycles}.
    Notice that by construction, the intersection of any $2$-cycle from the collection with any $1$-cycle from the collection is empty. For the supplementary $1$-cycle described in Section \ref{cyclessupp3var}, this is a consequence of the fact that the cycle is contained in the surface associated to the IV datum restricted to the hyperplane $\{ x_4=0 \}$. Given any $2$-cycle and any $1$-cycle from the collection, this implies that their $\mathbb{Z}_2$-intersection number is zero, \textit{i.e.} that they are orthogonal with respect to the $\mathbb{Z}_2$-intersection pairing. 
    Hence, each $2$-cycle from the collection defines a linear form belonging to the dual of $H_1(\Tilde{\Gamma}_m^{IV};\mathbb{Z}_2)$ whose evaluation on any $1$-cycle from the collection is zero. The space spanned by these linear forms is a subspace of the orthogonal of the subspace of $H_1(\Tilde{\Gamma}_m^{IV};\mathbb{Z}_2)$ spanned by the $1$-cycles from the collection. Together with the linear independence of the cycles from the collection, this implies that the dimension of $H_1(\Tilde{\Gamma}_m^{IV};\mathbb{Z}_2)$ is greater than or equal to the sum of the number of $2$-cycles and the number of $1$-cycles from the collection. By Poincaré duality, this inequality also holds for the dimension of $H_2(\Tilde{\Gamma}_m^{IV};\mathbb{Z}_2)$.
    \end{sloppypar}
    \begin{sloppypar}
    We can conclude that $\Tilde{\Gamma}_m^{IV}$ has a total $\mathbb{Z}_2$-Betti number greater than or equal to \mbox{$m^4 -5m^3 +10m^2 -10m +8$}. However the maximal possible number is  \mbox{$8 - \sum_{j=0}^{3} (-1)^j m^{j+1} {5 \choose j+2} = m^4 -5m^3 +10m^2 -10m +8$}.
    Hence, the degree $m$ real algebraic hypersurfaces constructed using either the odd or the even IV construction are maximal.\qedsymbol\\
\end{sloppypar}

\section{Proof of Proposition \ref{type} : Itenberg-Viro $3$-manifolds are of type $\chi = \sigma$.}

\label{prooftype}

\begin{sloppypar}
We denote by $X^m$ a non-singular real algebraic hypersurface of even degree $m$ in $\mathbb{RP}^4$ obtained by applying the combinatorial patchworking theorem to $\Delta_m^4$ endowed with the odd or even IV datum. Let $Y^m$ be a double covering of $\mathbb{CP}^4$ branched along $\mathbb{C}X^m$.
Retaining the notations from Section \ref{definitions}, we want to compare the Euler characteristic of $\mathbb{R}Y_-^m$ or $\mathbb{R}Y_+^m$, the fixed point sets of the involutions that lift the complex conjugation to $Y^m$, with the signature of $\mathbb{C}Y^m$.
Recall from Section \ref{definitions} that $\mathbb{R}Y_+^m \setminus \mathbb{R}X^m$ is a double covering of $X^m_+ \setminus \mathbb{R}X$. Hence, as $X^m$ is odd dimensional, we have $\chi(\mathbb{R}Y_+^m) = 2\chi(X_+^m)$. From \hyperref[orevkov]{[Ore02, Lemma 4.2]}, we know that computing $\chi(X_+^m)$ reduces to counting how many symmetric copies with only vertices of sign $+$ each of the simplices of the odd or even IV triangulation has. To do this, we borrow the notion of even and odd simplices from \cite{orevkov}.
\begin{defin}
A $k$-simplex $Q$ in $\mathbb{R}^n$ with integer vertices $v_1$, ..., $v_{k+1}$ is said to be \emph{even} if $\sum_{1 \leq i \leq k+1} v_i \in 2\mathbb{Z}^{n}$. If $Q$ is not even, then it is said to be \emph{odd}.
\end{defin}
\begin{lem}
(\cite{orevkov}) Let $Q \subset \mathbb{R}^n$ be a primitive $n$-simplex. Then $Q$ has exactly one non-empty even face. \hfill \qedsymbol
\end{lem}
In the odd and even IV triangulations, each $4$-simplex contains an even vertex. Hence the only even simplices of these triangulations are the even vertices. To analyze the signs of the empty copies of every simplex, we need the following lemma.
\begin{lem}
\label{invertsym}
Given a primitive simplex $Q$ with integer vertices and no non-empty even face, there exists a composition of symmetries with respect to coordinate hyperplanes that inverts the signs of all vertices of $Q$.
\end{lem}
\textbf{Proof :} We prove the propostion by induction on the dimension of $Q$. If $Q$ is a $0$-simplex, the statement is straightforward. For simplicity, we call any composition of symmetries with respect to coordinate hyperplanes a reflection.\end{sloppypar}
\begin{sloppypar}
Let $k$ be a positive integer. Suppose that $Q$ is a primitive $k$-simplex with no non-empty even face. Suppose that for each primitive $(k-1)$-simplex with no non-empty even face, there exists a reflection that inverts the signs of all its vertices. Then, for each $(k-1)$-face of $Q$, we can find such a reflection. If one of these reflections inverts simultaneously all the vertices of $Q$, we are done. If not, composing these reflections, we obtain new ones that invert the signs of any $2$ vertices of $Q$ and leave the signs of all the other vertices unchanged. Composing these new reflections, we obtain reflections that invert the signs of any $4$ vertices of $Q$ and no other sign. We can repeat this process to obtain reflections that invert the signs of any $2j$ vertices of $Q$ and no other sign for $2j \leq k+1$.\end{sloppypar}
\begin{sloppypar}
If $k$ is odd, we immediately obtain a suitable reflection. If $k$ is even (\textit{i.e.} the number of vertices of $Q$ is odd), as $Q$ is not an even simplex, there exists a coordinate labelled $i$ which is odd for an odd number of vertices of $Q$. The reflection with respect to the $i$-th coordinate hyperplane inverts exactly the signs of these vertices. The number of the remaining vertices is even and strictly less than $k$. Hence, there is another reflection that inverts exactly the signs of these vertices. Composing both reflections, we get a new one that inverts the signs of all vertices simultaneously.\hfill$\qedsymbol$\\
\end{sloppypar}

\begin{sloppypar}
Hence, half of the empty copies of any simplex in the triangulations described in Section \ref{construction3} with no non-empty even face have only vertices of sign $+$, while the other half of the empty copies have only vertices of sign $-$.
%As the triangulation is primitive, the number of empty copies of a $k$-simplex $v$ only depends on $k$ and on the number of $3$-faces of $\Delta_m^4$ the relative interior of $v$ is contained in.
The following lemma is a consequence of Lemma \ref{linsystem}. A proof can be found in \cite{tsurf}.
\begin{lem}
\label{emptycopies}
Let $m$ be a positive even integer. A primitive $k$-simplex $Q$ from the odd or even IV triangulation of $\Delta_m^4$ with signed integer vertices has exactly $2^{n-p-k}$ empty copies where $p$ is the number of $3$-faces of $\Delta_m^4$ the relative interior of $Q$ is contained in. \hfill \qedsymbol
\end{lem}
In the triangulation of $\Delta_m^4$ of the \emph{odd} IV construction, counting the simplices with no even vertex and whose relative interior is contained in the interior of $\Delta_m^4$, we obtain
\begin{itemize}
    \item $\frac{m^4}{2}$ simplices of dimension $3$; each of these has $1$ empty symmetric copy with all vertices of sign $+$;
    \item $\frac{7m^4}{8} - \frac{5m^3}{4}$ simplices of dimension $2$; each of these has $2$ empty symmetric copies with all vertices of sign $+$;
    \item $\frac{5m^4}{12} - \frac{25m^3}{16} + \frac{35m^2}{24}$ simplices of dimension $1$; each of these has $4$ empty symmetric copies with all vertices of sign $+$;
    \item $\frac{5m^4}{128} - \frac{35m^3}{16} + \frac{35m^2}{32} - \frac{24m}{4}$ simplices of dimension $0$; each of these has $8$ empty symmetric copies with all vertices of sign $+$.
\end{itemize}
Counting the simplices with no even vertex and whose relative interior is contained in the relative interior of the $3$-faces of $\Delta_m^4$, we obtain
\begin{itemize}
    \item $\frac{20m^3}{8}$ simplices of dimension $2$; each of these has $1$ empty symmetric copy with all vertices of sign $+$;
    \item $\frac{25m^3}{8} - 5m^2$ simplices of dimension $1$; each of these has $2$ empty symmetric copies with all vertices of sign $+$;
    \item $\frac{35m^3}{48} - \frac{15m^2}{4} + \frac{55m}{12}$ simplices of dimension $0$; each of these has $4$ empty symmetric copies with all vertices of sign $+$.
\end{itemize}
Counting the simplices with no even vertex and whose relative interior is contained in the relative interior of the $2$-faces of $\Delta_m^4$, we obtain
\begin{itemize}
    \item $5m^2$ simplices of dimension $1$; each of these has $1$ empty symmetric copy with all vertices of sign $+$;
    \item $\frac{15m^2}{4} - \frac{15m}{2}$ simplices of dimension $0$; each of these has $2$ empty symmetric copies with all vertices of sign $+$.
\end{itemize}
Counting the simplices with no even vertex and whose relative interior is contained in the relative interior of the $1$-faces of $\Delta_m^4$, we obtain $5m$ simplices of dimension $0$. Each of these has $1$ empty symmetric copy with all vertices of sign $+$.\end{sloppypar}
%\begin{itemize}
%    \item $5m$ simplices of dimension $0$; each of these has $1$ empty symmetric copy with all vertices of sign $+$.
%\end{itemize}
\begin{sloppypar}
    We can now compute the Euler characteristic of $X^m_+$: $\chi(X^m_+) = -\frac{5m^4}{48} + \frac{5m^2}{12}$.
%We can now compute the Euler characteristic of $X^m_+$: $\chi(X^m_+) = -\frac{5m^4}{48} + \frac{5m^2}{12}$.
\end{sloppypar}
\begin{sloppypar}
%\begin{center}$\chi(X^m_+) = -\frac{5m^4}{48} + \frac{5m^2}{12}$.\end{center}
Hence $\chi(\mathbb{R}Y^m_+) = -\frac{5m^4}{24} + \frac{5m^2}{6}$ and $\chi(\mathbb{R}Y^m_-) = \frac{5m^4}{24} - \frac{5m^2}{6} + 2$ as \mbox{$\chi(\mathbb{R}Y^m_-) + \chi(\mathbb{R}Y^m_+) = 2\chi(\mathbb{RP}^4) = 2$}.
Recall from Section \ref{definitions} that $\sigma(Y^m) = \frac{1}{24}(5m^4 - 20m^2 +48)$. We can now conclude: \mbox{$\sigma(Y^m) = \chi(\mathbb{R}Y^m_-)$} and thus $X^m$ is of type $\chi = \sigma$ for each even $m$. 
%Notice that here, the real structure on $Y^m$ that makes $Y^m$ a maximal real algebraic hypersurface is the one whose real part is $\mathbb{R}Y^m_-$.
One can repeat this count with the even IV triangulation. The final result is the same.
This ends the proof of Proposition \ref{type}. \hfill$\qedsymbol$
\end{sloppypar}

\section{Constructing maximal $3$-manifolds of type $\chi \neq \sigma$}
\label{smalldeviation}
\begin{sloppypar}
In this section, we modify the odd Itenberg-Viro construction to obtain maximal even degree hypersurfaces that have more connected components than the IV hypersurfaces of the same degree and that are not of type $\chi = \sigma$. We first describe a construction of maximal surfaces that are not of type $\chi = \sigma$. The goal is to "insert" these surfaces in the odd IV-construction while preserving their interesting properties.
\end{sloppypar}

\subsection{Surfaces of type $\chi \neq \sigma$}

\begin{sloppypar}
We modify an IV triangulation of $\Delta_m^3$ (described in Section \ref{prelim}) to construct non-singular maximal real algebraic surfaces with the Euler characteristic of their real part different from the signature of their complex part. The idea is to glue a part of the original construction from Section \ref{prelim} together with maximal surfaces which appear as double coverings of $\mathbb{CP}^2$ branched along maximal real algebraic curves that are not of type $\chi = \sigma$. We describe the construction when the degree of the surface is equal to $2k+1$ for some integer $k \geq 4$.
%for convenience ank+d brevity. The construction can be adapted to some other cases.
\end{sloppypar}

\subsubsection{Triangulation}
\label{trimodif}
\begin{sloppypar}
Let $k \geq 4$ be an integer and let $m=2k+1$. We triangulate the simplex $\Delta_m^3$. First, triangulate $\Delta_m^3 \bigcap \{ x_3 \geq 2 \}$ using the IV triangulation mentioned at the last step of the description of the odd IV triangulation of $\Delta_m^4$ in Section \ref{construction3}. Retaining the notations from Section \ref{prelim}, this triangulation satisfies the following properties.
    \begin{itemize}
        \item The $2$-face chosen at the first step of the subdivision in Section \ref{triangulation} is \mbox{ $\Delta_m^3 \bigcap \{ x_3=2 \}$}.
        \item The segment $E_1$ chosen at the second step of the subdivsion in Section \ref{triangulation} is \mbox{$\Delta_m^3 \bigcap \{ x_3 \geq 2 \} \bigcap  \{ x_1=0, x_1+x_2+x_3=k \}$}.
        \item The segment $E_2$ chosen at the second step of the subdivsion in Section \ref{triangulation} is \mbox{$\Delta_m^3 \bigcap \{ x_3 \geq 2 \} \bigcap  \{ x_2=0, x_1+x_2+x_3=k \}$}.
    \end{itemize}
Then divide $\Delta_m^3 \bigcap \{ x_3 \leq 2 \}$ into two parts using the plane \mbox{$\{ x_1+x_2+ kx_3 = 2k \}$}. We denote by $\Bar{T}$ the part of the subdivision that is a cone over the triangle with vertices $(0,0,0)$, $(2k,0,0)$, $(0,2k,0)$. We denote by $C$ the remaining part of $\Delta_m^3 \bigcap \{ x_3 \leq 2 \}$.
To triangulate $\Bar{T}$, first triangulate the triangle lying on the hyperplane $\{ x_3=0\}$. The following segments must appear in the triangulation as edges or as unions of edges (see Figure $12$ for an example):
\begin{itemize}
    \item the segments defined as the intersection of $\{ x_1 + x_2 = r \}$ with the basis of $\Bar{T}$, where $r$ is an integer such that $2k-2 \leq r \leq 2k+1$;
    \item the segments $[(2k-1,0,0),(0,2k-2,0)]$, $[(0,2k,0),(2k-1,0,0)]$, $[(0,2k,0),(2k+1,0,0)]$;
    \item a segment $A:=[a,b]$ of lattice length $1$ such that $a$ and $b$ are non-even integer points and $a$ (respectively $b$) lies in the relative interior of $[(0,0,0), (2k-2,0,0)]$ (respectively in the relative interior of $[(0,0,0), (0,2k-2,0)]$).
\end{itemize}
The triangulation of this triangle can then be arbitrarily completed into a primitive convex triangulation.
\begin{center}
    \includegraphics[scale=0.35]{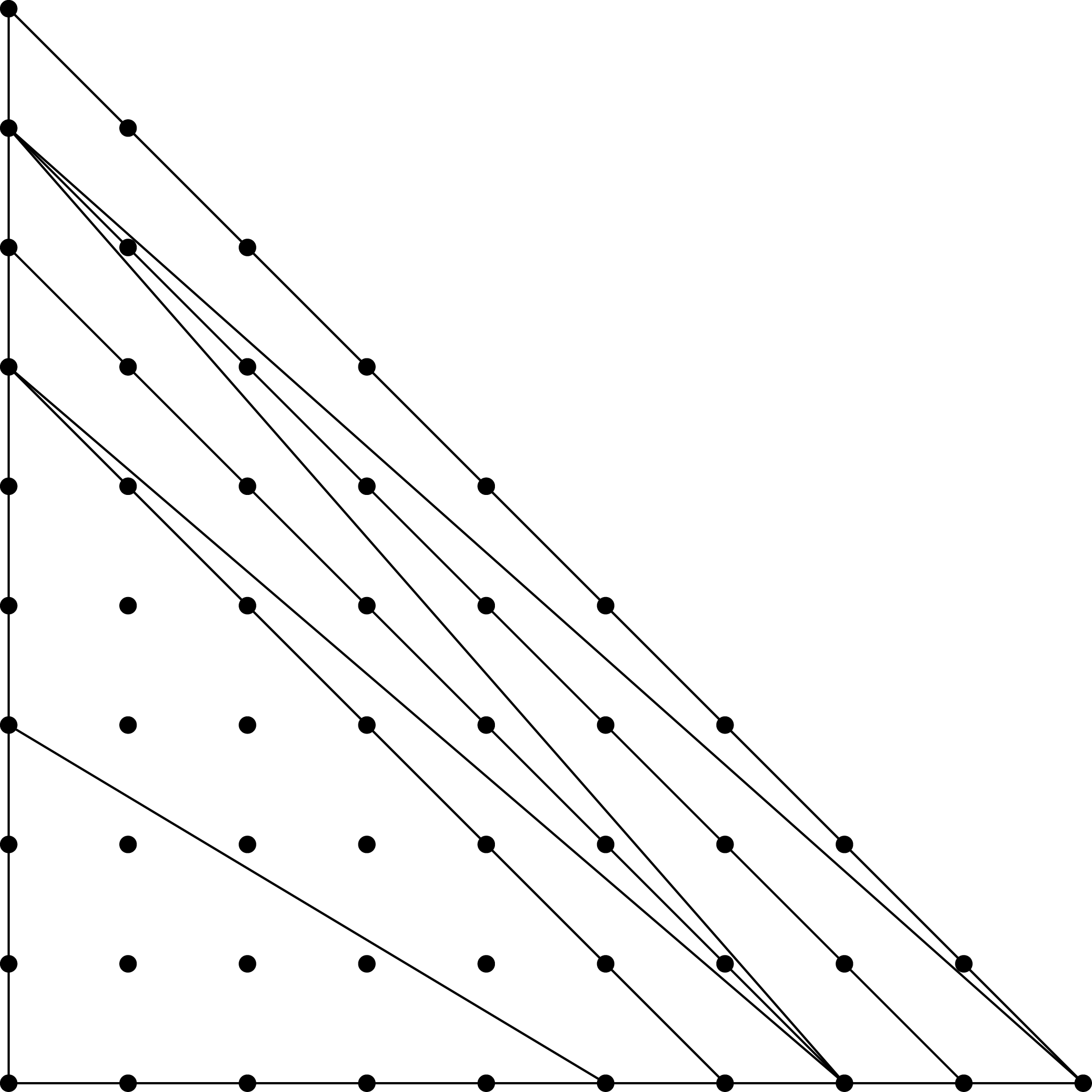}\\
    \begin{minipage}{16cm}
    \begin{it}
    \small{Fig. 12:  The important segments in the $((5,0,0),(0,3,0))$-SD triangulation of $\Delta_9^3 \bigcap \{ x_3=0 \}$, following the description in Section \ref{trimodif}}
    \end{it}
    \end{minipage}
\end{center}
The tetrahedron $\Bar{T}$ is triangulated using the cones of vertex $(0,0,2)$ over the simplices of the triangulation of the basis of $\Bar{T}$. We then divide the tetrahedra obtained in the only possible way to obtain a triangulation such that all integer points contained in $\Bar{T}$ are vertices of the triangulation.\end{sloppypar}
\begin{sloppypar}
To triangulate $C$, take the cone of vertex $(0,k+1,1)$ over the triangle with vertices $(2k,0,0)$, $(0,2k,0)$, $(0,0,2)$. It is naturally triangulated by taking the cone over the triangulation of its basis. Subdivide the remaining part of $C$ using the quadrangle with vertices $(0,k+1,1)$, $(k,0,1)$, $(0,2k,1)$, $(2k,0,1)$, denoted by $C_1$. These subdivision steps are represented on Figure $13$. We triangulate $C_1$ as follows. First, subdivide $C_1$ into quadrangles and one triangle using the segments parallel to $[ (2k,0,1), (0,2k,1) ]$ and whose extremities are integer points lying in the relative interior of the segments $[(k,0,1), (2k,0,1)]$ and $[(0,k+1,1), (0,2k,1)]$. Then, choose a diagonal of one of the quadrangles. Its extremities are points of parity $p_1$ and $p_2$ respectively. Divide each quadrangle using its diagonal whose extremities are integer points of parity $p_1$ and $p_2$ respectively. Finally, complete this subdivision of $C_1$ into a primitive triangulation in the only possible way.
%\begin{itemize}
%    \item Subdivide $C_1$ into quadrangles and one triangle using the segments parallel to $[ (2k,0,1), (0,2k,1) ]$, and whose extremities are integer points lying in the relative interior of the segments $[(k,0,1), (2k,0,1)]$ and $[(0,k+1,1), (0,2k,1)]$.
%    \item Choose a diagonal of one of the quadrangles. Its extremities are points of parity $p_1$ and $p_2$ respectively. Divide each quadrangle using its diagonal whose extremities are integer points of parity $p_1$ and $p_2$ respectively. 
%    \item Complete this subdivision of $C_1$ into a primitive triangulation in the only possible way.
%\end{itemize}
We now consider
\begin{itemize}
    \item the cone with vertex $(0,0,2)$ over $C_1$;
    \item the cone with vertex $(2k,0,1)$ over the triangle with vertices $(0,0,2)$, $(2k-1,0,2)$, $(0,2k-1,2)$;
    \item the join of the segments $[(0,0,2), (0,2k-1,2)]$ and $[(2k,0,1), (0,2k,1)]$;
    \item the cone with vertex $(0,2k,1)$ over the quadrangle with vertices $(2k,0,0)$, $(0,2k,0)$, $(0,2k+1,0)$, $(2k+1,0,0)$;
    \item the cone over $C_1$ with vertex $(2k,0,1)$;
    \item the join of the segments $[(2k,0,0), (2k+1,0,0)]$ and $[(2k,0,1), (0,2k,1)]$;
    \item the join of the segments $[(2k,0,0), (0,2k,0)]$ and $[(0,k+1,1),(0,2k,1)]$.
\end{itemize}
The pieces described above are naturally triangulated as cones over triangulated polygons or joins of triangulated segments. Lastly, remove from the triangulation the vertices which lie strongly inside $\Bar{T}$ and the vertices lying in the relative interior of $(\bar{T} \bigcap \{ x_1=0 \}) \bigcup (\bar{T} \bigcap \{ x_2=0 \})$. The convexity of the resulting triangulations follows from the convexity of the triangulations of the triangle $\Delta_m^3 \bigcap \{ x_3=0 \}$ and of the quadrangle $C_1$ (see \cite[Section~6]{tsurf}).
\begin{center}
    \includegraphics[scale=0.285]{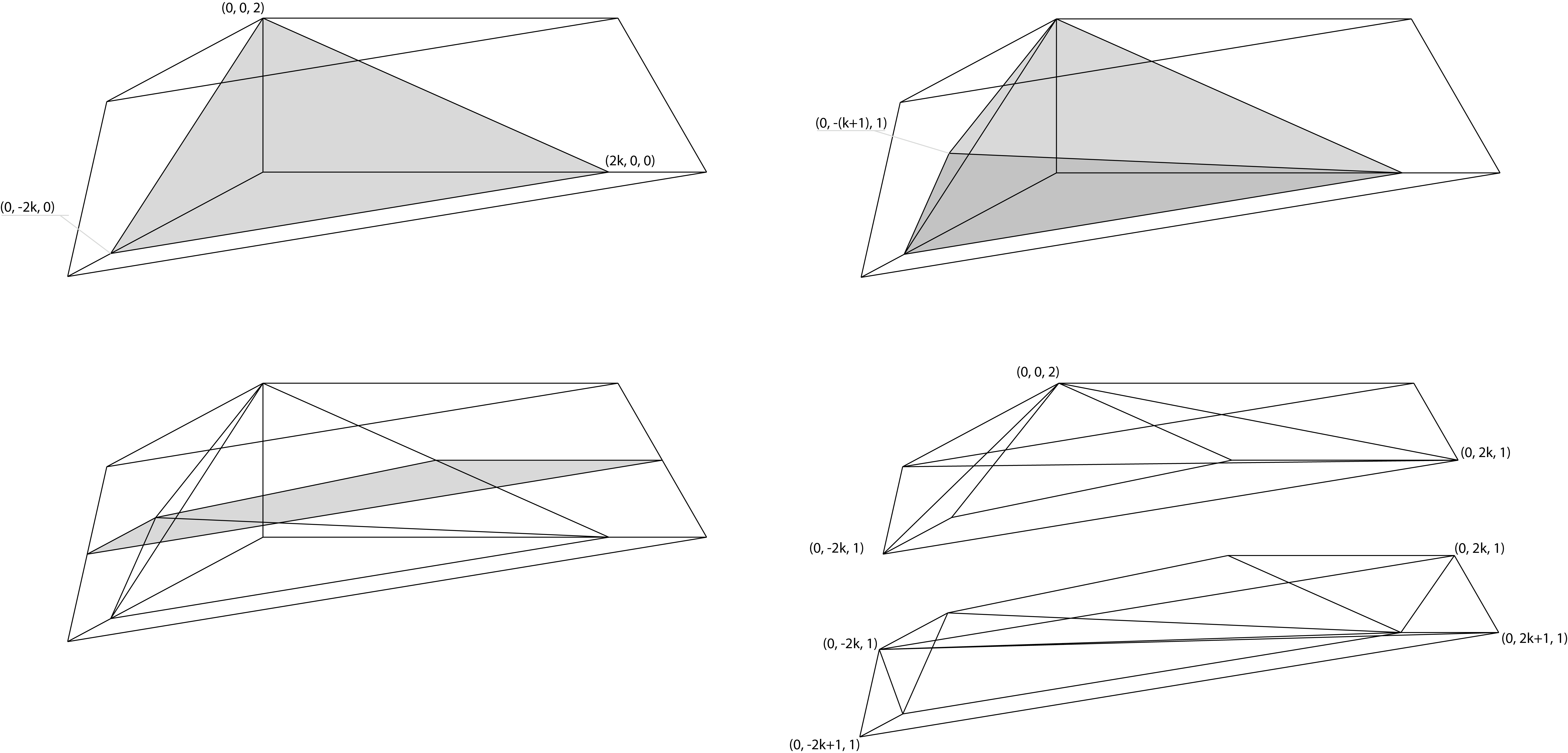}\\
    \begin{minipage}{16cm}
    \begin{it}
    \small{Fig. 13: Subdividing $C$ following step by the description in Section \ref{trimodif}. On the picture in the bottom right corner, $C$ is split along $C_1$ for readability purposes.}
    \end{it}
    \end{minipage}
\end{center}
\end{sloppypar}

\subsubsection{Sign distribution and topology of the patchworked surfaces}
\label{gudkov}
\begin{sloppypar}
We describe a sign distribution on the vertices of the triangulations described above. Let $(x_1,x_2,x_3)$ be a vertex of one these triangulations. \end{sloppypar}
\begin{sloppypar}
    If $x_3 \geq 2$, then $(x_1,x_2,x_3)$ gets the sign $+$ if it is of parity $(0,0,1)$ or $(1,0,0)$ and the sign $-$ otherwise.
\end{sloppypar}
\begin{sloppypar}
    If $x_3=1$, then $(x_1,x_2,x_3)$ receives the sign $+$ if $(x_1,x_2,x_3)$ does not belong to $C_1$, the sign $+$ if $(x_1,x_2,x_3)$ belongs to $C_1$ and $x_1$ and $x_2$ are both even and the sign $-$ if $(x_1,x_2,x_3)$ belongs to $C_1$ but $x_1$ and $x_2$ are not both even.
\end{sloppypar}
%If $z=1$, then $(x,y,z)$ receives
%\begin{itemize}
%    \item the sign $+$ if $(x,y,z)$ does not belong to $C_1$;
%    \item the sign $+$ if $(x,y,z)$ belongs to $C_1$ and $x$ and $y$ are both even;
%    \item the sign $-$ if $(x,y,z)$ belongs to $C_1$ but $x$ and $y$ are not both even.
%\end{itemize}
\begin{sloppypar}
    If $x_3=0$ and $(x_1,x_2,x_3)$ does not lie in the triangle with vertices $(0,0,0)$, $a$, $b$, then $(x_1,x_2,x_3)$ receives the sign $+$ if $x_1+x_2 \leq 4k+2$ and $x_1$ and $x_2$ are both even, the sign $+$ if $x_1+x_2 = 4k+5$ and $(x_1,x_2,x_3) \equiv (1,0,0)$ mod$2$ and the sign $-$ otherwise.
\end{sloppypar}
%If $z=0$ and $(x,y,z)$ does not lie in the triangle with vertices $(0,0,0)$, $a$, $b$, then $(x,y,z)$ receives
%\begin{itemize}
%    \item the sign $+$ if $x+y \leq 4k+2$ and $x$ and $y$ are both even;
%    \item the sign $+$ if $x+y = 4k+5$ and $(x,y,z) \equiv (1,0,0)$ mod$2$;
%    \item the sign $-$ otherwise.
%\end{itemize}
\begin{sloppypar}
    If $x_3=0$ and $(x_1,x_2,x_3)$ lies in the triangle with vertices $(0,0,0)$, $a$, $b$, then $(x_1,x_2,x_3)$ receives the sign $+$ if $(x_1,x_2,0) \equiv (0,1,0) ($mod $2)$, the sign $-$ otherwise.
\end{sloppypar}
%If $z=0$ and $(x,y,z)$ lies in the triangle with vertices $(0,0,0)$, $a$, $b$, then $(x,y,z)$ receives
%\begin{itemize}
%    \item the sign $+$ if $(x,y,0) \equiv (0,1,0) ($mod $2)$.
%    \item the sign $-$ otherwise.
%\end{itemize}
\begin{sloppypar}
For given vertices $a$ and $b$, we call \textit{the $(a,b)$-SD datum on $\Delta_m^3$} (where SD stands for small deviation) the datum of the signed triangulation described above. An example of curve obtained by restricting an $(a,b)$-SD datum to a polygon lying in the plane $\{x_3=0\}$ is shown on Figure $14$.
%If we apply the combinatorial patchworking theorem to the SD datum restricted to the triangle with vertices $(0,0,0)$, $(2k-2,0,0)$, $(0,2k-2,0)$, we obtain a curve $C_k$ of degree $4k+2$. A double covering $Y^k$ of $\mathbb{CP}^2$ ramified over $\mathbb{C}C_k$ (as defined in Section \ref{definitions}) verifies the equality $\chi(\mathbb{R}Y_+^k) = \sigma(\mathbb{C}Y^k) - 2 + \frac{(4k+2)^2}{2} = 0$. This is a heuristical justification for the following proposition.
\begin{prop}
\label{res1}
Let $a$ and $b$ be fixed triangulation vertices satisfying the conditions presented in Section \ref{trimodif} and let $p$ (respectively $n$) be the number of even integer points (respectively, of non-even integer points) lying in the relative interior of the triangle with vertices $(0,0,0)$, $a$, $b$. Applying the combinatorial patchworking theorem to the $(a,b)$-SD datum on $\Delta_{m=2k+1}^3$ produces a piecewise linear surface $\Tilde{\Gamma}_m^{(a,b)-SD}$ homeomorphic to the real part of a maximal non-singular real algebraic surface $X_m^{(a,b)-SD}$ of degree $m$ in $\mathbb{RP}^3$ and to the disjoint union of $\frac{m^3}{6}-m^2+11\frac{m}{6}-1 + (n-p)$ spheres and a projective plane with $\frac{m^3}{3}-m^2+\frac{7m-3}{6} - (n-p)$ handles attached to it. In particular, $\chi(\mathbb{R}X_m^{(a,b)-SD}) = \sigma(\mathbb{C}X_m^{(a,b)-SD}) + 4(n-p)$.
\end{prop}
\end{sloppypar}

\subsubsection{Counting cycle-axis pairs for an $(a,b)$-SD datum on $\Delta_m^3$}
\begin{sloppypar}
Given two integer points $a$ and $b$ satisfying the conditions presented in Section \ref{trimodif}, we describe cycles representing $\mathbb{Z}_2$-linearly independent homology classes of the surface $\Tilde{\Gamma}_m^{(a,b)-SD}$ to prove the maximality of $X_m^{(a,b)-SD}$. The approach is the same as in Section \ref{prelim}, but we have to consider the particularities of the $(a,b)$-SD triangulation of $\Delta_m^3 \bigcap \{ x_3 \leq 2 \}$.
\end{sloppypar}

\paragraph{{Cycle-axis pairs in the cone of vertex $(0,0,2)$ over the triangle with vertices $(0,0,0)$, $(2k,0,0)$, $(0,2k,0)$}}

\begin{sloppypar}
Let $T_1$ be the triangle with vertices $(0,0,0)$, $a$, $b$. The cone over this triangle is responsible for the Euler characteristic deviation. Let $v:=(x_1,x_2,0)$ be an integer point lying in the relative interior of $T_1$. The choice of sign distribution on $T_1$ ensures that there exists a symmetric copy $v'$ of $v$ such that the intersection of the star of $v'$ with $\{ x_3=0 \}$ contains a cell complex of the patchworked surface homemorphic to a circle. If $x_1$ and $x_2$ are not both even, then $v'$ bears the sign $+$ and the star of $v'$ contains a sphere belonging to the patchworked surface. Otherwise, $v'$ bears the sign $-$. In this case, we form two $(1,1)$-pairs associated to $v$. The cycle of the first pair is the cell complex homeomorphic to a circle contained in the intersection of the star of $v'$ with $\{ x_3=0 \}$. The dual axis is the suspension with vertices $(0,0,2)$ and $(0,0,-2)$ over the pair of points $(v', a)$. The axis of the second pair is the boundary of the intersection of the star of $v'$ with $\{ x_3=0 \}$. The dual cycle is formed as follows: there exists a simple rectilinear path $C$ contained in $\{x_3=0\} \bigcap \{ |x_1|+|x_2| \leq 2k-1 \}$ from $v'$ to $a$ that intersects the patchworked surface in exactly two points $e_1$ and $e_2$, where $e_1$ is contained in the star of $v'$. Remove from $C$ the segments $[v',e_1]$, $[e_2, a]$ and denote by $C'$ the resulting rectilinear path. The intersection of the suspension with vertices $(0,0,2)$ and $(0,0,-2)$ over $C$ and the patchworked surface is a double covering of $C'$ ramified at $e_1$ and $e_2$. This double covering is the desired $1$-cycle.
\end{sloppypar}
\begin{sloppypar}
Let $T_1^c$ be the closure of $(\{ x_3=0 \} \bigcap \{ x_1+x_2 \leq 2k \}) \ T_1$. Let $v:=(x_1,x_2,0)$ be an integer point lying in the relative interior of $T_1^c$. We can apply the arguments above, but the roles of the parities of the points are reversed. If $x_1$ and $x_2$ are both even, there is a sphere belonging to the patchworked surface in the star of a symmetric copy $v'$ of $v$. Otherwise, we can associate two $(1,1)$-pairs to $v$.
\end{sloppypar}
%\begin{sloppypar}
%Consider the boundary $D$ of the triangle $T_1$. In the intersection of $\{ z=0 \}$ with the union of the symmetric copies of the stars of the simplices belonging to $D$, there is a cell complex homeomorphic to a circle and belonging to the patchworked surface. Using the same description as above, we can exhibit two $(1,1)$-pairs associated to $D$.
%\end{sloppypar}
\begin{sloppypar}
In total, we have described $m^2-5m+6$ pairs in the symmetric copies of the cone of vertex $(0,0,2)$ over the triangle with vertices $(0,0,0)$, $(2k,0,0)$, $(0,2k,0)$. We call these pairs suspension pairs, like those described in Section \ref{surfsusp}.
\end{sloppypar}

\paragraph{Cycle-axis pairs associated to the two-dimensional walls with respective vertices $(2k,0,0)$, $(0,2k,0)$, $(0,k+1,1)$ and $(2k,0,0)$, $(0,2k,0)$, $(0,2k,1)$}

\label{otherwalls}

\begin{sloppypar}
The wall with vertices $(2k,0,0)$, $(0,2k,0)$, $(0,k+1,1)$ and the wall with vertices $(2k,0,0)$, $(0,2k,0)$, $(0,2k,1)$ can almost be treated as in Sections \ref{wallregu} and \ref{wallspecial}, using Lemma \ref{linsystem}. The only difference is the form of the axis of the special $(1,1)$-pair associated to each of these walls.
\end{sloppypar}
\begin{sloppypar}
Let $p$ be the integer point lying in the relative interior of $[(2k,0,0), (0,2k,0)]$ such that all the axes from every regular wall pairs whose essential part is a copy of $[v,(0,k+1,1)]$ where $v$ is an integer point lying in the relative interior of $[(2k,0,0), (0,2k,0)]$ and whose first two coordinates are not both even contain a copy of $[v,(0,k+1,1)]$. Lemma \ref{linsystem} ensures that there exists a unique orthant $O$ such that there is $1$-cycle in the star of the $O$-copy of $[p,(k+1,0,1)]$. The dual axis is the union of the $O$-copy of $[p,(0,k+1,1)]$, the $O$-copy of $[p,(2k+1,0,0)]$, the $O$-copy of $[(2k+1,0,0), (2k,0,1)]$ and a simple path from the $O$-copy of $(2k,0,1)$ to the $O$-copy of $(0,k+1,1)$, contained in the $O$-copy of $C_1$, formed of triangulation segments whose extremities bear the same sign.
%the following triangulation segments:
%\begin{itemize}
%    \item the $O$-copy of $[p,(0,k+1,1)]$;
%    \item the $O$-copy of $[p,(2k+1,0,0)]$;
%    \item the $O$-copy of $[(2k+1,0,0), (2k,0,1)]$;
%    \item a simple path from the $O$-copy of $(2k,0,1)$ to the $O$-copy of $(0,k+1,1)$, contained in the $O$-copy of $C_1$, formed of triangulation segments whose extremities bear the same sign.
%\end{itemize}
\end{sloppypar}
\begin{sloppypar}
Let $p'$ be the integer point lying in the relative interior of $[(2k,0,0), (0,2k,0)]$ such that all the axes from every regular wall pairs whose essential part is a copy of $[v,(0,2k,1)]$ where $v$ is an integer point lying in the relative interior of $[(2k,0,0), (0,2k,0)]$ and whose first two coordinates are not both even contain a copy of $[v,(0,2k,1)]$. Lemma \ref{linsystem} ensures that there exists a unique orthant $O'$ such that there is $1$-cycle in the star of the $O'$-copy of $[p',(0,2k,1)]$. The dual axis is the union of the $O'$-copy of $[p',(0,2k,1)]$, the $O'$-copy of $[p', (2k-1,0,0)]$, the $O'$-copy of $[(2k-1,0,0), (k,0,1)]$ and a simple path from the $O'$-copy $(k,0,1)$ to the $O'$-copy of $(0,2k,1)$, contained in the $O'$-copy of $C_1$, formed of triangulation segments whose extremities bear the same sign.
%the following triangulation segments:
%\begin{itemize}
%    \item the $O'$-copy of $[p',(0,2k,1)]$;
%    \item the $O'$-copy of $[p', (2k-1,0,0)]$;
%    \item the $O'$-copy of $[(2k-1,0,0), (k,0,1)]$;
%    \item a simple path from the $O'$-copy $(k,0,1)$ to the $O'$-copy of $(0,2k,1)$, contained in the $O'$-copy of $C_1$, formed of triangulation segments whose extremities bear the same sign.
%\end{itemize}
\end{sloppypar}
\begin{sloppypar}
We obtain two $(1,1)$-pairs for each integer point in the strongly inside the segment $[(2k,0,0),(0,2k,0)]$, that is $m-4$ pairs in total. Among those pairs, there are two special wall pairs.
\end{sloppypar}

\paragraph{Cycle-axis pairs in the star of the triangle with vertices  $(2k,0,0)$, $(0,2k,0)$, $(0,0,2)$}
\label{tampon}
\begin{sloppypar}
Consider the triangle $D$ with vertices $(2k,0,0)$, $(0,2k,0)$, $(0,0,2)$. There are two cones over $D$, with respective vertices $(2k-1,0,0)$ and $(0,k+1,1)$. Let $p$ be an integer point lying in the relative interior of $D$. In the intersection of the copies of the star of $p$ with the copies of $D$, there are two $1$-cycles. For one of them, we can find a suitable $1$-axis that contains the suspension over the corresponding copy of $p$ whose vertices are copies of the vertices of the cones over $D$ in the triangulation. The suspension over the other $1$-cycle is a sphere.
\end{sloppypar}
\begin{sloppypar}
In the star of each segment $S$ whose extremities are $(0,0,2)$ and a point on the segment $[(2k,0,0),(0,2k,0)]$, there is a $1$-cycle. A suitable $1$-axis can be constructed as the union of $4$ copies of $S$.
In total, we exhibit $2m-5$ pairs. We call these pairs suspension pairs.  
\end{sloppypar}

\paragraph{Cycle-axis pairs associated to one-dimensional walls between two-dimensional walls}
\label{1wall}

\begin{sloppypar}
Consider the segments $[(2k,0,0), (0,k+1,1)]$, $[(k,0,1), (0,k+1,1)]$ and $[(2k,0,0), (0,2k,1)]$.
In the case of $[(2k,0,0), (k+1,0,1)]$ and $[(k,0,1), (0,k+1,1)]$, analyzing the signs of the vertices in the union of the symmetric copies of their star, we can exhibit one $1$-cycle that is linked with a union of four symmetric copies of the segment whose endpoints bear the same sign. In the case of $[(2k,0,0), (0,2k,1)]$, one can construct a $1$-axis containing a copy of $[(2k,0,0), (0,2k,1)]$, a copy of $[(2k,0,0),(2k,0,1)]$, passing through one symmetric copy of the quadrangle $C_1$ and linked with a $1$-cycle which is not contractible in $\mathbb{RP}^3$ and which lies in the union of the symmetric copies of the star of $[(2k,0,0), (0,2k,1)]$. This provides three supplementary pairs that we call $1$-wall pairs.
\end{sloppypar}

\paragraph{Remaining parts of the triangulation}
\begin{sloppypar}
The suspension over the quadrangle $C_1$ can be treated in the same way as the suspensions from Section \ref{surfsusp} (using Lemma \ref{linsystem}). For each point in the relative interior of $C_1$, it is possible to exhibit two suspension $(1,1)$-pairs and one suspension $(2,0)$-pair. As in \ref{surfsusp}, we also obtain one special suspension $(1,1)$-pair. This amounts to $\frac{3m^2}{2}-8m+\frac{21}{2}+1$ pairs.   
\end{sloppypar}
\begin{sloppypar}
We can treat the join of the segments $[(2k,0,0), (0,2k,0)]$ and $[(0,k+1,1), (0,2k,1)]$ and the join of the segments $[(0,4k+4,0), (0,4k+5,0)]$ and $[(4k+4,0,1), (0,4k+4,1)]$ in the same way as the joins from Section \ref{join}, using Lemma \ref{linsystem}. We obtain $m^2-5m+6$ pairs and $\frac{m^2}{2}-\frac{7m}{2}+5$ pairs respectively. These pairs are called join pairs, like those described in Section \ref{join}.
\end{sloppypar}
\begin{sloppypar}
The following walls can be treated as in Sections \ref{wallregu} and \ref{wallspecial} (using Lemma \ref{linsystem}): the wall with vertices $(2k,0,1)$, $(0,2k,1)$, $(0,0,2)$, the wall with vertices $(0,0,2)$, $(0,2k-1,2)$, $(2k,0,1)$, the wall with vertices  $(0,k+1,1)$, $(0,2k,1)$, $(2k,0,0)$ and the wall with vertices $(2k,0,1)$, $(0,2k,1)$, $(2k,0,0)$. %For each segment lying in the relative interior of one of these walls, we exhibit one $(1,1)$-pair. 
In the star of the copies of the wall with vertices $(2k,0,1)$, $(0,2k,1)$, $(0,0,2)$, we exhibit $m-3$ pairs.
We exhibit $m-3$ pairs associated to segments inside the wall with vertices $(0,0,2)$, $(0,2k-1,2)$, $(2k,0,1)$. One of these pairs is a special wall pair, while the others are regular wall pairs.
We exhibit $\frac{m-1}{2}-2$ pairs associated to segments inside the wall with vertices $(0,k+1,1)$, $(0,2k,1)$, $(2k,0,0)$. One of these pairs is a special wall pair.
We exhibit $m-2$ pairs associated with segments inside the wall with vertices $(2k,0,1)$, $(0,2k,1)$, $(2k,0,0)$. One of these pairs is a special wall pair.
\end{sloppypar}
\begin{sloppypar}
The rest of the construction can be treated exactly as in Section \ref{prelim}, as it is the same as an IV construction for surfaces. We exhibit $m^3-8m^2+20m-12$ pairs. Among those pairs there are suspension pairs, special suspension pairs, join pairs, regular wall pairs, special wall pairs.   
\end{sloppypar}

\paragraph{Linking relations and partial order on the collection of $(1,1)$-pairs}
\label{ordermodif}
\begin{sloppypar}
We want to repeat the argument from Section \ref{indep} with the collections of $(1,1)$-pairs and $(2,0)$-pairs that we have described. It is necessary to order the collection of $(1,1)$-pairs. The order is essentially the same as the one described in Section \ref{order}. The $\mathbb{Z}_2$-linking number of the cycle and the axis from a pair is clearly $1$. The only other $\mathbb{Z}_2$-linking relations that are not mentioned in Section \ref{link} are the following ones.
\begin{itemize}
    \item The axes from the special wall pairs described in Section \ref{otherwalls} may be linked with some suspension pairs associated with the suspension over $C_1$. Indeed, these axes are completed using copies of triangulation segments that lie on $C_1$.
    \item For the same reason, the axis of the $1$-wall pair constructed starting with a copy of the segment $[(2k,0,0), (0,2k,1)]$ may be linked with some suspension pairs associated with the suspension over $C_1$.
    \item Some axes from regular wall pairs described in Section \ref{otherwalls} contain a copy of $[(2k,0,0), (0,k+1,1)]$, $[(k,0,1), (0,k+1,1)]$ or $[(2k,0,0), (0,2k,1)]$. Hence these wall axes may be linked with a cycle from a $1$-wall pair.
\end{itemize}
We consider the partial order generated by the following relations.
\begin{itemize}
    \item Let $(c,a)$ be a join pair and let $(c',a')$ be a regular wall pair. Then, we have the relation $(c,a)<(c',a')$.
    \item Let $(c,a)$ be a regular wall pair and let $(c',a')$ be a special wall pair. Then, we have the relation $(c,a)<(c',a')$.
    \item Let $(c,a)$ be a special wall pair and let $(c',a')$ be a $1$-wall pair. Then, we have the relation $(c,a)<(c',a')$.
    \item Let $(c,a)$ be a $1$-wall pair and let $(c',a')$ be a special suspension pair. Then, we have the relation $(c,a)<(c',a')$.
    \item Let $(c,a)$ be a special suspension pair. Let $(c',a')$ be a regular suspension pair. Then, we have the relation $(c,a)<(c',a')$.
\end{itemize}
\end{sloppypar}

\paragraph{Proof of Proposition \ref{res1}}
\begin{sloppypar}
Remove from the collection of $(1,1)$-pairs the $1$-wall pair whose cycle in not contractible in $\mathbb{RP}^3$. The class in $H_1(\Tilde{\Gamma}_m^{(a,b)-SD};\mathbb{Z}_2)$ represented by the cycle from this pair is independent from the classes represented by all the other described $1$-cycles, as all the other $1$-cycles are null-homologous in $\mathbb{RP}^3$.
As in Section \ref{indep}, we can then use the partial order defined in \mbox{Section \ref{ordermodif}} to prove the $\mathbb{Z}_2$-linear independence in $H_1(\Tilde{\Gamma}_m^{(a,b)-SD};\mathbb{Z}_2)$ of the homology classes represented by the other $1$-cycles we have described. The homology classes in $H_2(\Tilde{\Gamma}_m^{(a,b)-SD};\mathbb{Z}_2)$ represented by the described $2$-cycles are clearly independent, as the $2$-cycles represent the classes of different connected components of the patchworked surface. Hence, the exhibited cycles represent $m^3-4m^2+6m-3$ independent $\mathbb{Z}_2$-homology classes in $\Tilde{\Gamma}_m^{(a,b)-SD}$.
There is one more $2$-cycle representing a class that is independent from the others: the $2$-cycles exhibited above are homeomorphic to $2$-spheres but as the degree is odd, all connected components cannot be orientable (see \cite{wilson}). The Smith congruence (see Theorem \ref{smith}) then guarantees the existence of one more $1$-cycle independent from the others. Hence, the surface is maximal.\end{sloppypar}
\begin{sloppypar}
Recall that $p$ (respectively $n$) is the number of even (respectively non-even) integer points lying in the relative interior of the triangle with vertices $(0,0,0)$, $a$, $b$. The surface $\mathbb{R}X_m^{(a,b)-SD}$ has $n-p$ more spheres and $n-p$ less handles than the IV surface of same degree described in Section \ref{prelim}. The surfaces from Section \ref{prelim} are obtained using primitive triangulations, which implies that they are of type $\chi = \sigma$. Thus, the surface $X_m^{(a,b)-SD}$ we constructed in this section verifies $\chi(\mathbb{R}X_m^{(a,b)-SD}) = \sigma(\mathbb{C}X_m^{(a,b)-SD}) + 4(n-p)$. Notice that $n-p \geq 0$. \hfill$\qedsymbol$

\begin{center}
    \includegraphics[scale=0.75]{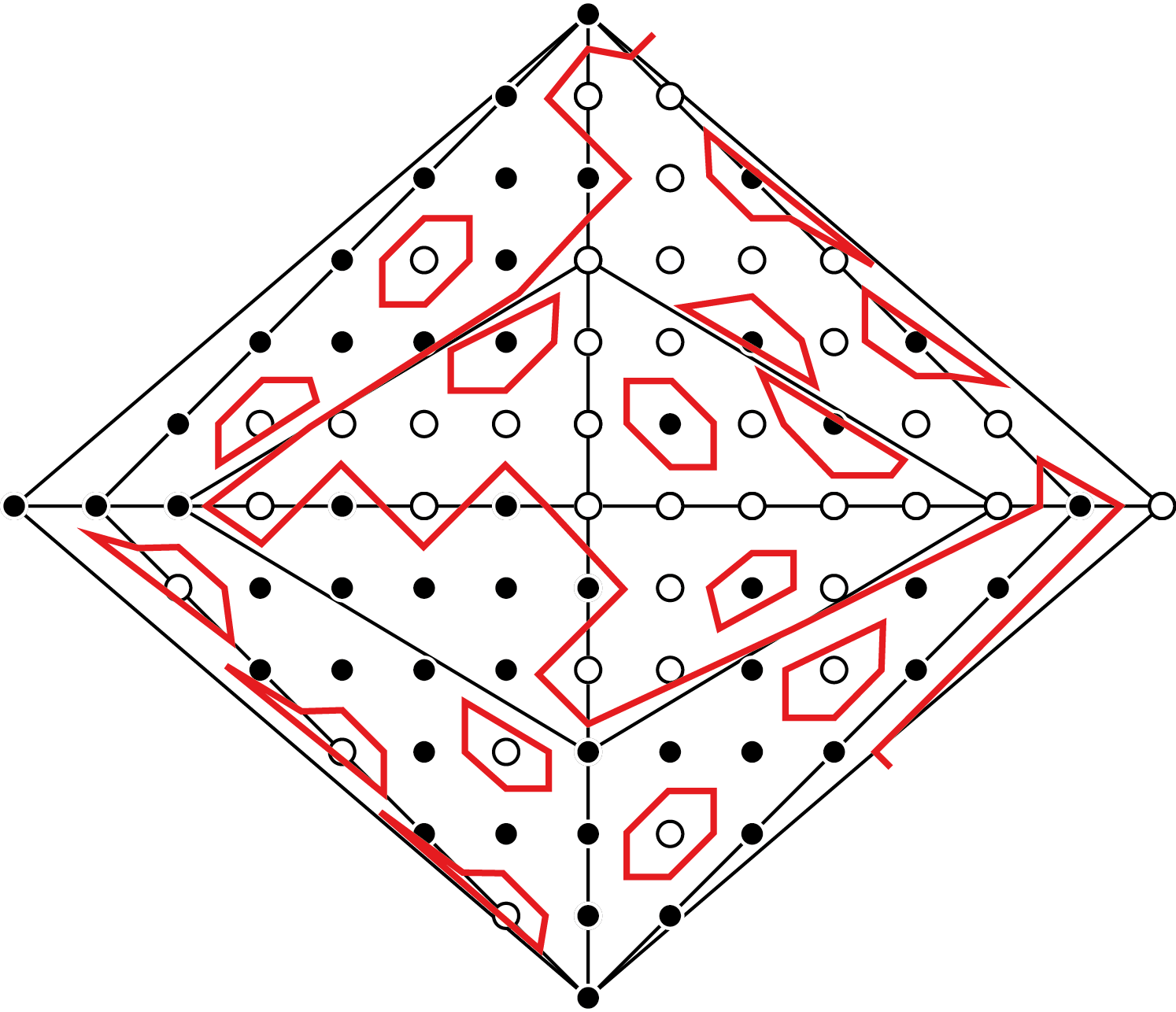}\\
    \begin{minipage}{16cm}
    \begin{it}
    \small{Fig. 14:  The curve obtained by restricting the $((5,0,0),(0,3,0))$-SD datum on $\Delta_9^3$ (here, $k=4$) to the triangle with vertices (0, 0, 0), (7, 0, 0), (0, 6, 0). White dots represent integer points bearing the sign $-$, black dots represent integer points bearing the sign $+$. Some important edges of the triangulation are drawn in black. The curve is sketched in red.}
    \end{it}
    \end{minipage}
\end{center}
For each cycle-axis pair exhibited above, either the cycle or the axis is connected when viewed in $(\Delta_m^3)^*$, without performing the identification of pairs of antipodal points lying on the boundary. This allows us to take the suspension over these cycles (or over their axis) later on, as we did in Section \ref{maximality3manif} for the even and the odd IV constructions.
\end{sloppypar}

\subsection{Propagating the Euler characteristic difference}

\begin{sloppypar}
Let $m$ be an even integer such that $m \geq 10$. In the odd Itenberg-Viro triangulation of $\Delta_m^4$, we replace the triangulation and sign distribution on each tetrahedron of the form \mbox{$\Delta_m^4 \bigcap \{ x_4= m-(2k+1) \}$} for some integer $k$ such that $k\geq 4$ and $m-(2k+1)>0$ by a certain SD datum on $\Delta_{2k+1}^3$. We prove that the resulting $3$-manifolds are maximal and have more connected components homeomorphic to $3$-spheres than the IV hypersurfaces of same degree. Furthermore, we prove that the modified hypersurfaces are not of type $\chi = \sigma$, which completes the proof of Theorem \ref{theorem1}.
\end{sloppypar}

\subsubsection{Triangulation, sign distribution and topology of the resulting $3$-manifold}
\label{construction3modif}
\begin{sloppypar}
Let $k$ be an integer such that $k \geq 4$ and $2k+1 < m$. On the triangle with vertices $(2k+1,0,0)$, $(0,2k+1,0)$, $(0,0,2k+1)$, any SD datum on $\Delta_{2k+1}^3$ coincides with the IV datum on $\Delta_{2k+1}^3$ used to triangulate $\Delta_m^4 \bigcap \{ x_4=m-(2k+1) \}$ in the odd IV construction for $3$-manifolds (see Section \ref{construction3}).
We triangulate $\Delta_m^4$ by following the same subdivision steps as for the odd IV triangulation with the following exception: for any positive integer $k$ such that $9 \leq 2k+1 \leq m-1$, we triangulate $\Delta_m^4 \bigcap \{ x_4=m-(2k+1) \}$ using the $( (2k-3,0,0), (0,2k-5,0))$-SD triangulation of $\Delta_{2k+1}^3$ instead of the IV triangulation of $\Delta_{2k+1}^3$. As $m \geq 10$, the resulting triangulation is different from the odd IV triangulation. Using Remark \ref{rmkconvexity}, the convexity of the resulting triangulation of $\Delta_m^4$, called \textit{the SD triangulation of} $\Delta_m^4$ for simplicity, follows from the convexity of the SD triangulation of $\Delta_{2k+1}^3$ for every $k$.\end{sloppypar}
\begin{sloppypar}
We define a sign distribution $\mu_{SD}$ on the vertices of the SD triangulation of $\Delta_m^4$ as follows. For any positive integer $k$ such that $9\leq 2k+1 \leq m-1$, the triangulation vertices lying in $\Delta_m^4 \bigcap \{ x_4=m-(2k+1) \}$ receive the signs given by the $((2k-3,0,0),(0,2k-5,0))$-SD sign distribution on the vertices of the SD triangulation of $\Delta_{2k+1}^3$. The remaining vertices receive the signs given by the odd IV sign distribution.
%\begin{itemize}
%    \item For any integer $k \geq 4$ such that $1\leq m-(2k+1) \leq m-1$, the triangulation vertices lying in $\Delta_m^4 \bigcap \{ t=m-(2k+1) \}$ receive the signs given by the $((2k-3,0,0),(0,2k-5,0))$-SD sign distribution on the vertices of the SD triangulation of $\Delta_{2k+1}^3$.
%    \item The remaining vertices receive the signs given by the odd IV sign distribution.
%\end{itemize}
By extension, we call this new triangulation of $\Delta_m^4$ and this new sign distribution \textit{the SD datum on $\Delta_m^4$}.
\label{maximalitytop}
\begin{prop}
The hypersurface $X_m^{SD}$ of even degree $m$ in $\mathbb{RP}^4$ obtained by applying the combinatorial patchworking theorem to the SD datum on $\Delta_m^4$ is maximal. This manifold has \mbox{$\frac{m^3}{24}-\frac{5m^2}{8}+\frac{17m}{6}-4+\frac{1}{4}\bar{m}$} more connected components homeomorphic to $3$-spheres than the odd IV hypersurface of the same degree, where $\bar{m}$ denotes the remainder of the Euclidean division of $m$ by $4$.
\end{prop}
\textbf{Proof:} We can count cycles representing linearly independent $\mathbb{Z}_2$-homology classes exactly as in Section \ref{maximality3manif}. The only thing that changes is the number of cycles of each dimension we recover in the symmetric copies of the suspensions over the tetrahedra of the form \mbox{$\Delta_m^4 \bigcap \{ x_4=m-(2k+1)  \}$} for some integer $k \geq 4$ such that $2k+1 \leq m-1$: there are $\sum_{k\geq4, 2k+1<m, k \: even} k^2 -4k + 4 + \sum_{k\geq4, 2k+1<m, k \: odd} k^2 -4k + 3$ more $3$-cycles homeomorphic to $3$-spheres and $\sum_{k\geq4, 2k+1<m, k \: even} k^2 -4k + 4 + \sum_{k\geq4, 2k+1<m, k \: odd} k^2 -4k + 3$ less $2$-cycles. The total number of cycles does not change, as we have exhibited the maximal number of cycles in the surfaces produced using a SD datum. Hence the resulting $3$-manifold is maximal. \hfill\qedsymbol
\end{sloppypar}

\subsubsection{Proof of Theorem \ref{theorem1}}

\begin{sloppypar}
\begin{prop}
The $3$-manifold $X_m^{SD}$ of degree $m$ obtained by applying the combinatorial patchworking theorem to the SD datum on $\Delta_m^4$ is not of type $\chi = \sigma$: it verifies the equality $\chi(\mathbb{R}Y^m_-) = \sigma(Y^m) -  8(\frac{m^3}{24}-\frac{5m^2}{8}+\frac{17m}{6}-4+\frac{1}{4}\bar{m})$,  where $Y^m$ is a double covering of $\mathbb{CP}^4$ ramified along $X_m^{SD}$.
\end{prop}
\textbf{Proof:} We repeat the computation of Section \ref{maximality3manif} with the SD datum and compare the result with what we obtained in \ref{maximality3manif} for the odd IV datum. In order to compare the number of copies with only vertices of sign $+$ of simplices in the SD and the odd IV construction, we notice that the odd IV datum and the SD datum are identical on the complementary of the suspensions over the tetrahedra of the form $(\Delta_m^4)^* \bigcap \{ x_4=m-(2k+1) \}$ for positive integers $k$ such that $9 \leq 2k+1 \leq m-1$. Moreover, the vertices of such suspensions are even integer points bearing the sign $-$. Hence, it suffices to count the number of simplices with only vertices of sign $+$ that are a copy of a simplex lying in a tetrahedron of the form $(\Delta_m^4)^* \bigcap \{ x_4=m-(2k+1) \}$ for some positive integer $k\geq4$ such that $2k+1 \leq m-1$ in the SD and in the odd IV construction.
Let $k$ be an integer such that $k \geq 4$ and $2k+1 \leq m-1$. The integer $m-(2k+1)$ is odd. Hence, two vertices lying in $(\Delta_m^4)^* \bigcap \{ x_4= \pm (m-2k+1) \}$ that are symmetric with respect to the coordinate hyperplane $\{ x_4=0 \}$ bear opposite signs. Thus, the problem reduces to counting the empty simplices (disregarding the signs) in the extended signed triangulations of $(\Delta_{2k+1}^3)^*$ induced by the $((2k-3,0,0),(0,2k-5,0))$-SD datum and IV datum on $\Delta_{2k+1}^3$ for each possible $k$.\end{sloppypar}
\begin{sloppypar}
%As we are computing an Euler characteristic, we do not need to know the number of simplices of each dimension, but only their alternate sum. This allows us to use the following trick.\\
The Euler characteristic of the real part of a real algebraic surface $X$ of degree $d$ produced by applying the combinatorial patchworking theorem to a triangulation $\tau$ and a sign distribution $\mu$ can be computed using the cell structure induced on $\mathbb{R}X$ by $\tau$ and $\mu$. We have \mbox{$\chi(\mathbb{R}X) = c_3 - c_2 + c_1$}, where $c_i$ is the number of non-empty $i$-simplices in the signed triangulation of $\widetilde{(\Delta_m^3)}^*$ induced by $\tau$ and $\mu$.
We want to compute $\bar{\chi}(\mathbb{R}X) = -\bar{c}_3 + \bar{c}_2 - \bar{c}_1 +\bar{c}_0$ where $\bar{c}_i$ is the number of empty simplices in the triangulation of $\Tilde{(\Delta_m^3)}^*$ induced by $\tau$. But $\bar{\chi}(\mathbb{R}X) =\chi(\mathbb{R}X)$. Indeed, \mbox{$\bar{\chi}(\mathbb{R}X) - \chi(\mathbb{R}X) = -\bar{c}_3 + \bar{c}_2 - \bar{c}_1 +\bar{c}_0 - c_3 + c_2 - c_1 = \chi(\mathbb{RP}^3) = 0$}.
A surface $X_{2k+1}^{IV}$ of degree $2k+1$ produced using the IV datum on $\Delta_{2k+1}^3$ is primitive. Hence $\chi(\mathbb{R}X_{2k+1}^{IV}) = \sigma(\mathbb{C}X_{2k+1}^{IV})$. From Proposition \ref{res1}, we know that a surface $X_{2k+1}^{((2k-3,0,0),(0,2k-5,0))-SD}$ of degree $2k+1$ produced using the $((2k-3,0,0),(0,2k-5,0))$-SD datum on $\Delta_{2k+1}^3$ satisfies $\chi(\mathbb{R}X_{2k+1}^{SD}) = \sigma(\mathbb{C}X_{2k+1}^{SD}) + 4k^2 - 16k + 16 - 4\bar{k}$, where $\bar{k}$ is the remainder of the Euclidean division of $k$ by $2$.
This implies that a hypersurface $X_m^{SD}$ of even degree $m$ in $\mathbb{RP}^4$ produced with the SD construction satisfies $\chi(\mathbb{R}Y^m_+) = \sigma(Y^m) -2 + 2\sum_{k,  9 \leq 2k+1 \leq m-1} 4k^2 - 16k + 16 - 4 \bar{k} = \sigma(Y^m) -2 + 8(\frac{m^3}{24}-\frac{5m^2}{8}+\frac{17m}{6}-4+\frac{1}{4} \bar{m})$ where $Y^m$ is a double covering of $\mathbb{CP}^4$ branched along $X_m^{SD}$, \textit{i.e.} $\chi(\mathbb{R}Y^m_-) = \sigma(Y^m) -  8(\frac{m^3}{24}-\frac{5m^2}{8}+\frac{17m}{6}-4+\frac{1}{4}\bar{m})$. If $m>10$, this means that $\chi(\mathbb{R}Y^m_-) < \sigma(Y^m)$. This completes the proof of Theorem \ref{theorem1}. \hfill\qedsymbol
\end{sloppypar}

\section{Stirring asymptotically away from $\chi = \sigma$}
\label{asymptoticalmodification}
\begin{sloppypar}
    In this last section, we modify the even IV construction to construct a family of maximal hypersurfaces that is not asymptotically of type $\chi = \sigma$ in order to prove \hyperref[theorem2]{Theorem 1.4}.
\end{sloppypar}

\subsection{Asymptotical deviation construction and Betti numbers of the AD hypersurfaces}

\subsubsection{Triangulation}

\label{triAD}

\begin{sloppypar}
    Let $m$ be an even integer such that $m \geq 8$. We retain the notations associated to the even IV triangulation from Section \ref{construction3}. We triangulate $\Delta_m^4$ following the subdivision steps that produce the even IV triangulation, but we use different triangulations of some of the triangles labeled $1$. In the even IV triangulation, each triangle labeled $1$ is sliced by the segments parallel to its edge of type $(11, 12)$. The resulting quadrangles are subdivided using one of their diagonals. The triangulation is then completed into a primitive one in the only possible way. Let $k$ be an even integer such that $8 \leq k \leq m$. We introduce different triangulations of the triangle with vertices $(k,0,0,m-k)$, $(0,k,0,m-k)$, $(0,0,k, m-k)$.
    \begin{itemize}
        \item Subdivide the triangle using the segments \mbox{$[(k-1,0,1,m-k), (0,k,0,m-k)]$}, \mbox{$[(1,0,k-1,m-k), (0,k,0,m-k)]$}.
        %$[(k-1,0,1,m-k), (2,k-5,3,m-k)]$, $[(1,0,k-1,m-k), (2,k-5,3,m-k)]$.
        \item Choose three distinct odd integer vertices $a$, $b$, $c$ such that $a$ and $b$ lie in the relative interior of $[(1,0,k-1,m-k), (k-1,0,1,m-k)]$, the vertex $c$ lies in the interior of the triangle with vertices $(k,0,0,m-k)$, $(0,k,0,m-k)$, $(0,0,k, m-k)$ and $[a,c]$ and $[b,c]$ are primitive. Subdivide the triangle using the segments $[a,c]$ and $[b,c]$.
        \item Complete this subdivision into a primitive convex triangulation.
    \end{itemize}
    We consider a triangulation of $\Delta_m^4$ induced by such a modification. It is characterized by a list $L$ of $\frac{m}{2}-3$ triplets of points $L_8$, $L_{10}$, ...,$L_m$ where $L_k$ is the triplet of points $a$, $b$, $c$ chosen to perform the subdivision of the triangle with vertices $(k,0,0,m-k)$, $(0,k,0,m-k)$, $(0,0,k, m-k)$ as described above. We call this triangulation \textit{the $L$-AD triangulation of $\Delta_m^4$} (where AD stands for asymptotical deviation). The $L$-AD triangulation of $\Delta_m^4$ is different from the even IV triangulation of $\Delta_m^4$. The $L$-AD triangulation is primitive, as the completion chosen above is primitive. As the completion is also convex, the $L$-AD triangulation of $\Delta_m^4 \bigcap \{x_4=m-k \}$ where $k$ is an even positive integer gives a triangulation of $\Delta_3^k$ which corresponds to the description from \mbox{\cite[Section~5]{tsurf}}. This ensures that there exists a family of convex functions certifying the convexity of the considered triangulations of $\Delta_k^3$ for any even integer $k$ such that $0 \leq k \leq m$. This family and Remark \ref{rmkconvexity} ensure the convexity of the $L$-AD triangulation of $\Delta_m^4$.
\end{sloppypar}

\subsubsection{Sign distribution}

\begin{sloppypar}
    We fix a list $L = \{ (a_k, b_k, c_k) \}_{8\leq k\leq m,\:k\:even}$ of triplets of points satisfying the conditions described in Section \ref{triAD} and we define a sign distribution on the vertices of the $L$-AD triangulation of $\Delta_m^4$ as follows.
    The vertices that do not belong to a triangle labeled $1$ with edges of length at least $8$ get the same signs as in the even IV construction.
    Let $k$ be an even integer such that $8 \leq k \leq m$. We denote by $T_k$ the triangle with vertices $(k,0,0,m-k)$, $(0,k,0,m-k)$, $(0,0,k, m-k)$. We describe the sign distribution on the vertices lying in $T_k$. 
    \begin{itemize}
        \item The vertices of the triangulation lying on the segment  $[(0,k,0,m-k), (0,0,k,m-k)]$ and of parity $(0,1,1,0)$ receive the sign $+$.
        \item The vertices of the triangulation of parity $(1,1,0,0)$ lying in the relative interior of $T_k$ or on the segments $[(k,0,0,m-k), (0,k,0,m- k)]$ and $[(k,0,0,m-k), (0,0,k,m-k)]$ but not in the triangle with vertices $a_k$, $b_k$, $c_k$ receive the sign $+$.
        \item The vertices of the triangulation of $T_k$ lying in the triangle with vertices $a_k$, $b_k$, $c_k$ and of parity $(0,0,0,0)$ receive the sign $+$.
        \item The other vertices of the triangulation of $T_k$ receive the sign $-$.
    \end{itemize}
    %\begin{itemize}
    %    \item The vertices that do not belong to a triangle labeled $1$ with edges of length at least $8$ get the same signs as in the even IV construction.
    %    \item Let $k$ be an even integer such that $8 \leq k \geq m$. We denote by $T_k$ the triangle with vertices $(k,0,0,m-k)$, $(0,k,0,m-k)$, $(0,0,k, m-k)$. We describe the sign distribution on the vertices lying in $T_k$. 
    %    \begin{itemize}
    %        \item The vertices of the triangulation lying on the segment  $[(0,k,0,m-k), (0,0,k,m-k)]$ and of parity $(0,1,1,0)$ receive the sign $+$.
    %        \item The vertices of the triangulation of parity $(1,1,0,0)$ lying in the relative interior of $T_k$ or on the segments $[(k,0,0,m-k), (0,k,0,m- k)]$ and $[(k,0,0,m-k), (0,0,k,m-k)]$ but not in the triangle with vertices $a_k$, $b_k$, $c_k$ receive the sign $+$.
    %        \item The vertices of the triangulation of $T_k$ lying in the triangle with vertices $a_k$, $b_k$, $c_k$ and of parity $(0,0,0,0)$ receive the sign $+$.
    %        \item The other vertices of the triangulation of $T_k$ receive the sign $-$.
    %    \end{itemize}
    %\end{itemize}
    \end{sloppypar}
    \begin{center}
    \includegraphics[scale=0.6]{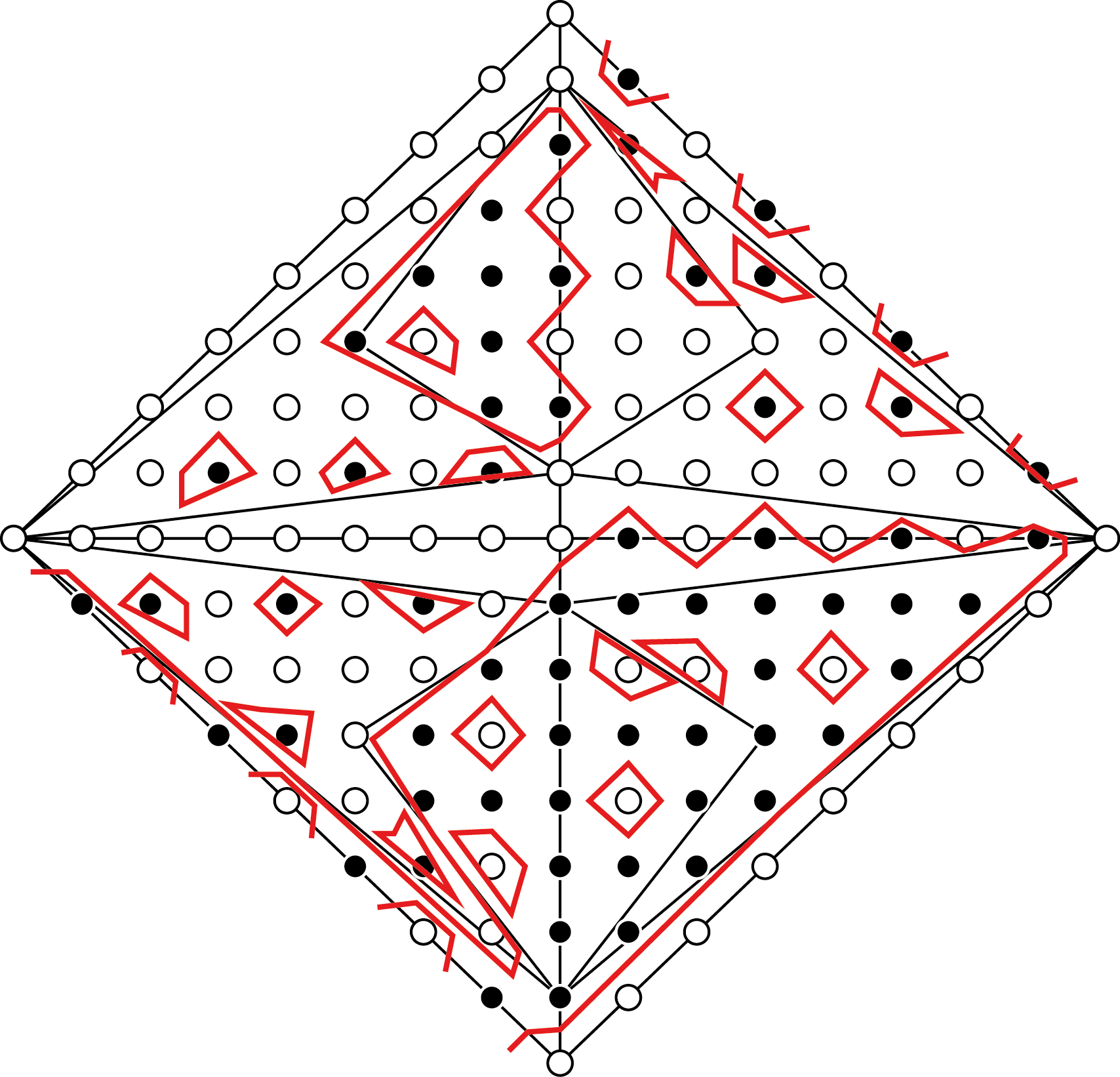}\\
    \begin{minipage}{16cm}
    \begin{it}
    \small{Fig. 15: The curve associated to a certain $L$-AD datum restricted to a triangle of a triangle labeled $1$ whose edges are of lattice length $8$. Here, we have $a_8=(7,0,1,m-8)$, $b_8=(1,0,7,m-8)$, $c_8=(2,3,3,m-8)$. We have taken a copy in the plane $\mathbb{R}^2$ of the triangle. Black dots represent integer points bearing the sign $+$. White dots represent integer points bearing the sign $-$. The curve appears in red. It is a maximal curve, as there are $22$ connected components.}
    \end{it}
    \end{minipage}
\end{center}
    \begin{sloppypar}
    We call the \textit{$L$-AD datum on $\Delta_m^4$} the datum composed of the $L$-AD triangulation and the above described sign distribution. The hypersurface of degree $m$ in $\mathbb{RP}^4$ obtained by applying the combinatorial patchworking theorem to $\Delta_m^4$ endowed with the $L$-AD datum is called \textit{the $L$-AD hypersurface of degree} $m$.
    The curve associated to a certain $L$-AD datum restricted to a triangle labeled $1$ whose edges are of lattice length $8$ is shown on Figure $15$.
    
\end{sloppypar}

\subsubsection{Betti numbers of the $L$-AD hypersurfaces}
\begin{sloppypar}
    We begin the proof of Theorem \ref{theorem2} by the following proposition.
    \begin{prop}
        \label{maxAD}
        Let $m$ be an even integer. Let $L=\{ (a_k,b_k,c_k) \}_{8 \leq k \leq m, \: k \: even}$ be a list of triplets of points as described in Section \ref{triAD}. The $L$-AD hypersurface of degree $m$ in $\mathbb{RP}^4$ is maximal and has the same $\mathbb{Z}_2$-Betti numbers as the even IV hypersurface of the same degree.
    \end{prop}
    \textbf{Proof:} We repeat the computation and the arguments from Section \ref{maximality3manif} with the $L$-AD datum on $\Delta_m^4$. There is only one difference: for each even integer $k$ such that $8 \leq k \leq m$, one cannot find four cell complexes homeomorphic to circles, contained in the intersection of the copies of $T_k$ with the copies of the star of $c_k$ and with the patchworked hypersurface. Instead, the curve associated to the AD datum restricted to the triangle $T_k$ with vertices $(k,0,0,m-k)$, $(0,k,0,m-k)$, $(0,0,k,m-k)$ has four connected components that are contained in the union of exactly two copies of $T_k$ sharing a copy of an edge of type $(10,12)$. These connected components are said to be \textit{large circles} and appear in the count only when dealing with parts of type $(00,10,11,12,22)$ and specifically when forming a $(2,1)$-pair whose axis contains a copy of the segment whose extremities are the point $c_k$ and a point of type $00$ (respectively, of type $22$). The union of two parts of type $(00,10,11,12,22)$ sharing a tetrahedron of type $(00,10,11,12)$ (respectively, of type $(10,11,12,22)$) can be viewed as the union of two joins of a triangle of type $(10,11,12)$ with a segment of type $(00,22)$. The $2$-cycles described in Section \ref{0010111222} and contained in the union of the symmetric copies of such a union can be taken to be the join of a circle that is the intersection of the patchworked hypersurface with some symmetric copy of the star of an integer point lying in the relative interior of the triangle of type $(10,11,12)$ and of the pair of points formed by the middle points of the copies of the two segments of type $(00,22)$ lying in the same orthant. This can also be done starting with a large circle. Indeed, the copies of the segments of type $(00,22)$ lie in the same coordinate hyperplane as the common face (of type $(10,12)$) of the pair of copies of $T_k$ that contain a fixed large circle. Let $Q$ be a triangulation $4$-simplex such that one of the segments of type $(00,22)$ is an edge of $Q$, the segment $[a_k, c_k]$ is an edge of $Q$ and the last vertex of $Q$ lies on $T_k$ but not in the relative interior of the triangle with vertices $a_k$, $b_k$, $c_k$.
    %\begin{itemize}
    %    \item one of the segments of type $(00,22)$ is an edge of $Q$;
    %    \item the segment $[a_k, c_k]$ is an edge of $Q$;
    %    \item the last vertex of $Q$ lies on $T_k$ but not in the relative interior of the triangle with vertices $a_k$, $b_k$, $c_k$.
    %\end{itemize}
    Applying \mbox{Lemma \ref{linsystem}} to $Q$, we obtain that there is a unique orthant $O$ such that there is a large circle partially contained in the $O$-copy of $T_k$, this large circle intersects $O$-copy of the intersection of star of $c_k$ with $T_k$, the middle points of the $O$-copy of the segments of type $(00,22)$ belong to the patchworked hypersurface and the signs of the $O$-copies of $c_k$ and the point of type $00$ (respectively $22$) are the same.
    %\begin{itemize}
    %    \item there is a large circle contained in the $O$-copy of $T_k$ and this large circle intersects $O$-copy of the intersection of star of $c_k$ with $T_k$;
    %    \item the middle points of the $O$-copy of the segments of type $(00,22)$ belong to the patchworked hypersurface;
    %    \item the signs of the $O$-copies of $c_k$ and the point of type $00$ (respectively $22$) are the same.
    %\end{itemize}
    Hence, the join of the considered large circle and the middle points of the $O$-copies of the segment of type $(00,22)$ is a $2$-cycle of the patchworked hypersurface. The dual axis contains the $O$-copy of the segments with extremities $c_k$ and the point of type $00$ (respectively of type $22$) and is formed as in Section \ref{0010111222}. The rest of the proof is unchanged and we can conclude that the $L$-AD hypersurface of degree $m$ in $\mathbb{RP}^4$ is maximal and has the same $\mathbb{Z}_2$-Betti numbers as the even IV hypersurface of the same degree. \hfill \qedsymbol
\end{sloppypar}
%\begin{sloppypar}
%    Notice that the $L$-AD hypersurface of degree $m$ has the same $\mathbb{Z}_2$-Betti numbers as the even IV hypersurface of degree $m$.
%\end{sloppypar}

\subsection{Euler characteristic computations for $L$-AD hypersurfaces}
\label{eulercharcompLAD}
\begin{sloppypar}
Let $m$ be an even positive integer such that $m \geq 8$. Let $L=\{ (a_k,b_k,c_k) \}_{8 \leq k \leq m, \: k \: even}$ be a list of triplets of points as described in Section \ref{triAD}. We consider the even IV datum and the $L$-AD datum on $\Delta_m^4$. In order to prove Theorem \ref{theorem2}, we need to compare the quantities \mbox{$\chi_{L-AD,m}^+ := c_{L-AD,m}^4 - c_{L-AD,m}^3 + c_{L-AD,m}^2 - c_{L-AD,m}^1 + c_{L-AD,m}^0$} and \mbox{$\chi_{IV,m}^+ := c_{IV,m}^4 - c_{IV,m}^3 + c_{IV,m}^2 - c_{IV,m}^1 + c_{IV,m}^0$}, where $c_{ST,m}^i$ is the number of simplices of dimension $i$ and with only vertices bearing the sign $+$ in the signed triangulation of $\widetilde{(\Delta_m^4)^*}$ induced by the ST datum, where ST $\in$ \mbox{$\{L$-AD, IV$\}$}, for each $i \in \{0,1,2,3,4 \}$. As we have proved in Section \ref{prooftype} that the IV-manifolds are of type $\chi = \sigma$, proving \hyperref[theorem2]{Theorem 1.4} is equivalent to finding a family of lists $\{L_m\}_{m \geq 8}$ such that $\chi_{L_m - AD,m}^+ \sim_{\infty} \chi_{IV,m}^+ + am^4$, where $a$ is a non-zero real number.
\end{sloppypar}

\subsubsection{Preliminaries : studying the $L$-AD and IV data in the interior of the triangles $T_k$}
\label{EulerTriangles}
\begin{sloppypar}
    Let $k$ be an even integer such that $8 \leq k \leq m$. Recall that $T_k$ is the triangle with vertices $(k,0,0,m-k)$, $(0,k,0,m-k)$, $(0,0,k,m-k)$. The $L$-AD and the even IV constructions differ only in the star of such triangles. We begin by comparing the restrictions of the $L$-AD datum and of the even IV datum to the relative interior of $T_k$. Consider a polygon $P$ contained in $T_k$ and whose vertices are integer points. Endow the polygon $P$ with a signed primitive triangulation $ST$. We denote by
    \begin{itemize}
        \item $\mathcal{A}(P)$ the lattice volume of $P$, \textit{i.e.} the euclidean area of $P$ divided by the euclidean area of the quadrangle with vertices $(k,0,0,m-k)$, $(k-1,1,0,m-k)$, $(k-1,0,1,m-k)$, $(k-2,1,1,m-k)$;
        \item $\Gamma_{P,ST}$ the collection of simplices of the triangulation whose relative interior is contained in the interior of $P$;
        \item $t_P$ the number of triangles in $\Gamma_{P,ST}$; by $t_{adj,P,ST}^+$ the number of triangles of $\Gamma_{P,ST}$ adjacent to an even vertex bearing the sign $+$; by $t_{even,P,ST}^+$ the number of even triangles in $\Gamma_{P,ST}$ whose empty copies only have vertices of sign $+$;
        \item $e_P$ the number of edges in $\Gamma_{P,ST}$; by $e_{adj,P,ST}$ the number of edges of $\Gamma_{P,ST}$ adjacent to an even vertex; by $e_{adj,P,ST}^+$ the number of edges of $\Gamma_{P,ST}$ adjacent to an even vertex bearing the sign $+$; by $e_{nadj,P,ST}$ the number of edges of $\Gamma_{P,ST}$ that are not adjacent to an even vertex;
        \item $v_P$ the number of vertices in $\Gamma_{P,ST}$; by $v_{even,P}$ the number of even vertices in $\Gamma_{P,ST}$; by $v_{odd,P}$ the number of odd vertices in $\Gamma_{P,ST}$; by $v_{bdry,P}$ the number of even vertices lying on the boundary of $P$; by $v_{even,P,ST}^+$ the number of even vertices in $\Gamma_{P,ST}$ bearing the sign $+$; by $v_{bdry,P,ST}^+$ the number of even vertices bearing the sign $+$ and lying on the boundary of $P$;
        \item $\chi_{P,ST}^+$ the quantity $t_{even,P,ST}^+ + t_{adj,P,ST}^+ -2e_{adj,P,ST}^+ + 4v_{even,P,ST}^+ - e_{nadj,P,ST} + 2v_{odd,P}$.
    \end{itemize}
    As the considered triangulation is primitive, the quantities $t_P$, $e_P$, $v_P$, $v_{even,P}$, $v_{odd,P}$, $v_{bdry,P}$ do not depend on $ST$. Moreover, we have \mbox{$t_P = 2\mathcal{A}(P)$} and \mbox{$t_P+v_P-e_P = 1$}, as well as \mbox{$e_{adj,P,ST}^+ = t_{adj,P,ST}^+ - v_{bdry,P,ST}^+$}.
    %\begin{itemize}
    %    \item $t_P = 2\mathcal{A}(P)$;
    %    \item $t_P+v_P-e_P = 1$;
    %    \item $e_{adj,P,ST}^+ = t_{adj,P,ST}^+ - v_{bdry,P,ST}^+$.
    %\end{itemize}
    The last equalities follow from the Euler formula applied to $P$ and to the union of the stars of the even points in $P$.
    In particular, when $P = T_k$, we obtain
    \begin{center}
        $\mathcal{A}(T_k) = \frac{k^2}{2}$; \hspace{20pt} $t_{T_k}=k^2$; $e_{T_k}=\frac{3k^2}{2} - \frac{3k}{2}$;\\ $v_{even,T_k}=\frac{k^2}{8} - \frac{3k}{4}+1$; \hspace{20pt} $v_{odd,T_k} = \frac{3k^2}{8}-\frac{3k}{4}$; $v_{bdry,T_k} = \frac{3k}{2}$.
    \end{center}
    %\begin{itemize}
    %    \item $\mathcal{A}(T_k) = \frac{k^2}{2}$;
    %    \item $t_{T_k}=k^2$;
    %    \item $e_{T_k}=\frac{3k^2}{2} - \frac{3k}{2}$;
    %    \item $v_{even,T_k}=\frac{k^2}{8} - \frac{3k}{4}+1$;
    %    \item $v_{odd,T_k} = \frac{3k^2}{8}-\frac{3k}{4}$;
    %    \item $v_{bdry,T_k} = \frac{3k}{2}$.
    %\end{itemize}
    First, we consider $P = T_k$ endowed with the restriction of the even IV datum.
    \begin{lem}
        \label{contributionIV}
        Let $m$ and $k$ be even integers such that $8 \leq k \leq m$. Then, one has $\chi_{T_k,IV}^+ = \frac{k^2}{4} - \frac{3k}{2}$.
    \end{lem}
    \textbf{Proof:} In the even IV construction, all the triangles of the triangulation of $T_k$ are adjacent to an even vertex and all the even vertices of the triangulation of $T_k$ bear the sign $-$. Hence $e_{nadj,T_k,IV} = e - e_{adj,T_k, IV} = \frac{3k^2}{2} - \frac{3k}{2} - k^2 + 3k$. Indeed, as the the star of any even point $v$ is homeomorphic to a disk, there are as many edges as triangles in such a star when $v$ lies in the relative interior of $T_k$ and as many edges as the number of triangles $-1$ when $v$ lies on the boundary of $T^k$.
    Now, we have:
    \begin{center}
        $\chi_{T_k,IV}^+ = -e_{nadj,T_k,IV} + 2v_{odd,T_k} = -(\frac{3k^2}{2} - \frac{3k}{2} - k^2 + \frac{3k}{2}) + \frac{3k^2}{8}-\frac{3k}{4} = \frac{k^2}{4} - \frac{3k}{2}.$
    \end{center}
    \hfill \qedsymbol
    \end{sloppypar}
    \begin{sloppypar}
    In a $L$-AD triangulation of $T_k$, the number of even triangles and the number of edges that are not adjacent to an even point can vary. Plus, some of the even vertices lying on $T_k$ bear the sign $+$. We treat separately the triangle $T_k^1$ with vertices $(1,0,k-1,m-k)$, $(0,0,k,m-k)$, $(0,k,0,m-k)$, the triangle $T_k^2$ with vertices $a_k$, $b_k$, $c_k$ and $T_k \setminus (T_k^1 \bigcup T_k^2)$.
    %the three following zones of the triangulation of $T_k$:
    %\begin{itemize}
    %    \item the triangle $T_k^1$ with vertices $(1,0,k-1,m-k)$, $(0,0,k,m-k)$, $(0,k,0,m-k)$;
    %    \item the triangle $T_k^2$ with vertices $a_k$, $b_k$, $c_k$;
    %    \item $T_k \setminus (T_k^1 \bigcup T_k^2)$.
    %\end{itemize}
    To simplify the computations, we borrow the notion of Harnack zone from \cite{haas}.
    \begin{defin}
        Assume that $T_k$ is equipped with a signed primitive triangulation $ST$. A lattice polygon $P \subset T_k$ is called a \textit{Harnack zone} of $T_k$ if the edges of $P$ are edges of the triangulation of $T_k$ and the triangulation vertices lying in $P$ of three different parities bear the same sign while the triangulation vertices lying in $P$ of the last parity bear the opposite sign.
        %it satisfies the following conditions:
        %\begin{itemize}
        %    \item the edges of $P$ are edges of the triangulation of $T_k$;
        %    \item the triangulation vertices lying in $P$ of three different parities bear the same sign; the triangulation vertices lying in $P$ of the last parity bear the opposite sign.
        %\end{itemize}
    \end{defin}
    Notice that $T_k^1$, $T_k^2$ and $T_k \setminus (T_k^1 \bigcup T_k^2)$ are Harnack zones of $T_k$ endowed with the restriction of the $L$-AD datum.
    \begin{lem}
    \label{lemharnack}
        Let $P$ be a Harnack zone of $T_k$ endowed with a signed triangulation $ST$.\\
        If the even vertices lying in $P$ bear the sign $+$, then $\chi_{P,AD}^+ = -2\mathcal{A}(P) + 3 v_{even,P} + v_{odd,P} + v_{bdry,P} + 1$.\\
        If the even vertices lying in $P$ bear the sign $-$, then $\chi_{P,AD}^+ = -v_{even,P} + v_{odd,P} - v_{bdry,P} + 1$.
    \end{lem}
    \textbf{Proof:}
    As the triangulation of $P$ is primitive, one has $t_{P} = 2\mathcal{A}(P)$.
    Suppose that the even vertices lying in $P$ bear the sign $+$. Notice that this means that the empty symmetric copies of the even triangles of the triangulation of $P$ only have vertices of sign $-$. Thus, we have \mbox{$\chi_{P,AD}^+ = t_{adj,P,AD} - 2e_{adj,P,AD} - e_{nadj,P,AD} + 4v_{even,P} + 2v_{odd,P}$}. We also have  \mbox{$t_{adj,P,ST} = 2\mathcal{A}(P) - t_{even,P,ST}$} and \mbox{$e_{adj,P,ST} = t_{adj,P,ST} - v_{bdry,P} = 2\mathcal{A}(P) - t_{even, P,AD} -v_{bdry,P}$}, as well as \mbox{$e_{nadj,P,ST} = 2\mathcal{A}(P) + v_{even,P} + v_{odd,P} - 1 - e_{adj,P} = v_{even,P} + v_{odd,P} - 1 + t_{even, P,ST} + v_{bdry,P}$}.
    %\begin{itemize}
    %    \item $t_{adj,P,ST} = 2\mathcal{A}(P) - t_{even,P,ST}$;
    %    \item $e_{adj,P,ST} = t_{adj,P,ST} - v_{bdry,P} = 2\mathcal{A}(P) - t_{even, P,AD} -v_{bdry,P}$;
    %    \item $e_{nadj,P,ST} = 2\mathcal{A}(P) + v_{even,P} + v_{odd,P} - 1 - e_{adj,P} = v_{even,P} + v_{odd,P} - 1 + t_{even, P,ST} + v_{bdry,P}$.
    %\end{itemize}
    %The second relation is again an application of the Euler formula to the union of stars of even points. Each such star is homeomorphic to a disk.\\
    Thus, we obtain
    \begin{align*}
        \chi_{P,ST}^+ &= t_{adj,P,ST} - 2e_{adj,P,ST} - e_{nadj,P,ST} + 4v_{even,P} + 2v_{odd,P}\\
        & = 2\mathcal{A}(P) - t_{even, P, ST} - 2(2\mathcal{A}(P) - t_{even,P,ST} -v_{bdry,P}) \\
        & \;\;\;\;\;  - (v_{even,P} + v_{odd,P} - 1 + t_{even, P, ST} + v_{bdry,P}) + 4v_{even,P} + 2v_{odd,P}\\
        &= -2\mathcal{A}(P) + 3v_{even,P} + v_{odd,P} + v_{bdry,P} + 1.
    \end{align*}
    Now suppose that the even vertices lying in $P$ bear the sign $-$. Notice that this means that the empty symmetric copies of the even triangles of the triangulation of $P$ only have vertices of sign $+$. We then have
    \begin{align*}
        \chi_{P,ST}^+ &=  t_{even,P,ST} - e_{nadj,P,ST} + 2v_{odd,P}\\
        & = t_{even,P,ST} - (v_{even,P} + v_{odd,P} - 1 + t_{even,P,ST} + v_{bdry,P}) +  2v_{odd,P} \\
        &= -v_{even,P} + v_{odd,P} - v_{bdry,P} + 1. & & \specialcell{\hfill\qedsymbol}
    \end{align*}
    \end{sloppypar}
    \begin{sloppypar}
    %In Fig. 8, we illustrate the simple observation that is at the heart of the proof of Lemma \ref{lemharnack}. The picture shows two primitive triangulations of the triangle with vertices $(0,0)$, $(4,0)$, $(0,4)$. Blue edges are adjacent to an even vertex, while red edges are not. Likewise, blue triangles are adjacent to an even vertex, while red triangles are not (and hence are even triangles). The even vertices appear in blue, the odd vertices appear in red. The triangulation on the right contains one more red triangle than the triangulation on the left, and hence is bound to contain one more red edge, one less blue edge and one less blue triangle than the triangulation on the left.
%    \begin{center}
%    \includegraphics[scale=1.25]{images/quantities.png}\\
%    \begin{minipage}{16cm}
%    \begin{it}
%    \label{quantities}
%    \small{Fig. 16: Two primitive triangulations of the triangle with vertices $(0,0)$, $(4,0)$, $(0,4)$.}
%    \end{it}
%    \end{minipage}
%\end{center}
\end{sloppypar}

\subsubsection{Relating $\chi_{L-AD,m}^+$ to the numbers $\chi_{T_k,L-AD}^+$ and $\chi_{IV,m}^+$ to the numbers $\chi_{T_k,IV}^+$}

\label{EulerSimplification}
\begin{sloppypar}
To compute $\chi_{L-AD,m}^+ - \chi_{IV,m}^+$, it suffices to count in both the $L$-AD and the even IV constructions the number of copies with only vertices of sign $+$ of simplices that lie in the star of the triangles $T_k$, when $8 \leq k \leq m$. Indeed, the rest of the data is identical for both constructions.
\end{sloppypar}
\begin{sloppypar}
Let $k$ be an even integer such that $8 \leq k \leq m$. We recall that the vertices of $T_k$ are the points $(k,0,0,m-k)$, $(0,k,0,m-k)$, $(0,0,k,m-k)$. The star of $T_k$ in the $L$-AD or IV triangulation of $\Delta_m^4$ is composed of the following parts:
\begin{itemize}
    \item all the parts of type $(00,10,11,12,22)$, for all triangulation edges of type $(00,22)$ with extremities lying either on the hyperplane $\{x_4=m-k+1 \}$ or, if $m \neq k$, on the hyperplane $\{ x_4= m-k-1 \}$;
    \item a cone $C_k$ over $T_k$ whose vertex lies on the hyperplane $\{ x_4=m-k \}$;
    \item a cone over $C_k$ with vertex $(0,0,0,m-k+1)$ and, if $m \neq k$, a cone over $C_k$ with vertex $(0,0,0,m-k-1)$.
\end{itemize}
\begin{defin}
    Let ST $\in \{L-$AD, IV $\}$. Let $P$ be a polytope whose vertices are integer points contained in $\Delta_m^4$. We call \textit{the contribution of $P$ to $\chi_{ST,m}^+$} the quantity $c_4-c_3+c_2-c_1+c_0$, where $c_i$ is the number of simplices of the triangulation of $\widetilde{(\Delta_m^4)^*}$ of dimension $i$ whose vertices all bear the sign $+$ and whose relative interior is contained in the relative interior of a symmetric copy of $P$ for each $i \in \{0,1,2,3,4\}$.
\end{defin}
We first simplify the computation by getting rid of unnecessary terms.
\begin{lem}
    Let $k$ be an even integer such that $8 \leq k < m$. Let $l$ be an integer such that $0 \leq l \leq 3$. Let $Q$ be a primitive $l$-simplex of the $L$-AD or IV triangulation lying in the hyperplane $\{ x_4=m-k \}$ and in exactly $l'$ facets of $\Delta_m^4$. The contribution to $\chi_{L-AD,m}^+$ or $\chi_{IV,m}^+$ of a suspension over $Q$ whose vertices lie respectively in $\{ x_4=m-k-1 \}$ and $\{ x_4=m-k+1 \}$, have the same parity, bear the same sign and lie in the same $3$-facets of $\Delta_m^4$ as $Q$ is zero.
\end{lem}
\textbf{Proof:} Let $Q_{-1}$ and $Q_{+1}$ be the cones over $Q$ with respective vertices $(0,0,0,m-k-1)$ and $(0,0,0,m-k+1)$. Recall from Lemma \ref{emptycopies} that $Q$ has $2^{4-l-l'}$ empty copies. The primitive $(l+1)$-simplices $Q_{-1}$ and $Q_{+1}$ have $2^{4-l-l'-1}$ empty copies.\end{sloppypar}
\begin{sloppypar}
If $Q$ has an even face, then so do $Q_{-1}$ and $Q_{+1}$. In this case, the vertices of the empty copies of $Q$, $Q_{-1}$ and $Q_{+1}$ all bear the same sign. If this sign is $-$, then the contribution of $Q$, $Q_{-1}$ and $Q_{+1}$ is zero. If it is $+$, then the contribution of $Q$, $Q_{-1}$ and $Q_{+1}$ is $\pm (2^{4-l-l'} - 2^{4-l-l'-1} -2^{4-l-l'-1})=0$.\end{sloppypar}
\begin{sloppypar}
If $Q$ does not have any even face, then neither do  $Q_{-1}$ and $Q_{+1}$. Hence, by Lemma \ref{invertsym} the contributions of $Q$, $Q_{-1}$ and $Q_{+1}$ is $\pm \frac{1}{2}(2^{4-l-l'} - 2^{4-l-l'-1} -2^{4-l-l'-1})=0$. \hfill\qedsymbol\end{sloppypar}
\begin{sloppypar}
 The argument above can be easily adapted to the case in which $T_k$ lies on the hyperplane $\{ x_4=0 \}$. Hence, the contributions to $\chi_{L-AD,m}^+$ and $\chi_{IV,m}^+$ of  $C_k \setminus T_k$, of the suspension over $C_k \setminus T_k$, of $T_k$ and of the cone over $T_k$ lying in the hyperplane $\{x_1+x_2+x_3+x_4=m\}$ are zero.
%\begin{itemize}
    %\item $C_k \setminus T_k$ and the suspension over $C_k \setminus T_k$;
    %\item $T_k$ and the cones over $T_k$ lying in the hyperplane $\{x+y+z+t =m\}$.
%\end{itemize}
 The remaining parts of the star of $T_k$ are joins of $T_k$ with segments and with their endpoints.
\begin{lem}
    The contribution to $\chi_{L-AD,m}^+$ of the join of the boundary of $T_k$ with a segment or a point is the same as its contribution to $\chi_{IV,m}^+$.
\end{lem}
\textbf{Proof:} There are only three types of simplices in the boundary of $T_k$: even vertices, odd vertices and edges adjacent to an even vertex. The number of simplices of each of these three types in the $L$-AD construction is the same as in the even IV construction. An even vertex lying in the boundary of $T_k$ is adjacent to two edges lying in the boundary of $T_k$. If an even vertex has the sign $-$, then the contribution of all the simplices in the interior of its star is zero. This is the case in the even IV construction. In the $L$-AD construction, some of the even vertices bear the sign $+$. The vertices of the triangle $T_k$ still bear the sign $-$. Consider an even vertex bearing the sign $+$ and the two adjacent edges lying in the boundary of $T_k$. They lie in the intersection of a coordinate hyperplane with the hyperplane $\{x_1+x_2+x_3+x_4=m-k\}$. Their contribution to $\chi_{L-AD,m}^+$ is equal to $4-2\times 2 = 0$. Now the join of the vertex and its adjacent segments with any $l$-simplex $Q$ lying on the hyperplane $\{ x_4=m-k+1 \}$ or $\{ x_4=m-k-1\}$ consists of one $l$-simplex and two $(l+1)$-simplices. All of these simplices are adjacent to an even vertex. Hence their contribution to $\chi_{L-AD,m}^+$ is a multiple of $\pm (2^{4-l}-2\times 2^{4-l-1})=0$. \hfill\qedsymbol
\end{sloppypar}
\begin{sloppypar}
We can now study the remaining parts of the star of $T_k$.
\begin{lem}
    \label{contributiontriangle}
    Let ST $\in \{L-$AD, IV $\}$. Let $k$ be an even integer such that $8\leq k <m$. Then the contribution to $\chi_{ST,m}^+$ of the star of $T_k$ is equal to $-((k+1)\chi_{T_k,ST}^+ + (k-1)\chi_{T_k, ST}^+)$.
\end{lem}
\textbf{Proof:} In the star of $T_k$, it remains to investigate:
\begin{itemize}
    \item the join of the interior of $T_k$ with the $k+1$ segments of type $(00,22)$ lying on the hyperplane $\{x_4=m-k-1\}$ and with the $k-1$ segments of type $(00,22)$ lying on the hyperplane $\{x_4=m-k+1\}$;
    \item the join of the interior of $T_k$ with the $k+1$ vertices of type $(00)$ or $(22)$ that lie on the hyperplane $\{x_4=m-k-1\}$ but not on the hyperplane $\{ x_1+x_2+x_3+x_4= m \}$ and with the $k-1$ vertices of type $(00)$ or $(22)$ that lie on the hyperplane $\{x_4=m-k+1\}$ but not on the hyperplane $\{ x_1+x_2+x_3+x_4= m \}$.
    %\item the join of the interior of $T_k$ with the $k-1$ segments of type $(00,22)$ lying on the hyperplane $\{x_4=m-k+1\}$;
    %\item the join of the interior of $T_k$ with the $k-1$ vertices of type $(00)$ or $(22)$ that lie on the hyperplane $\{x_4=m-k+1\}$ but not on the hyperplane $\{ x+y+z+x_4= m \}$.
\end{itemize}
First, we recall that every triangle of the triangulation of $T_k$ has an even face. Hence, taking the join with a segment on the hyperplane $\{ x_4=m-k\pm 1 \}$ cannot produce any even triangle or $4$-simplex, as each $4$-simplex has exactly one even face. In other words, the join $Q \ast Q'$ of a simplex $Q$ of the triangulation of $T_k$ with a simplex $Q'$ of dimension $0$ or $1$ lying on the hyperplane $\{ x_4=m-k\pm 1 \}$ has an even $l$-face if and only if $Q$ has an even $l$-face. The result then follows from Lemmas \ref{invertsym} and \ref{emptycopies}: the contribution of the join of the interior of $T_k$ with a vertex is equal to $-2\chi_{T_k, ST}^+$ and the contribution of the join of the interior of $T_k$ with a segment is equal to $\chi_{T_k, ST}^+$. These arguments can easily be adapted to the case $k=m$. 
\hfill\qedsymbol
\end{sloppypar}

\subsection{The family of the AD hypersurfaces is not asymptotically of type $\chi = \sigma$}

\begin{sloppypar}
    Let $m$ be an even positive integer such that $m \geq 8$. Consider the list $\Lambda_m=\{ (a_k,b_k,c_k) \}_{8 \leq k \leq m, \: k \: even}$ where $a_k = (k-1,0,1,m-k)$, $b_k=(1,0,k-1,m-k)$, $c_k = (2,k-5,3,m-k)$. We denote by $X_{AD,m}$ the $\Lambda_m$-AD hypersurface and we call it \textit{the AD hypersurface of degree} $m$. The AD hypersurface of degree $m$ is maximal, as stated in Proposition \ref{maxAD}. It remains to prove the following proposition.
    \begin{prop}
    \label{propAsymptoDev}
        The family of maximal $3$-manifolds $(X_{AD,2k})_{k \geq 4}$ satisfies the relation $\chi(\mathbb{R}Y_{AD,2k}^-) \sim_{+\infty} \sigma(\mathbb{C}Y_{AD,2k}) - \frac{(2k)^4}{4}$, where $Y_{AD,2k}$ is a double covering of $\mathbb{CP}^4$ branched along the complex part of $X_{AD,2k}$.
    \end{prop}
    \begin{rmk}
        There are other families of lists of triplets of points that give rise to families of maximal hypersurfaces that are not asymptotically of type $\chi = \sigma$. We only present the family that produces the largest asymptotical deviation.
    \end{rmk}
\end{sloppypar}

\subsubsection{Preliminary Euler characteristic computations}

\begin{sloppypar}
     Retaining the notations from Section \ref{eulercharcompLAD}, we apply Lemma \ref{lemharnack} in order to compute $\chi_{T_k^1,AD}^+$, $\chi_{T_k^2,AD}^+$ and $\chi_{T_k \setminus (T_k^1 \bigcup T_k^2),AD}^+$.
    \begin{lem}
        For fixed positive even integers $m$ and $k$ such that $8 \leq k \leq m$, one has $\chi_{T_k^1,AD}^+ = -\frac{k}{2}$.
        %\begin{center}$\chi_{T_k^1,AD}^+ = -\frac{k}{2}$.\end{center}
    \end{lem}
    \textbf{Proof:} The interior of the $T_k^1$ does not contain any integer point. The even points lying on the segment $[(0,k,0,m-k), (0,0,k,m-k)]$ bear the sign $-$.\\ 
    We obtain $\chi_{T_k^1,AD}^+ = -v_{bdry,T_k^1} + 1 = -(\frac{k}{2}+1)+1$. \hfill\qedsymbol
    \begin{lem}
        Let $m$ and $k$ be even integers such that $8 \leq k \leq m$. Then, one has
        \begin{center} $\chi_{T_k^2,AD}^+ = -\frac{k^2}{4} + k + 12 + 4 \lfloor \frac{k-6}{4} \rfloor$. \end{center}
    \end{lem}
    \textbf{Proof}: The lattice volume of $T_k^2$ is equal to $\frac{k^2}{2}-\frac{7k}{2}+5$ and the interior of $T_k^2$ contains $\frac{k^2}{2}-4k+6$ integer points.
     All the even vertices in the interior of $T_k^2$ and on its boundary have the sign $+$.\\
     We have $v_{even,T_k^2} = \frac{k^2}{8}-\frac{10k}{8}+3 + \lfloor \frac{k-6}{4} \rfloor$ and $v_{bdry,T_k^2} = \frac{k}{2}-1$, as well as $v_{odd,T_k^2} = \frac{3k^2}{8} - \frac{22k}{8} + 3 -\lfloor \frac{k-6}{4} \rfloor$.\\
    %\begin{itemize}
    %    \item $v_{even,T_k^2} = \frac{k^2}{8}-\frac{10k}{8}+3 + \lfloor \frac{k-6}{4} \rfloor$;
    %    \item $v_{bdry,T_k^2} = \frac{k}{2}-1$;
    %    \item $v_{odd,T_k^2} = \frac{3k^2}{8} - \frac{22k}{8} + 3 -\lfloor \frac{k-6}{4} \rfloor$.
    %\end{itemize}
    Therefore, $\chi_{T_k^2,AD}^+ = -2\mathcal{A}(T_k^2) + 3 v_{even,T_k^2} + v_{odd,T_k^2} + v_{bdry,T_k^2} + 1 = -\frac{k^2}{4} + k + 12 + 4 \lfloor \frac{k-6}{4} \rfloor$.\hfill\qedsymbol
    \begin{lem}
         For fixed positive even integers $m$ and $k$ such that $8 \leq k \leq m$, one has
         \begin{center}$\chi_{T_k \setminus (T_k^1 \bigcup T_k^2),AD}^+ = k-2+2\lfloor \frac{k-6}{4}\rfloor$.\end{center}
    \end{lem}
    \textbf{Proof:} A primitive triangulation of $T_k \setminus (T_k^1 \bigcup T_k^2)$ can contain even triangles. The even vertices lying in $T_k \setminus (T_k^1 \bigcup T_k^2)$ bear the sign $-$. We have \mbox{$v_{even, T_k \setminus (T_k^1 \bigcup T_k^2)} = \frac{k}{2} - 2 - \lfloor \frac{k-6}{4} \rfloor$} and \mbox{$v_{odd, T_k \setminus (T_k^1 \bigcup T_k^2)} = 2k-4+ \lfloor \frac{k-6}{4} \rfloor$}, as well as \mbox{$v_{bdry, T_k \setminus (T_k^1 \bigcup T_k^2)} = \frac{k}{2} + 1$}.
    %We have:
    %\begin{itemize}
    %    \item $v_{even, T_k \setminus (T_k^1 \bigcup T_k^2)} = \frac{k}{2} - 2 - \lfloor \frac{k-6}{4} \rfloor$;
    %    \item $v_{odd, T_k \setminus (T_k^1 \bigcup T_k^2)} = 2k-4+ \lfloor \frac{k-6}{4} \rfloor$;
    %    \item $v_{bdry, T_k \setminus (T_k^1 \bigcup T_k^2)} = \frac{k}{2} + 1$.
    %\end{itemize}
    Thus, \mbox{$\chi_{T_k \setminus (T_k^1 \bigcup T_k^2),AD}^+ = k-2+2\lfloor \frac{k-6}{4}\rfloor$}. 
    \begin{flushright}
        \qedsymbol
    \end{flushright} 
    \begin{prop}
        \label{contributionAD}
        Let $m$ and $k$ be even integers such that $8 \leq k \leq m$. Then, one has
        \begin{center}
            $\chi_{T_k,AD}^+ = -\frac{k^2}{4}+\frac{3k}{2}+10+6 \lfloor \frac{k-6}{4} \rfloor$.
        \end{center}
    \end{prop}
    \begin{align*}
    \textbf{Proof:} \; \chi_{T_k,AD}^+ &= \chi_{T_k^1,AD}^+ + \chi_{T_k^2,AD}^+ + \chi_{T_k \setminus (T_k^1 \bigcup T_k^2),AD}^+ + \chi_{[(1,0,k-1,m-k), (2, k-5,3,m-k)],AD}^+ \\
    &\;\;\;\;+ \chi_{[(k-1,0,1,m-k),(2, k-5,3,m-k)],AD}^+ + \chi_{[(1,0,k-1,m-k),(0,k,0,m-k)], AD}^+ + \chi_{(2, k-5,3,m-k),AD}^+\\
    &= \chi_{T_k^1,AD}^+ + \chi_{T_k^2,AD}^+ + \chi_{T_k \setminus (T_k^1 \bigcup T_k^2),AD}^+ + (-1) + (-1) + 0 + 2\\
    %& = \chi_{T_k^1,AD}^+ + \chi_{T_k^2,AD}^+ + \chi_{T_k \setminus (T_k^1 \bigcup T_k^2),AD}^+\\
    & = -\frac{k^2}{4}+\frac{3k}{2}+10+6 \lfloor \frac{k-6}{4} \rfloor. &\specialcell{\hfill\qedsymbol}
    \end{align*}
    
\end{sloppypar}

\newpage

\subsubsection{Proof of Theorem \ref{theorem2}}

\begin{sloppypar}
    Now, we use the computations from the previous section to prove Proposition \ref{propAsymptoDev} and hence Theorem \ref{theorem2}.
    \end{sloppypar}
    \begin{sloppypar}
    \textbf{Proof:}
    Let $m$ be an even integer such that $m \geq 8$. Putting together the results of Lemmas \ref{contributiontriangle}, \ref{contributionIV} and Proposition \ref{contributionAD}, we obtain:
    \begin{align*}
        \chi_{AD,\: m}^+ - \chi_{IV,\: m}^+ &= -(\chi_{T_m,AD}^+ - \chi_{T_m,IV}^+) (m-1) + \sum_{k=4}^{\frac{m}{2}-1} -(\chi_{T_{2k},AD}^+ - \chi_{T_{2k},IV}^+) (2k-1 + 2k+1)\\
         &= -(-\frac{m^2}{2} + 3m + 10 + 6 \lfloor \frac{(m-6)}{4} \rfloor)(m-1) \\
         &\:\:\:\:\:\: +   \sum_{k=4}^{\frac{m}{2}-1} -(-\frac{(2k)^2}{2} + 3(2k) + 10 + 6 \lfloor \frac{(2k-6)}{4} \rfloor ) 4k\\
         &\sim_{+\infty} \frac{m^4}{8}.
    \end{align*}
    Let $Y_{AD,m}$ (respectively, $Y_{IV,m}$) be a double covering of $\mathbb{CP}^4$ branched along the AD hypersurface of degree $m$ (respectively, along the even IV hypersurface of degree $m$). Using the notations from Section \ref{definitions}, we have
    \begin{align*}
        \chi(\mathbb{R}Y_{AD,m}^-) - \sigma(\mathbb{C}Y_{AD,m}) = \chi(\mathbb{R}Y_{AD,m}^-) - \chi(\mathbb{R}Y_{IV,m}^-)
        %&= \chi_{AD,\: m}^- - \chi_{IV,\: m}^-
        = 2(1 - \chi_{AD,\: m}^+) - 2(1 -  \chi_{IV,\: m}^+)
         \sim_{+\infty} -\frac{m^4}{4}.
    \end{align*}
    Together with Proposition \ref{maxAD}, this finishes the proof of Theorem \ref{theorem2}. \hfill \qedsymbol
    %Hence, the family of maximal AD hypersurfaces obtained is not asymptotically of type $\chi = \sigma$. 
\end{sloppypar}

%\input{figures.tex}

%BIBLIO

\newpage
%\bibliographystyle{alpha}
%\bibliography{biblio}

\end{document}